\documentclass[12pt]{article}

\usepackage[letterpaper,margin=1in,footskip=0.25in]{geometry}
\usepackage{amssymb}
\usepackage{amsthm}
\usepackage{amsfonts}
\usepackage{mathtools}
\usepackage{amsmath}
\usepackage[makeroom]{cancel}
\usepackage{mathrsfs}
\usepackage{tikz}
\usetikzlibrary{arrows}
\usepackage{stmaryrd}
\usepackage[bottom]{footmisc}
\usepackage{enumitem}
\usepackage{todonotes}
\usepackage[hidelinks]{hyperref}
\usepackage[mathscr]{eucal} 
\usepackage[title]{appendix}
\usepackage[nottoc]{tocbibind}
\usepackage[all,cmtip]{xy}
\usepackage{float} 
\usepackage{lipsum}
\usepackage{setspace}

\theoremstyle{definition}
\newtheorem{definition}{Definition}[section]
\newtheorem{example}[definition]{Example}
\newtheorem{notation}[definition]{Notation}

\theoremstyle{remark}
\newtheorem{conjecture}[definition]{Conjecture}
\newtheorem{remark}[definition]{Remark}

\theoremstyle{plain}
\newtheorem{theorem}[definition]{Theorem}
\newtheorem{lemma}[definition]{Lemma}
\newtheorem{proposition}[definition]{Proposition}
\newtheorem{corollary}[definition]{Corollary}

\newcommand\abs[1]{\lvert #1 \rvert}

\makeatletter
\newcommand*\bigcdot{ {\mathpalette\bigcdot@{.5}} }
\newcommand*\bigcdot@[2]{\mathbin{\vcenter{\hbox{\scalebox{#2}{$\m@th#1\bullet$}}}}}
\makeatother

\newcommand\NN{ \mathbb{N} }

\newcommand\RR{ \mathbb{R} }

\newcommand\ZZ{ \mathbb{Z} }

\newcommand\cV{ \mathcal{V} }

\newcommand\sC{ \mathscr{C} }
\newcommand\sD{ \mathscr{D} }

\newcommand\sN{ \mathscr{N} }

\newcommand\sS{ \mathscr{S} }

\newcommand\Int{ {\operatorname{Int} } }
\newcommand\Lie[1]{ { \mathcal{L} }_{#1} }

\newcommand\pre{ {\operatorname{pre}} }
\newcommand\Rp{ { \mathbb{R}_{>0} } }

\newcommand\str{\operatorname{star}}

\newcommand\cat{\operatorname{cat}}
\newcommand\Cat{ {\operatorname{Cat}} }

\newcommand\Fun{ {\operatorname{Fun}} }
\newcommand\id{\operatorname{id}}

\newcommand\Ind{ {\operatorname{Ind} } }
\newcommand\Hom{\operatorname{Hom}}

\newcommand\Mor{ {\operatorname{Mor}} }

\newcommand\PrLcg[1][\omega]{  {  \operatorname{Pr}^{\operatorname{L} }_#1 } }
\newcommand\PrLcs{  {  \operatorname{Pr}^{\operatorname{L} }_{\omega,st} } }
\newcommand\PrLst{  {  \operatorname{Pr}^{\operatorname{L} }_{st} } }

\newcommand\PrRst{  {  \operatorname{Pr}^{\operatorname{R} }_{st} } }

\newcommand\colim{\operatorname{colim}}
\newcommand\clmi[1]{  \underset{ {#1} }{\operatorname{colim}} }

\newcommand\lmi[1]{  \underset{ {#1} }{\lim} }

\newcommand\cof{\operatorname{cof}}

\newcommand\dMod{\operatorname{-Mod}}
\newcommand\fib{\operatorname{fib}}

\newcommand\st{   {\operatorname{st}} }

\newcommand\Loc{ {\operatorname{Loc}} }

\newcommand\msh{\operatorname{\mu sh}}
\newcommand\ms{\operatorname{SS}}
\newcommand\msif{\operatorname{SS}^{\infty}}
\newcommand\msnz{ \dot{\operatorname{SS} }}
\newcommand\Op{\operatorname{Op}}
\newcommand\PSh{\operatorname{PSh}}
\newcommand\Sh{\operatorname{Sh}}

\newcommand\sHom{\operatorname{\mathscr{H}om}}

\newcommand\supp{ {\operatorname{supp}} }
\newcommand\wsh{\operatorname{\mathfrak{w} sh}}

\newcommand\Fuk{\operatorname{Fuk}}
\newcommand\HF{ {\operatorname{HF}} }

\newcommand\WF{ {\mathcal{W} } }
\newcommand\wrap{\mathfrak{W}}

\newcommand\Perf{\operatorname{Perf}}

\title{\textbf{Wrapped sheaves}}
\date{\today}
\author{Christopher Kuo}

\newcommand\constr{ {\operatorname{constr}  } }
\newcommand\PrLV[1][\cV]{  {  \operatorname{Pr}^{\operatorname{L} }_{\cV,st} } }

\newcommand\dT{ {\dot{T} } }

\begin{document}

\maketitle

\begin{abstract}

We construct a sheaf-theoretic analogue of the wrapped Fukaya category in Lagrangian Floer theory, by localizing 
a category of sheaves microsupported away from some given $\Lambda \subset S^*M$ along 
continuation maps constructed using the Guillermou-Kashiwara-Schapira sheaf quantization. 
%
%

When $\Lambda$ is a subanalytic singular isotropic, we also construct 
a comparison map 
to the category of compact objects in the category of unbounded sheaves microsupported in $\Lambda$, 
and show that it is an equivalence.
The last statement can be seen as a sheaf-theoretical incarnation of the sheaf-Fukaya comparison theorem of 
Ganatra-Pardon-Shende.
\end{abstract}

\tableofcontents

\section{Introduction}
The microlocal sheaf theory of Kashiwara and Schapira \cite{Kashiwara-Schapira1} 
relates sheaves -- which are topological structures
-- on manifolds with the symplectic geometry of their cotangent bundles.  The basic construction is the {\em microsupport}:
to a sheaf $F$ on a $C^1$ manifold $M$, one associates a  closed conic subset $\ms(F)$ in the cotangent bundle $T^* M$.
The intersection of this set with the zero section recovers the support $\supp(F)$; 
the projectivization $\msif(F) \coloneqq (\ms(F) \setminus 0_M) / \Rp$ indicates the codirections along which the sheaf changes. 
A key indicator of the symplectic nature of the theory is the involutivity theorem \cite[Thm. 6.5.4]{Kashiwara-Schapira1}, which
asserts that $\ms(F)$ is always a singular coisotropic subset with respect to the 
canonical symplectic structure of $T^* M$.  Under an appropriate tameness assumption,
$\ms(F)$ is a singular Lagrangian if and only if $F$ is constructible, i.e.,
there exists a stratification $\{ X_s \}$ of $M$ such that $F|_{X_s}$ is a local system for all $X_s$.

Deeper relationships between microlocal sheaf theory and symplectic geometry began to emerge in the mid 2000s.
Nadler and Zaslow related constructible sheaves to `infinitesimally wrapped' Floer theory \cite{Nadler-Zaslow, Nadler-brane}. 
Meanwhile, Tamarkin introduced purely sheaf-theoretical methods into symplectic topology in his study of
non-displaceability, a problem previously studied largely by Floer theoretic methods \cite{Tamarkin1}. 
The subsequent 
Guillermou-Kashiwara-Schapira sheaf quantization of contact isotopies \cite{Guillermou-Kashiwara-Schapira} --- 
i.e. the highly nonobvious statement that contact isotopies of $S^*M$ act on sheaves on $M$ --- led to a 
host of further incursions by sheaf theorists into symplectic topology 
\cite{Guillermou1, Guillermou2, Guillermou3, Shende-Treumann-Zaslow, Shende-Treumann-Williams-Zaslow, Chiu1,
Chiu2, Shende-conormal, Treumann-Zaslow, Casals-Gao, Casals-Zaslow, Ike, Asano-Ike1, Asano-Ike2, Bai-Cote, Li1, Li2}
 and vice versa \cite{Viterbo, Zhang, Rutherford-Sullivan}.

Meanwhile on the Floer theoretic side, Abouzaid and Seidel formulated a way to incorporate contact dynamics into 
Fukaya categories for noncompact symplectic manifolds \cite{Abouzaid-Seidel}.  Their construction is roughly
to localize a partially-defined infinitesimally wrapped Fukaya category along `continuation morphisms' associated to 
positive isotopies.  The resulting notion of wrapped
Fukaya category (and its later `partially wrapped' generalizations \cite{Sylvan1, Ganatra-Pardon-Shende1, Ganatra-Pardon-Shende2}) 
provides the correct mirrors to coherent sheaf categories on certain singular spaces.  
While such wrapped categories are nontrivial to compute directly (the simplest case was \cite{Abouzaid-Auroux-Efimov-Katzarkov-Orlov}), 
Nadler conjectured that they matched categories of compact objects inside categories of unbounded sheaves
with prescribed microsupport \cite{Nadler-pants}.  This conjecture was later established by the work of Ganatra, Pardon, and Shende \cite{Ganatra-Pardon-Shende3}; 
as a result, sheaf-theoretic methods 
(e.g.  \cite{Nadler-pants, Nadler-LG, Fang-Liu-Treumann-Zaslow, Kuwagaki1, Gammage-Shende1, Gammage-Shende2}) 
can be used to establish homological mirror symmetry in these settings.

While Nadler termed his category the `wrapped sheaves', the name did not entirely reflect its construction:
there is no wrapping in their definition.  In this paper we study the category which this name manifestly suits --- the localization
of the category of sheaves along the continuation morphisms of \cite{Guillermou-Kashiwara-Schapira} --- and give an entirely 
sheaf-theoretic proof that it is 
equivalent to Nadler's category.   

\begin{figure}[h]
  \centering
    \includegraphics[width=1\textwidth]{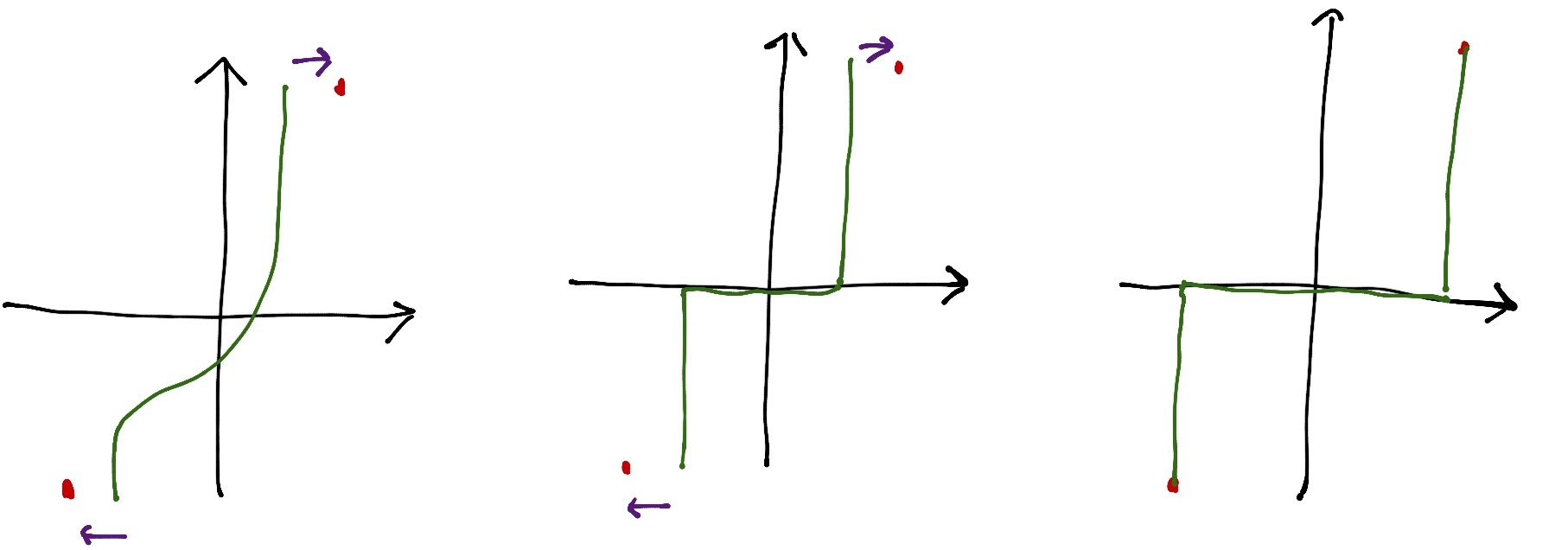}
  \caption{\label{alternative} An alternative to \cite{Ganatra-Pardon-Shende3}, illustrated for $M = \mathbb{R}^1$. 
  The purple arrow indicates the positive Reeb direction and the red dots are the stop $\Lambda$. 
  The green objects, from left to right, are a smooth curve conic at infinity 
  in $T^* \mathbb{R}^1$, the microsupport of a constant sheaf on an open interval,
  and the microsupport of the same sheaf after being pushed to the stop.}
\end{figure}

\begin{remark}
Let us compare and contrast this article with the work \cite{Ganatra-Pardon-Shende3}.   In that article, the authors construct
an equivalence between the partially wrapped Fukaya category of a cotangent bundle, stopped along
some subanalytic isotropic $\Lambda$, with the category of compact objects in the category of unbounded
sheaves microsupported in $\Lambda$: 
$$\Perf \WF(T^* M,\Lambda)^{op} \cong \Sh_\Lambda(M)^c$$
Their approach was to introduce an abstract axiomatic characterization (`microlocal Morse theatre')
and verify that both sides satisfy it. 

Here is another approach which is illustrated in the following diagram:
$$
\begin{tikzpicture}
\node at (-3.2,2) {$\Fuk_\epsilon(T^* M)^{op}$};
\node at (-1.7,2) {$\supseteq$};
\node at (0,2) {$\Fuk_\epsilon(T^* M,\Lambda)^{op}$};

\node at (6,2) {$\WF(T^* M,\Lambda)^{op}$};

\node at (-3.2,0) {$\Sh_\constr(M)$};
\node at (-1.9,0) {$\supseteq$};
\node at (0,0) {$\Sh_{\constr,S^*M \setminus \Lambda}(M)$};

\node at (6,0) {$\wsh_\Lambda(M)$};

\node at (9.8,0) {$\Sh_\Lambda(M)^c$};

\draw [->, thick, dashed] (1.6,2) -- (4.7,2) node [midway, above] {$ $};
\draw [->, thick] (1.7,0) -- (5,0) node [midway, below] {$ $};
\draw [->, thick] (7,0) -- (8.8,0) node [midway, left] {$ $}; 

\draw [double equal sign distance, thick] (-3,0.4) -- (-3,1.6) node [midway, below] {$ $};
\draw [double equal sign distance, thick] (0,0.4) -- (0,1.6) node [midway, below] {$ $};
\draw [<->, thick, dashed] (6,0.4) -- (6,1.6) node [midway, below] {$ $};

\draw [double equal sign distance, thick] (7.2,1.8) -- (9,0.3) node [midway, below] {$ $};

\node at (-2.2,1) {$\cite{Nadler-Zaslow, Nadler-brane}$};
\node at (0.8,1) {$\cite{Nadler-Zaslow, Nadler-brane}$};
\node at (8.6,1.2) {$\cite{Ganatra-Pardon-Shende3}$};

\node at (7.9,0.3) {$\wrap_\Lambda^+$};
\node at (3.3,0.3) {$l_{\Sh}$};

\node at (3.3,2.3) {$l_{\Fuk}$};

\end{tikzpicture}
$$ 

Assume given an equivalence of infinitesimally wrapped categories such as asserted in \cite{Nadler-Zaslow, Nadler-brane},  
localize both sides along the continuation morphisms,
and then show purely on the sheaf-theoretic side that the resulting localized category of sheaves is in fact equivalent to $\Sh_\Lambda(M)^c$.
See Figure \ref{alternative} for an illustration.
The present article establishes the last step, indicated by $\wrap_\Lambda^+$, in this argument.  
We remark that the dotted arrows, as far as the author is aware of, are still conjectural.
While this route of proof would be logically independent of \cite{Ganatra-Pardon-Shende1, Ganatra-Pardon-Shende2, Ganatra-Pardon-Shende3},
the proof of the main theorem of the present article follows a strategy adapted from \cite{Ganatra-Pardon-Shende3}.  

Note however that, in the literature, the category $\WF(T^* M,\Lambda)$ is not constructed as a localization of an infinitesimally wrapped category.
In fact, constructions of $\WF(T^* M, \Lambda)$ as in \cite{Abouzaid-Seidel, Sylvan1, Ganatra-Pardon-Shende1},
always avoid infinitesimal wrappings.
For example, the authors in \cite{Ganatra-Pardon-Shende1} begin with the sub-semi-category $\Fuk^\pre(T^* M,\Lambda)$ 
with morphisms defined only between transversal Lagrangians, and roughly speaking localize along continuation maps 
so they only ever have to work with large wrappings.

\end{remark}

\vspace{2mm}
We turn to a more precise discussion of our results.  Let us first discuss the continuation maps in the sheaf-theoretical setting.
Recall that for a closed conic subset $X \subseteq T^\ast M$, the inclusion $\iota_\ast: \Sh_X(M) \hookrightarrow \Sh(M)$ of
sheaves microsupported in $X$ to all sheaves has a left adjoint $\iota^\ast$ and right adjoint $\iota^!$.
Now consider $F \in \Sh(M \times \mathbb{R})$ as a one-parameter family of sheaves on $M$ and set $F_a \coloneqq F|_{ M \times \{a\}}$.
An explicit description of $\iota^\ast$ given in \cite[Proposition 4.8]{Guillermou-Kashiwara-Schapira} shows
when $F$ satisfies the condition $\ms(F) \subseteq T^\ast M \times T^\ast_{\leq 0} \mathbb{R}$, there is a 
\textit{continuation map} $F_a \rightarrow F_b$ for $a \leq b$.
(See \cite[2.2.2]{Tamarkin1} and \cite[(77)]{Guillermou-Schapira} for the dual construction.)
Now pick a contact form $\alpha$ on $S^\ast M$ coorienting the contact structure induced from the symplectic structure on $T^\ast M$.
We say a $C^\infty$ map $\Phi: S^\ast M \times \mathbb{R} \rightarrow S^\ast M$ is an isotopy if the induced map 
$\phi_t \coloneqq \Phi(-,t)$ is a contactomorphism for all $t \in \mathbb{R}$ and $\phi_0 = \id_{S^\ast M}$.
If $\Phi$ is a positive isotopy ($\alpha(\partial_t \phi_t) \geq 0$)
then the corresponding GKS sheaf kernel $K(\Phi)$ and hence its convolution $K(\Phi) \circ F$ with $F \in \Sh(M)$ 
will satisfy this condition and hence admit continuation maps.

Now fix an open set $\Omega \subseteq S^\ast M$.
Homotopy classes of compactly supported isotopies with fixed ends can be organized to an $\infty$-category $W(\Omega)$ 
whose morphisms are given by concatenating with positive isotopies.
We refer this category as the category of \textit{positive wrappings}.
The discussion on continuation maps will imply that there is a \textit{wrapping kernel} functor $w: W(\Omega) \rightarrow \Sh(M \times M)$
which sends isotopies to the end point of the GKS sheaf kernels and positive isotopies to continuation maps.
One can use this functor to the define the \textit{infinite wrapping} 
functors $\wrap^\pm(\Omega): \Sh(M) \rightarrow \Sh_{S^\ast M \setminus \Omega}(M)$
by sending $F$ to the colimit $\colim_\Phi \left( w(\Phi) \circ F \right)$ or limit $\lim_\Phi \left( w(\Phi) \circ F \right)$ over $\Phi \in W(\Omega)$.
Geometrically, we push $F$ with increasingly positive (resp. negative) isotopies and take colimit (resp. limit) over them.
These functors give a geometric description for the adjoints of the inclusion 
$\iota_\ast: \Sh_{S^\ast M \setminus \Omega}(M) \hookrightarrow \Sh(M)$.

\begin{theorem}\label{w=ad}
Let $\iota_\ast: \Sh_{S^\ast M \setminus \Omega}(M) \hookrightarrow \Sh(M)$ denote the tautological inclusion.
Then the functor $\wrap^+(\Omega)$ (resp. $\wrap^-(\Omega)$) is the left (resp. right) adjoint of $\iota_\ast$.
\end{theorem}
See subsection \ref{Sheaf-theoretical wrappings} for the \hyperlink{proof_w=ad}{proof}.

The main construction of this paper is the category of \textit{wrapped sheaves} 
$\wsh_\Lambda(M)$ where
$M$ is a real analytic manifold and $\Lambda$ is a closed subset in $S^\ast M$.
It is a stable category defined by first collecting sheaves which have subanalytic singular isotropic microsupport away from $\Lambda$ 
and those which are compactly isotopic to them in $S^\ast M \setminus \Lambda$, 
and then inverting continuation maps which come from positive isotopies satisfying similar conditions.
One effect of this localization is that objects which can be connected through an isotopy on $S^\ast M \setminus \Lambda$
will be identified.
We show that $\Hom$'s in $\wsh_\Lambda(M)$ can be computed as colimits of 
$\Hom$'s between ordinary sheaves over $W(S^\ast M \setminus \Lambda)$.
Finally, when $\Lambda$ is a subanalytic singular isotropic, by using the infinite wrapping functor $\wrap^+(S^\ast M \setminus \Lambda)$, 
we define a canonical \textit{comparison} functor $\wrap_\Lambda^+(M): \wsh_\Lambda(M) \rightarrow \Sh_\Lambda(M)^c$. 
The main theorem of this paper is that $\wrap_\Lambda^+(M)$ is an equivalence.

\begin{theorem}\label{main}
Let $\Lambda \subseteq S^\ast M$ be a subanalytic singular isotropic.
The comparison functor $ \wrap_\Lambda^+(M): \wsh_\Lambda(M) \rightarrow \Sh_\Lambda(M)^c$ is an equivalence.
\end{theorem}
See subsection \ref{proof_of_the_main_theorem} for the \hyperlink{proof_main}{proof}.

\begin{remark}
Note that, unlike the analogous isomorphism in \cite{Ganatra-Pardon-Shende3}, our isomorphism is induced by an explicit functor. 
\end{remark}

Since all the above constructions are functorial on the inclusion of open sets of $M$, we obtain a precosheaf $\wsh_\Lambda$
and we refer its objects as the \textit{wrapped sheaves}.
The corollary of the above theorem is that this precosheaf is a cosheaf.
\begin{corollary}\label{maincor}
Let $\Lambda \subseteq S^\ast M$ be a subanalytic singular isotropic.
The comparison morphism $\wrap_\Lambda^+: \wsh_\Lambda \rightarrow \Sh_\Lambda^c$ between precosheaves is an isomorphism.
In particular, the precosheaf $\wsh_\Lambda$ is a cosheaf.
\end{corollary}

The proof of Theorem \ref{main} follows the same strategy as \cite{Ganatra-Pardon-Shende3}.
In short, subanalytic geometry implies that, for a subanalytic singular isotropic $\Lambda$, 
there exists a $C^1$ Whitney triangulation $\mathscr{S}$
such that $\Lambda$ is contained in $N^\ast_\infty \mathscr{S} \coloneqq \cup_{s \in \mathscr{S}} \, N^\ast_\infty X_s$.
For this special case, the two categories are natural identified as $\Perf \mathscr{S}$, 
the category of perfect $\mathscr{S}$-modules,
and hence admit a preferred set of generators which are matched under $\wrap_{N^\ast_\infty \mathscr{S}}^+(M)$.
We then apply the nearby cycle technology developed in \cite{Nadler-Shende} to conclude that $\wrap_{N^\ast_\infty \mathscr{S}}^+(M)$
induces an equivalence on the $\Hom$'s for these generators and hence finished the proof for this case.
To conclude the theorem for the general case, we note that the construction on both sides are contravariant on $\Lambda$.
Thus we study the fiber of the canonical maps $\wsh_{N^\ast_\infty \mathscr{S}} (M) \rightarrow \wsh_\Lambda(M)$ and 
$\Sh_{N^\ast_\infty \mathscr{S}}(M)^c \rightarrow \Sh_\Lambda(M)^c$, and show that they are generated by a sheaf-theoretical
version of the linking disks and microstalks at the smooth points of $N^\ast_\infty \mathscr{S} \setminus \Lambda$.
Finally, we show that $\wrap_\Lambda^+ (M)$ matches those objects and thus conclude the general case.


\subsection{Acknowledgement}
The author would like to thank his advisor Vivek Shende for pointing out the main ideas of the constructions performed in this paper as well as helpful guidance,
Germ\'an Stefanich for helpful discussions concerning higher category theory, 
and Pierre Schapira for carefully reading the first version of this paper on arXiv. 
The author would also like to thank Laurent C\^{o}t\'{e}, Wenyuan Li, and Harold Williams for useful discussions during the revision of the paper.
Finally, the author would like to thank the referee for their careful reading of this paper.
This project was supported by NSF CAREER DMS-1654545.

\newpage

\section{Microlocal sheaf theory}
We will consider microlocal sheaf theory in the higher categorical setting over a fixed coefficient category $\cV$.
We begin this section with recalling standard higher categorical notations from \cite{Lurie1, Lurie2},
and classical definitions and results of microlocal sheaf theory from the foundational work in \cite{Kashiwara-Schapira1}.
Then we recall previous results concerning the category $\Sh_\Lambda(M)$ of sheaves microsupported in a fixed subanalytic
isotropic closed subset from \cite{Ganatra-Pardon-Shende3}.

\subsection{Stable categories}

We will work in the higher categorical setting developed in \cite{Lurie1} and \cite{Lurie2}.
The main purpose of this subsection is to fix notations.
A thorough beginner guide as well as the set-up needed to work over a field $k$ of characteristic $0$
can be found in \cite[Chapter I.1]{Gaitsgory-Rozenblyum}.
To work over more general coefficients, one requires further the theory of rigid categories 
from \cite{Hoyois-Scherotzke-Sibilla}.
Since the relative case enjoys the same formal properties which are needed for the purpose of 
as the absolute case, we fix once and for all a rigid symmetric monoidal category $(\cV_0,\otimes,1_\cV)$
and its Ind-completion $\cV \coloneqq \Ind (\cV_0)$. 
Unless specified, we will assume without mentioning that all the categories we consider will be tensored over
and thus enriched in $\cV$ and functors between those are $\cV$-enriched as well.

One of the main advantages for working in the higher categorical setting is that
there is an abundance of limits and colimits (in an appropriate sense).
As a result, many constructions can be performed formally as universal constructions which greatly simplifies the situation. 
Because of the higher categorical nature of this paper, 
we will refer an $\infty$-category $\sC$ simply as a category
and when we need to emphasis that it is in particular an ordinary category, we will refer it as a $1$-category.

Recall that a \textit{presentable category} is a category with certain cardinality assumptions.
Roughly speaking, such a category is large enough to contain (small) colimits 
but is controlled by a small category. 
A main consequence of these assumptions is that the adjoint functor theorem holds.
See \cite[Corollary 5.5.2.9]{Lurie1}.
In addition, up to set-theoretic issues, the totality of such categories forms a (very large) category itself
which has nice properties concerning limits and colimits. 
We will not consider the whole collection of such categories but a small portion of it which satisfies stronger finiteness conditions
which we now recall.

\begin{definition}
Let $\sC$ be a category.
An object $c \in \sC$ is \textit{compact} if $\Hom(c,-)$ preserve (small) filtered colimit.
That is for any (small) filtered index category $I$ and any functor $X: I \rightarrow \sC$, the canonical morphism
$$  \varinjlim_I \Hom(c, X_i) \rightarrow \Hom(c, \varinjlim_I X_i)$$
is an isomorphism.
Here, we use the notation $\varinjlim\limits_I$ instead of $\clmi{I}$ to emphasis the index category $I$ is filtered.
\end{definition}

\begin{definition}
A category $\sC$ is \textit{compactly generated} if there exists a small subcategory $\sC_0 \subseteq \sC$
consisting of compact objects of $\sC$ such that $\sC$ is generated by $\sC_0$ under filtered colimits.
That is, $\Ind (\sC_0) \cong \sC$ where $\Ind$ denotes the Ind-completion.
\end{definition}

\begin{definition}
Let $\Cat$ denote the (very large) category of categories.
We use $\PrLcg$ to denote the (non-full) subcategory of $\Cat$ whose objects are compactly
generated categories and morphisms are functors which preserve small colimits and compact objects. 
We also use $\cat$ to denote the subcategory of 
$\Cat$ consisting of idempotent complete small categories which admit finite colimits whose 
morphisms are given by functors which preserve finite colimits.
\end{definition}

\begin{proposition}[{\cite[Proposition 5.5.7.8]{Lurie1}}]
The funtor $\Ind: \cat \rightarrow \PrLcg$
taking $\sC_0$ to its ind-completion $\Ind (\sC_0)$ is an equivalence. 
Its inverse is given by the functor $\theta: \PrLcg \rightarrow \cat$ sending $\sC$ to $\sC^c$, 
the subcategory of $\sC$ consisting compact objects.
\end{proposition}

\begin{proposition}[{\cite[Proposition 5.5.7.11]{Lurie1}}]\label{clmbrad}
The category $\PrLcg$ and hence $\cat$ admits small colimits, which can be computed
in $\Cat$ as limits by passing to right adjoints. 
Here we use the fact that a morphism $F: \sC \rightarrow \sD$ in $\PrLcg$
is a left adjoint since it preserves colimits. 
\end{proposition}

Classically, one use the theory of triangulated categories to encode homological information.
They are $1$-categories with structures and can be used to remember a small portion of homotopies. 
However, limits and colimits are scarce in this setting.
For example, the kernel of the (unique) non-zero morphism $e: \ZZ/2 \rightarrow \ZZ/2[1]$ 
in $D(\ZZ)$ does not exist. See, for example, \cite[2.2.1]{Toen}.
Hence, we use the theory of stable categories \cite[Chapter 1]{Lurie2} in this paper instead.

\begin{definition}
A category $\sC$ is \textit{pointed} if there exists a zero object $0$, i.e., an object which is both initial and final.
\end{definition}

A sequence $X \rightarrow Y \rightarrow Z$ in a pointed category $\sC$ is a \textit{fiber} 
(resp. \textit{cofiber}) \textit{sequence} if the diagram
 
$$
\begin{tikzpicture}

\node at (0,2) {$X$};
\node at (2,2) {$Y$};
\node at (0,0) {$0$};
\node at (2,0) {$Z$};

\draw [->, thick] (0.3,0) -- (1.7,0) node [midway, above] {$ $};
\draw [->, thick] (0.3,2) -- (1.7,2) node [midway, above] {$ $};

\draw [->, thick] (0,1.7) -- (0,0.3) node [midway, left] {$  $}; 
\draw [->, thick] (2,1.7) -- (2,0.3) node [midway, left] {$ $};

\end{tikzpicture}
$$

is a pullback/pushforward. 
In this case, we say $X$ (resp. $Z$) is the \textit{fiber} (resp. \textit{cofiber}) 
of the corresponding morphism $Y \rightarrow Z$ (resp. $X \rightarrow Y$).

\begin{definition}
A pointed category $\sC$ is \textit{stable} if fibers and cofibers exist and
a diagram as above is a fiber sequence if and only if it is a cofiber sequence.
We say a functor $F: \sC \rightarrow \sD$ between stable categories
is \textit{exact} if $F$ preserves finite limits and finite colimits.
Note in a stable category $\sC$ finite limits are the same as finite colimits
so preserving one kind means preserving the other.
\end{definition}

\begin{remark}
We recall that in the stable case, an object $X$ is compact if and only if $\Hom(X,-)$ preserves coproducts.
\end{remark}

\begin{example}
A stable category $\sC$ admits a ``\textit{shifting} by $1$" automorphism $[1]: \sC \rightarrow \sC$
which can be defined by $X \mapsto \cof(X \rightarrow 0)$.
Its inverse is the shifting by $-1$ automorphism $[-1]: \sC \rightarrow \sC$
which can defined by $X \mapsto \fib(0 \rightarrow X)$.
In the case when $\sC = \cV$, these shiftings are given by shifting the degree of the chain complexes.
\end{example}

\begin{example}
Let $\sC$ be a stable category. 
For $X$, $Y \in \sC$, the direct sum $X \oplus Y$ in $\sC$
can be computed as $\cof(Y[-1] \xrightarrow{0} X) = \fib(Y \xrightarrow{0} X[1])$. 
\end{example}

We will use the following lemma: 
\begin{lemma}\label{stsl}
Let $\sC$ be a stable category and $X$, $Y$, $Z$, $X^\prime$, $Y^\prime$, and $Z^\prime \in \sC$.
Assume we have the following commutative diagram 
$$
\begin{tikzpicture}

\node at (0,2) {$X$};
\node at (2.5,2) {$Y$};
\node at (5,2) {$Z$};
\node at (0,0) {$X^\prime$};
\node at (2.5,0) {$Y^\prime$};
\node at (5,0) {$Z^\prime$};

\draw [->, thick] (0.4,2) -- (2.1,2) node {$$};
\draw [->, thick] (2.8,2) -- (4.7,2) node {$$};
\draw [->, thick] (0.4,0) -- (2.1,0) node {$ $};
\draw [->, thick] (2.8,0) -- (4.7,0) node {$$};

\draw [->, thick] (0,1.7) -- (0,0.3) node [midway, right] {$\alpha$};
\draw [->, thick] (2.5,1.7) -- (2.5,0.3) node [midway, right] {$\beta$};
\draw [->, thick] (5,1.7) -- (5,0.3) node [midway, right] {$\gamma$}; 

\end{tikzpicture}
$$
such that each row is a fiber sequence.
Let $X^{\prime \prime} = \cof(\alpha)$, $Y^{\prime \prime} = \cof(\beta)$ and $Z^{\prime \prime} = \cof(\gamma)$
be the corresponding cofibers of the vertical maps.
Then there exist a canonical fiber sequence $X^{\prime \prime} \rightarrow Y^{\prime \prime} \rightarrow Z^{\prime \prime}$.

In particular, for $f_1: X_1 \rightarrow X_2$ and $f_2: X_2 \rightarrow X_3$, we have a fiber sequence
$\cof(f_1) \rightarrow \cof(f_2 \circ f_1) \rightarrow \cof(f_2)$.
This special case is usually referred as the octahedral axiom in the setting of triangulated categories. 
\end{lemma}

\begin{proof}
The goal is to show that $Z^{\prime \prime}$ is the cofiber of $(X^{\prime \prime} \rightarrow Y^{\prime \prime})$.
Recall that cofibers are computed as a colimit of the diagram $I = [\cdot \leftarrow \cdot \rightarrow \cdot]$.
For example, the object $Z$ is computed as the colimit given by the diagram $0 \leftarrow X \rightarrow Y$.
Now we consider the following diagram.

$$
\begin{tikzpicture}

\node at (0,4) {$0$};
\node at (2,4) {$Y$};
\node at (4,4) {$Y^\prime$};
\node at (0,2) {$0$};
\node at (2,2) {$X$};
\node at (4,2) {$X^\prime$};
\node at (0,0) {$0$};
\node at (2,0) {$0$};
\node at (4,0) {$0$};

\draw [->, thick] (1.5,4) -- (0.3,4) node {$ $};
\draw [->, thick] (2.5,4) -- (3.6,4) node {$$};
\draw [->, thick] (1.5,2) -- (0.3,2) node {$$};
\draw [->, thick] (2.5,2) -- (3.6,2) node {$$};
\draw [->, thick] (1.5,0) -- (0.3,0) node {$ $};
\draw [->, thick] (2.5,0) -- (3.6,0) node {$$};

\draw [->, thick] (0,1.7) -- (0,0.3) node {$ $};
\draw [->, thick] (2,2.3) -- (2,3.7) node {$ $};
\draw [->, thick] (4,1.7) -- (4,0.3) node {$ $}; 
\draw [->, thick] (0,2.3) -- (0,3.7) node {$ $};
\draw [->, thick] (2,1.7) -- (2,0.3) node {$ $};
\draw [->, thick] (4,2.3) -- (4,3.7) node {$ $}; 

\end{tikzpicture}
$$
Taking colimit for the vertical arrows gives $0 \leftarrow Z \rightarrow Z^\prime$ and taking the colimit again gives $Z^{\prime \prime}$.
Similarly, taking first the horizontal arrows and then the vertical arrows gives $\cof( X^{\prime \prime} \rightarrow Y^{\prime \prime})$.
But colimits commute with each other and thus $Z^{\prime \prime} =  \cof( X^{\prime \prime} \rightarrow Y^{\prime \prime})$.

Now for the special case, we apply the above result to the commutative diagram
$$
\begin{tikzpicture}

\node at (0,2) {$X_1$};
\node at (3,2) {$X_2$};
\node at (0,0) {$X_1$};
\node at (3,0) {$X_3$};

\draw [->, thick] (0.3,2) -- (2.7,2) node [midway, above] {$f_1$};
\draw [->, thick] (0.3,0) -- (2.7,0) node [midway, above] {$f_2 \circ f_1$};

\draw [double equal sign distance, thick] (0,1.7) -- (0,0.3) node {$ $};
\draw [->, thick] (3,1.7) -- (3,0.3) node [midway, right] {$f_2$};

\end{tikzpicture}
$$
\end{proof}

\begin{example} 
Consider two short exact sequences in a Grothendieck abelian $1$-category $\mathscr{A}$ 
and compatible maps between them as the following.

$$
\begin{tikzpicture}
\node at (-2.5,2) {$0$};
\node at (0,2) {$X$};
\node at (2.5,2) {$Y$};
\node at (5,2) {$Z$};
\node at (7.5,2) {$0$};
\node at (-2.5,0) {$0$};
\node at (0,0) {$X^\prime$};
\node at (2.5,0) {$Y^\prime$};
\node at (5,0) {$Z^\prime$};
\node at (7.5,0) {$0$};

\draw [->, thick] (-2.1,2) -- (-0.4,2) node {$$};
\draw [->, thick] (0.4,2) -- (2.1,2) node {$$};
\draw [->, thick] (2.8,2) -- (4.7,2) node {$$};
\draw [->, thick] (5.4,2) -- (7.1,2) node {$$};
\draw [->, thick] (-2.1,0) -- (-0.4,0) node {$$};
\draw [->, thick] (0.4,0) -- (2.1,0) node {$ $};
\draw [->, thick] (2.8,0) -- (4.7,0) node {$$};
\draw [->, thick] (5.4,0) -- (7.1,0) node {$$};

\draw [->, thick] (0,1.7) -- (0,0.3) node [midway, right] {$\alpha$};
\draw [->, thick] (2.5,1.7) -- (2.5,0.3) node [midway, right] {$\beta$};
\draw [->, thick] (5,1.7) -- (5,0.3) node [midway, right] {$\gamma$}; 

\end{tikzpicture}
$$
Recall that $\mathscr{A}$ naturally embeds into a stable category $D(A)$ 
which is usually referred as the Derived category of $\mathscr{A}$.
See for example \cite[Section 1.3]{Lurie2}.
An application of the above lemma implies a special case of the snake lemma in classical homological algebra.
\end{example}

We denote $\PrLcs$ and $\st$ the subcategory of $\PrLcg$ and $\cat$ which consists of stable categories. 
We recall that the property of being stable is compatible with the finiteness condition we discussed earlier.
In particular, the Ind-completion $\Ind (\sC_0)$ of a (small) stable category $\sC_0$ is stable.
Similarly, the subcategory of compact objects $\sC^c$ of a compactly generated stable category $\sC$
is stable.

\begin{proposition}\label{cs=st}
The equivalence $\Ind: \cat \leftrightarrows \PrLcg: \theta$ 
restricts to the (very large) subcategories consisting of stable categories 
$$\Ind: \st \leftrightarrows \PrLcs: \theta.$$
\end{proposition}

\subsection{General sheaf theory}

For a topological space $X$, the category of $\cV$-valued presheaves is the category $\PSh(X) \coloneqq \Fun(Op_X^{op},\cV)$ 
of contravariant functors from the $1$-category of open sets in $X$ to $\cV$.
The category of sheaves on $X$, $\Sh(X)$, is the reflexive subcategory of $\PSh(X)$ 
consisting of those presheaves $F$ which turn colimits in $Op_X$ to limits.
In more concrete terms, $F$ is a sheaf if for any open cover $\mathscr{U}$ of an open set $U \subseteq X$, 
the canonical map $F(U) \xrightarrow{\sim} \lim_{U_I \in C(\mathscr{U})} F(U_I)$ is an isomorphism, 
i.e., sections over $U$ can be computed as the totalization of the sections of the corresponding \v{C}ech nerve $C(\mathscr{U})$.
Recall that the term `reflexive' means that the inclusion $\Sh(X) \hookrightarrow \PSh(X)$ admits a left adjoint, 
which is usually referred as sheafification $(-)^\dagger$.
Thus, the inclusion is limit-preserving.
Since $\PSh(X)$ inherits limits and colimits from $\cV$, so does $\Sh(X)$
where colimits in $\Sh(X)$ can be computed as the sheafification of the colimits in $\PSh(X)$.

We also recall the six-functor formalism.
First, there is a symmetric monoidal structure $(\Sh(X),\otimes)$ on $\Sh(X)$ 
which is induced from $\cV$.
The unit of this tensor product is the sheaf $1_X$ which is the sheafification of the presheaf $(U \mapsto 1_\cV)$
whose restrictions are given by the identity $\id_1$.
For a fixed sheaf $F$,  the functor $(-) \otimes F$ 
which is given by tensoring with $F$  has a right adjoint $\sHom(F,-)$. 
This provides $\Sh(X)$ with an internal Hom.
The global section of this sheaf $\sHom$ is the $\cV$-valued $\Hom$, i.e., 
$\Gamma(X;\sHom(G,F)) = \Hom(G,F) \in \cV$ for any $F, G \in \Sh(X)$.

Let $f: X \rightarrow Y$ be a continuous map. 
There is a pushforward functor $f_*: \Sh(X) \rightarrow \Sh(Y)$
induced by pulling back open sets $f^{-1}: Op_Y \rightarrow Op_X$, $V \mapsto f^{-1} (V)$.
This functor admits a left adjoint $f^*: \Sh(Y) \rightarrow \Sh(X)$ and the adjunction $(f^*,f_*)$
is usually referred as the \textit{star pullback/pushforward}.
When $X$ and $Y$ are both locally compact Hausdorff spaces,
there is another pair of adjunction
$(f_!,f^!)$ such that $f_!: \Sh(X) \rightarrow \Sh(Y)$ and $f^!: \Sh(Y) \rightarrow \Sh(X)$.
This adjunction is usually referred as the \textit{shriek pullback/pushforward}.
When $f$ is proper, $f_!$ coincides with $f_*$.

\begin{remark}
Since a large portion of this paper is a sheaf-theoretic parallel of \cite{Ganatra-Pardon-Shende3},
we mention that their setting corresponds to the choice $\cV = \ZZ \dMod$, 
the presentable stable category of modules over $\ZZ$.
It can be modified as the dg category of (possibly unbounded) chain complexes of abelian groups with quasi-isomorphisms
inverted. See for example \cite{Cohn}.
This category is compactly generated and one usually denotes the compact objects $(\ZZ \dMod )^c$ by $\Perf \ZZ$. 
When representing $\ZZ \dMod$ by chain complexes,
$\Perf \ZZ$ consists of objects which are quasi-isomorphic to bounded chain complexes whose cohomology groups are finite rank.
\end{remark}

\begin{example}
When $i:Z \hookrightarrow X$ is a locally closed subset of $X$, 
one usually uses $F |_Z$ to denote $i^* F$ for $F \in \Sh(X)$ and call it the restriction of $F$ on $Z$.
\end{example}

\begin{example}
Consider a closed set $i: Z   \hookrightarrow X$ and an open set $j: U \hookrightarrow X$.
In these cases, we have $i_* = i_!$ and $j^* = j^!$.
In addition, the functors $i_*$, $j_*$, $j_!$ are fully faithful with the corresponding adjoints being a left inverse.
\end{example}

\begin{example}
Let $x \in X$ be a point.
For $F \in \Sh(X)$, we call the object $F_x \coloneqq F|_{\{x\}} \in \Sh(\{x\}) = \cV$ the \textit{stalk} of $F$ at $x$.
A key property of the stalks is that a morphism $G \rightarrow F$ is an isomorphism if and only if
it induces isomorphism $G_x \rightarrow F_x$ on the stalk at $x$ for all $x \in X$.
We also use the convention that $F$ has perfect stalk at $x$ if $F_x \in \cV_0$
or, equivalently, $F_x$ is perfect. 
We say that $F$ has perfect stalks if $F_x$ is perfect for all $x \in X$.
\end{example}

\begin{example}
Let $a_X: X \rightarrow \{*\}$ be the projection to a point. 
The object $a_X^! 1$ is usually denoted as $\omega_X$ and is referred as the dualizing complex/sheaf.
When $X$ is a $C^0$-manifold, $\omega_X$ is a locally constant sheaf whose stalk is given by $1_\cV[\dim X]$.
\end{example}

Now fix a topological space $X$. 
We see that taking integer coefficient sheaves itself forms a presheaf $\Sh$ in $\Cat$:
For an open set $U \subseteq X$, we assign the category $\Sh(U)$. 
For an inclusion of open sets $i_{U,V}: U \hookrightarrow V$, 
we assign the pullback functor $i_{U,V}^*: \Sh(V) \rightarrow \Sh(U)$
which is the right adjoint of ${i_{U,V}}_!: \Sh(U) \rightarrow \Sh(V)$.

\begin{proposition}\label{shissh}
The presheaf $\Sh: \Op_X \rightarrow \Cat$ is a sheaf.
\end{proposition}

\begin{proof}
Let $U$ be an open set of $X$ and $\mathscr{U}$ an open cover of $U$.
The functor $$\lmi{U_I \in C(\mathscr{U}) } i_{U_I,U}^* :\Sh(U) \rightarrow  \lim_{U_I \in C(\mathscr{U})} \Sh(U_I)$$
is an equivalence and its inverse $\lim_{U_I \in C(\mathscr{U})} \Sh(U_I) \rightarrow \Sh(U)$ can be described by
$$ \ (F_{U_I})_{U_I \in C(\mathscr{U}) } \mapsto \clmi{U_I \in C(\mathscr{U}) } \, (i_{U_I,U})_! F_{U_I}.$$
\end{proof}

We recall some standard properties of the six-functor formalism.
First, the assignment $X \mapsto \Sh(X)$, viewed as a functor on locally compact Hausdorff spaces, has a structure of \textit{base change}.
(See \cite{Jin1} for the exact statement in the higher categorical setting.)

\begin{theorem}\label{bsch}
Consider a pullback diagram of locally compact Hausdorff spaces,
$$
\begin{tikzpicture}

\node at (0,2) {$X^\prime$};
\node at (3,2) {$Y^\prime$};
\node at (0,0) {$X$};
\node at (3,0) {$Y$};

\draw [->, thick] (0.3,0) -- (2.7,0) node [midway, above] {$f$};
\draw [->, thick] (0.3,2) -- (2.7,2) node [midway, above] {$f^\prime$};

\draw [->, thick] (0,1.7) -- (0,0.3) node [midway, left] {$g^\prime$}; 
\draw [->, thick] (3,1.7) -- (3,0.3) node [midway, left] {$g$};

\end{tikzpicture}.
$$
There is an equivalence $g^* f_!  = f^\prime_! {g^\prime}^*$.
\end{theorem}

The push/pull functors satisfy some compatibility properties with $\otimes$ and $\sHom$.
We list a few which we will use:

\begin{proposition}\label{cp6f}
Let $f: X \rightarrow Y$ be a  continuous map between locally compact Hausdorff spaces. 
Then:
\begin{enumerate}
\item $f^* (F \otimes G) = f^* F \otimes f^* G$, for $F$, $G \in \Sh(Y)$,
\item $(f_! G) \otimes F = f_! (G \otimes f^* F)$, for $F$, $G \in \Sh(X)$,
\item $\sHom(f_! G, F) = f_* \sHom(G, f^! F)$, for $F \in \Sh(X)$, $G \in \Sh(Y)$,
\item $f^! \sHom(G,F) = \sHom(f^* F, f^! F)$, for $F$, $G \in \Sh(Y)$.
\end{enumerate}
\end{proposition}

We recall the excision fiber sequences.
Let $X$ be a locally compact Hausdorff space, $i: Z   \hookrightarrow X$ be a close set, and
$j: U = X \setminus Z   \hookrightarrow X$ be its open complement, then $j^* i_* = 0$
by base change and there are fiber sequences
$$ j_! j^! F \rightarrow F \rightarrow i_* i^* F$$
$$ i_! i^! F \rightarrow F \rightarrow j_* j^* F$$
where the arrows are the units/counits of the shriek/star adjunction pairs.
Such a triple $(\Sh(X),\Sh(Z),\Sh(U))$ is usually referred as  a \textit{recollement} in homological algebra.
See for example \cite[Section A.8]{Lurie2}.
When $Z$ is locally closed, one denotes $F_Z = i_! i^* F$ and $\Gamma_Z (F) = i_* i^! F$
for $F \in \Sh(X)$.
Thus, one can write the above fiber sequences as
$$ F_U \rightarrow F \rightarrow F_Z, \  \Gamma_Z (F) \rightarrow F \rightarrow \Gamma_U (F).$$


Let $a_X: X \rightarrow \{ * \}$ denote the projection to a point.
We use $A_X$ to denote the pullback $a_X^* A$ for $A \in \Sh( \{ * \}) = \cV$.
When $X$ is a manifold, $\cV = \ZZ \dMod$, and $A$ is an abelian group regarded as a chain complex concentrated as $0$, 
a standard representative of $A_X$ is the $A$-coefficient singular cochains $( U \mapsto C^*(U;A) )$.
When $Z \subseteq X$ is a locally closed subset of $X$, we abuse the notation and
write $A_Z$ for both the sheaf in $\Sh(Z)$ or its shriek pushforward.
When $A = 1_\cV$, we abuse the notation and write it simply as $1_Z$.

\begin{definition}
We call $A_Z$ the \textit{constant sheaf} on $Z$ with stalk $A$.
In general, we say a sheaf $F$ is a \textit{locally constant sheaf} or a \textit{local system} if there exists an open cover $\mathcal{U}$ such that
$F|_U$ is constant for $U \in \mathcal{U}$ and we use $\Loc(X)$ to denote the subcategory spanned by such sheaves.
\end{definition}

\begin{example}
Let $i: Z \hookrightarrow X$ be a local closed subset and $F \in \Sh(X)$.
By the above (1) and (2) of Proposition \ref{cp6f}, $F \otimes 1_Z  
= F \otimes i_! i^* 1_X =  i_! (i^* F \otimes i^* 1_X) = i_! i^* (F \otimes 1_X) = F_Z$.
A similar statement holds for open inclusions and the fiber sequence  $F_U \rightarrow F \rightarrow F_Z$ can be obtained 
from tensoring $F$ with the canonical one
$$ 1_U \rightarrow 1_X \rightarrow 1_Z.$$
\end{example}

We also consider set-theoretic invariants associated to sheaves.

\begin{definition}
Let $X$ be a topological space and $F \in \Sh(X)$.
The \textit{support} of a sheaf $F$ is defined to be the closed subset
$$ \supp(F) = \overline{ \{x \in X| F_x \neq 0 \} }.$$
\end{definition}

\begin{example}
Let $i: Z \hookrightarrow X$ be a closed subset.
The pushforward $i_*$ identifies $\Sh(Z)$ as the subcategory of $\Sh(X)$ 
consisting of sheaves $F$ whose support $\supp(F)$ is contained in $Z$. 
\end{example}

Before leaving this section, we recall a fundamental lemma for microlocal sheaf theory,
which holds for general Hausdorff spaces.
\begin{lemma}[ {\cite[Proposition 2.7.2]{Kashiwara-Schapira1} }, {\cite[Theorem 4.1]{Robalo-Schapira} }]\label{ncdef}
Let $X$ be a Hausdorff space, $F \in \Sh(X)$. Let $\{ U_s \}_{s \in \mathbb{R} }$ be a family of open subsets of $X$. We assume

\begin{enumerate}[label=(\alph*)]
\item for all $t \in \mathbb{R}, U_t = \bigcup_{s < t} U_{s}$,
\item for all pairs $(s,t)$ with $s \leq t$, the set $\overline{U_t \setminus U_s} \cap \supp(F)$ is compact,
\item setting $Z_s = \cap_{t > s} \overline{U_t \setminus U_s}$, we have for all pairs $(s,t)$ with $s \geq t$,
and all $x \in Z_s \setminus U_t$,
 $$(\Gamma_{X \setminus U_t} F)_x = 0.$$
\end{enumerate}

Then we have the isomorphism in $\Sh(X)$, for all $t \in \mathbb{R}$,
$$\Gamma( \bigcup_s U_s;F) \xrightarrow{\sim} \Gamma(U_t,F).$$
\end{lemma}

\subsection{Microlocal sheaf theory}

Now let $M$ be a $C^\alpha$-manifold where $\alpha \in \ZZ_{>0} \cup \{\infty, \omega\}$.
The term `microlocal' usually refers to `local' in the cotangent bundle $T^* M$.
In \cite[Section 5.1]{Kashiwara-Schapira1}, Kashiwara and Schapira define the notion of microsupport, 
which is a set in $T^* M$ enhancing the support.
One description of it is the following:
Let $F$ be a sheaf and $\phi$ a $C^1$ function on $M$, let $m \in \phi^{-1}(t)$.
We denote by $i_{\phi,t}: \{ x | \phi(x) \geq t \} \hookrightarrow M$ the closed inclusion. 
We say $m$ is a cohomological $F$-critical point of $\phi$ if $(i_{\phi,t}^! F)_m \neq 0$.

\begin{definition}
The \textit{microsupport} of a sheaf $F$ is defined to be the closure of the locus of differentials
of the $C^1$ functions at their cohomological $F$-critical points.
That is,
$$ \ms(F) = \overline{ \bigcup_{\phi \in C^1(M)} \{(x,\xi)| \exists t \in \mathbb{R}, (i_{\phi,t}^! F)_x \neq 0, \xi = d\phi_x  \} }.$$
\end{definition}

Although the microsupport is defined as a $C^1$-invariant, it is sufficient to check a smaller class of functions.

\begin{proposition}[{\cite[Proposition 5.1.1]{Kashiwara-Schapira1}}]
The microsupport of a sheaf $F$ is the same as the closure of the locus of differentials
of the $C^\alpha$ functions at their cohomological $F$-critical points.
That is,
$$ \ms(F) = \overline{ \bigcup_{\phi \in C^\alpha(M)} \{(x,\xi)| \exists t \in \mathbb{R}, (i_{\phi,t}^! F)_x \neq 0, \xi = d\phi_x  \} }.$$
\end{proposition} 

It is straightforward to see that the microsupport $\ms(F)$ of a sheaf $F \in \Sh(M)$ is conic and closed, 
and its intersection with the zero section $\ms(F) \cap 0_M = \supp(F)$ recovers the support.
The involutivity theorem \cite[Theorem 6.5.4]{Kashiwara-Schapira1} of Kashiwara and Schapira 
states that $\ms(F)$ is always a singular coisotropic.
Since $\ms(F)$ is conic, it can be recovered from $\supp(F)$ 
and its projectivization $\msif(F) \coloneqq (\ms(F) \setminus 0_M)/\Rp$.

\begin{definition}
For a closed conic subset $X \subseteq T^* M$, we use $\Sh_X(M)$ to denote the subcategory of sheaves consisting 
of those $F$ such that $\ms(F) \subseteq X$.
Similarly, for a closed subset $X \subseteq S^* M$, we use $\Sh_X(M)$ to denote the subcategory of sheaves consisting
of those $F$ such that $\msif(F) \subseteq X$.
Note for a closed $X \subseteq S^* M$, $\Sh_X(M) = \Sh_{(\Rp X \cup 0_M)}(M)$.
\end{definition}

\begin{example}\label{loceg}
Let $M$ be a manifold.
Being a local system is a microlocal condition. 
More precisely, $\Loc(M) = \Sh_{0_M}(M)$.
\end{example}

\begin{example}[{\cite[Proposition 5.3.1]{Kashiwara-Schapira1}}]\label{conems}
More generally, let $M = \mathbb{R}^n$, $\gamma$  be a closed convex cone with vertex at $0$,
and denote by $\gamma^\circ \coloneqq \{\xi \in (\mathbb{R}^n)^* | \xi(v) \geq 0, v \in \gamma \}$ its dual cone.
One has $\ms( 1_\gamma) \cap T^*_0 \mathbb{R}^n = \gamma^\circ$.
As a corollary, if $M^\prime \subseteq M$ is a closed submanifold, 
then $\ms( 1_{M^\prime} ) = N^* M^\prime $ is the normal bundle of $M^\prime $.
\end{example}

One might want to assign an invariant similar to the stalks for points in $(x,\xi) \in S^* M$.
In general, the object $(i_{\phi,t}^! F)_x$ depends on $\phi$ and is not an invariant associated to the point $(x,\xi)$.
However, the situation is better when cartain transversality condition is satisfied.

\begin{definition}
Fix a singular isotropic $\Lambda \subseteq S^* M$, i.e, $\Lambda$ is stratified by isotropic submanifolds.
Let $f$ be a function defined on some open set $U$ of $M$.
We say that a point $x \in U$ is a $\Lambda$-critical point of $f$ if the graph of its differential $\Gamma_{df}$ intersects 
$\Rp \Lambda \cup 0_M$ at $(x,df_x)$.
A $\Lambda$-critical point $x$ is Morse if $(x,d f_x)$ is a smooth point of $\Rp \Lambda \cup 0_M$ 
and the intersection $\Gamma_{df} \cap \Lambda$ is transverse at $(x,d f_x)$.
A function $f$ is $\Lambda$-Morse if all its $\Lambda$-critical point is Morse.
\end{definition}

\begin{proposition}[{\cite[Proposition 7.5.3]{Kashiwara-Schapira1}}]\label{mstk}
Let $\Lambda$ be a singular isotropic.
Assume $\phi$ is $\Lambda$-Morse at a smooth point $(x,\xi) \in \Lambda$.
For $F \in \Sh(M)$ such that $\msif(F) \subseteq \Lambda$ in a neighborhood of $(x,\xi)$,
the object $(i_{\phi,t}^! F)_m \in \cV$ is, up to a shift, independent of $\phi$.
\end{proposition}

\begin{definition}
Let $\Lambda \subseteq S^* M$ be a singular isotropic and $(x,\xi) \in \Lambda$ a smooth point.
For $F \in \Sh_\Lambda(M)$, we call functors $\mu_{(x,\xi)}: \Sh_\Lambda(M) \rightarrow \cV$ of the form
$$\mu_{(x,\xi)} F \coloneqq (i_{\phi,t}^! F)_x$$
\textit{microstalk functors} where $\phi$ is any function which satisfies the assumption in the last Proposition.
Since this functor is well-defined up to a shift, which will not play a significant role in the paper, 
we will abuse notation and call $\mu_{(x,\xi)} F$ the \textit{microstalk} of $F$ at $(x,\xi)$.
\end{definition}

Let $X \subseteq T^* M$ be conic and closed, and let $\Lambda \subseteq \left( T^* M \setminus X \right)$
be a closed conic subanalytic isotropic.
Assume $\ms(F) \subseteq X \cup \Lambda$.
To determine whether $(x,\xi) \in \Lambda$ is in the microsupport of $F$,
by definition one has to check if $(i_{\phi,t}^! F)_x$ vanishes
over all functions $\phi$ such that $x \in \phi^{-1}(t)$ and $d \phi_x = \xi$.
However, since $\Lambda$ is isotropic, it is sufficient to check the Morse ones.

\begin{proposition}[{\cite[Proposition 4.9]{Ganatra-Pardon-Shende3}}]\label{msbms}
Let $X$ and $\Lambda$ be as above.
Then $\Sh_X(M) \subseteq \Sh_{X \cup \Lambda}(M)$ is the fiber of all microstalk functors $\mu_{(x,\xi)}$
for smooth Lagrangian points $(x,\xi) \in \Lambda$.  
\end{proposition}

In practice, it is hard to compute the microsupport of a sheaf directly and it is usually
sufficient to deduce desired conclusions by having an upper bound.
Here we collect some standard results for microsupport estimation.
Let $f: M \rightarrow N$ be a map between manifolds, we use the following notations 

$$
\begin{tikzpicture}

\node at (0,2) {$ T^* M$};
\node at (5,2) {$ M \times_N T^* N$};
\node at (10,2) {$ T^* N$};
\node at (5,0) {$M$};
\node at (10,0) {$N$};
\node at (7.5,1) {$\square$};

\draw [<-, thick] (0.6,2) -- (3.9,2) node [midway, above] {$df^*$};
\draw [->, thick] (6,2) -- (9.4,2) node [midway, above] {$f_\pi$};
\draw [->, thick] (5.5,0) -- (9.7,0) node [midway, above] {$f$};


\draw [->, thick] (5,1.7) -- (5,0.3) node [midway, left] {$\pi_M $};
\draw [->, thick] (10,1.7) -- (10,0.3) node [midway, left] {$  \pi_N$}; 

\end{tikzpicture}
$$
where the square on the right is the pullback of the cotangent bundle $T^* N$ of $N$ along $f$
and $df^*$ is given fiberwisely by the adjoint of the differential $df_x:T_x M \rightarrow T_{f(x)} N$.
Let  $T^*_M N$ denote the set 
$$\{(x,\alpha) \in M \times_N T^* N | df_x^* \alpha = \alpha \circ d f_x = 0 \}.$$
Note in case $M \subseteq N$ is a submanifold,
$T^*_M N$ is the conormal bundle of $M$ in $N$.

\begin{definition}[{\cite[Definition 5.4.12]{Kashiwara-Schapira1}}]
Let $A$ be a closed conic subset of $T^* N$.
We say $f$ is \textit{noncharacteristic} for $A$ if
$$ f_\pi^{-1}(A) \cap T_M^* N \subseteq M \times_N 0_N.$$
For a sheaf $F \in \Sh(M)$, we say $f$ is \textit{noncharacteristic} for $F$ if it is the case for $\ms(F)$.
\end{definition}

\begin{proposition}\label{mses}
We have the following results:
\begin{enumerate}
\item (\cite[Proposition 5.1.3]{Kashiwara-Schapira1}) If $F \rightarrow G \rightarrow H$ is a fiber sequence in $\Sh(M)$, then
$$ \left( \ms(F) \setminus \ms(H) \right) \cup \left( \ms(H) \setminus \ms(F) \right)
\subseteq \ms(G) \subset \ms(F) \cup \ms(H).$$
This is usually referred as the microlocal triangular inequalities.

\item (\cite[Proposition 5.4.1]{Kashiwara-Schapira1}) 
For $F \in \Sh(M)$, $G \in \Sh(N)$, $\ms(F \boxtimes G) \subseteq \ms(F) \times \ms(G)$.

\item (\cite[Proposition 5.4.4]{Kashiwara-Schapira1}) 
For $f: M \rightarrow N$ and $F \in \Sh(M)$, if $f$ is proper on $\supp(F)$, then $$\ms(f_* F) 
\subseteq f_\pi \left( (d f^*)^{-1} \ms(F) \right).$$

\item (\cite[Proposition 5.4.5 and Proposition 5.4.13]{Kashiwara-Schapira1}) 
For $f: M \rightarrow N$ and $F \in \Sh(N)$,  if $f$ is noncharacteristic for $F$, then
$$\ms(f^* F) \subseteq df^* \left( f_\pi^{-1} \left(\ms(F) \right) \right)$$
and the natural map $f^* F \otimes f^! 1_Y \rightarrow f^! F$ is an isomorphism.
If $f$ is furthermore smooth, the estimation is an equality.

\item  (\cite[Proposition 5.4.8]{Kashiwara-Schapira1}) 
Let $Z \subseteq M$ be closed. 
If $\ms(F) \cap N^*_{out} (Z) \subseteq 0_M$, then 
$$\ms(F_Z) \subseteq N^*_{in}(Z) + \ms(F).$$
Similarly, let $U \subseteq M$ be open. 
If $\ms(F) \cap N^*_{in} (U) \subseteq 0_M$, then 
$$\ms(F_U) \subseteq N^*_{out}(U) + \ms(F).$$

\item (\cite[Proposition 5.4.14]{Kashiwara-Schapira1}) 
For $F$ and $G \in \Sh(M)$. If $\ms(F) \cap -\ms(G) \subseteq 0_M$, then
$$ \ms(F \otimes G) \subseteq \ms(F) + \ms(G).$$

\item  (\cite[Proposition 5.4.14 and Exercise V.13]{Kashiwara-Schapira1}) 
For $F$ and $G \in \Sh(M)$. If $\ms(F) \cap \ms(G) \subseteq 0_M$, then
$$ \ms( \sHom(G,F) ) \subseteq \ms(F) - \ms(G).$$
If moreover $G$ is cohomological constructible, then
the natural map $$\sHom(G,1_M) \otimes F \rightarrow \sHom(G,F)$$ is an isomorphism.
If furthermore $F = \omega_M$, i.e., when $\sHom(G,F) \eqqcolon D_M (G)$ is the Verdier dual,
then $\ms(D_X (G) ) = - \ms(G)$.
\end{enumerate}

\end{proposition}

Sometimes, when the noncharacteristic condition is absent, there is still a less refined upper bounds for pullbacks.
\begin{definition}[{\cite[Definition 6.2.3]{Kashiwara-Schapira1}}]\label{ntb}
We define two constructions of closed conic subsets of cotangent bundles:
\begin{enumerate}
\item Given closed conic subsets $A, B \subseteq T^* M$, we define $A \hat{+} B$ to be the closed subset
consisting of points $(x,\xi) \in T^* M$ such that, in some local coordinate, there exist sequence
$\{(x_n,\xi_n)\}$ in $A$ and $\{(y_n,\eta_n)\}$ in $B$ such that $x_n, y_n \rightarrow x$, 
$\xi_n + \eta_n \rightarrow \xi$, and $\abs{x_n - y_n} \abs{\xi_n } \rightarrow 0$.  
\item Let $i: X \hookrightarrow M$ be a closed submanifold and choose a local coordinate $(x,y,\xi,\eta)$ of $T^* M$
such that $X$ is given by $\{ y = 0\}$.
Given a closed conic subset $A \subseteq T^* M$, we define $i^\# (A)$ to be the closed subset of $T^* X$
consisting of points $(x,\xi)$ such that there exists $\{ (x_n,y_n,\xi_n,\eta_n) \}$ in $A$ such that
$y_n \rightarrow 0$, $x_n \rightarrow x$, $\xi_n \rightarrow \xi$, and $\abs{y_n} \abs{\eta_n} \rightarrow 0$.
\end{enumerate}
\end{definition}

In general for $f: M \rightarrow N$ and closed conic subset $A \subseteq T^* N$,
$f^\#(A)$ can be defined as a special case of a more general construction 
which also includes $A \hat{+} B$ as a special case.
Moreover, the definition can be made free of choice of local coordinates using the technique of
\textit{deformation to the normal cone} \cite[Section 4.1, Section 6.2]{Kashiwara-Schapira1}.
Also, notice that the above construction $A \hat{+} B$ and $i^\# (A)$ contain
$A + B$ and $di^* ( i_\pi^{-1}(\ms(F)) )$ as closed subsets.
It can be shown that they are equal when the noncharacteristic condition is satisfied.

\begin{proposition}\label{ncmses}
We have the following results:
\begin{enumerate}
\item (\cite[Theorem 6.3.1]{Kashiwara-Schapira1}) 
Let $j: U \hookrightarrow M$ be open and $F \in \Sh(U)$, then $$\ms(j_* F) \subseteq \ms(F) \hat{+} N^*_{in} U.$$

\item (\cite[Corollary 6.4.4]{Kashiwara-Schapira1}) 
Let $i: Z \rightarrow M$ be a closed submanifold and $F \in \Sh(M)$, 
then $$\ms(i^* F) \subseteq i^\#  \ms(F).$$

\item (\cite[Corollary 6.4.5.]{Kashiwara-Schapira1})
Let $F$, $G \in \Sh(M)$, then
$$ \ms(F \otimes G) \subseteq \ms(F) \hat{+} \ms(G), 
\ \ms\left( \sHom(G,F) \right) \subseteq \ms(F) \hat{+} \left(-\ms(G) \right).$$
\end{enumerate}
\end{proposition}

The notion of $\hat{+}$ can also be used to measure compatibility between $!$-pullback and $\otimes$.
\begin{proposition}[{\cite[Exercise VI.4]{Kashiwara-Schapira1}}]\label{shktn}
Let $f: M \rightarrow N$ be a morphism of manifolds and let $F, G \in \Sh(N)$.
If $f$ is non-characteristic for $\ms(F) \hat{+} \ms(G)$, then
$$f^*F \otimes f^! G = f^! (F \otimes G).$$
\end{proposition}

As a Corollary, we mention a computational tool:
\begin{lemma}\label{omegap}
Let $M$ and $N$ be manifolds. Then
\begin{enumerate}
\item $p_N^! 1_N = p_M^* \omega_M$,
\item $\omega_{M \times N} = \omega_M \boxtimes \omega_N$.
\end{enumerate}
As a Corollary, $\omega_M$ is invertible whose inverse can be given as $\Delta^! 1_{M \times M}$.
\end{lemma}

\begin{proof}
Consider the pullback diagram:
$$
\begin{tikzpicture}
\node at (0,2) {$M \times N$};
\node at (4,2) {$N$};
\node at (0,0) {$M$};
\node at (4,0) {$\{*\}$};

\draw [->, thick] (0.7,2) -- (3.7,2) node [midway, above] {$p_N$};
\draw [->, thick] (0.3,0) -- (3.6,0) node [midway, above] {$a_M$};

\draw [->, thick] (0,1.7) -- (0,0.3) node [midway, right] {$p_M$}; 
\draw [->, thick] (4,1.7) -- (4,0.3) node [midway, right] {$a_N$};
\end{tikzpicture}
$$ 

For (1), base change $a_M^* {a_N}_!  = {p_N}_!  p_M^*$ implies that there exists
a canonical map $p_M^* a_N^! \rightarrow p_N^! a_M^*$.
This map is in general not an isomorphism but in our case,
we may assume $M$ and $N$ are Euclidean spaces by checking the map locally.
Then the isomorphism follows from the isomorphism $1_\cV = \Gamma_c(\RR^k;\cV)[k]$
and $\omega_{\RR^k} = 1_{\RR^k}[k]$.
For (2), we can use (1) and (4) of Proposition \ref{mses} and compute that 
$$\omega_M \boxtimes \omega_N = p_M^* \omega_M \otimes p_N^* \omega_N 
= p_M^* \omega_M \otimes p_M^! 1_M = p_M^! \omega_M = \omega_{M \times N}.$$
To obtain the Corollary, we apply Proposition \ref{shktn} and compute that 
$$ 
\Delta^!( 1_{M \times M}) \otimes \omega_M 
= \Delta^!( 1_{M \times M}) \otimes \Delta^* (p_1^* \omega_M )
= \Delta^! (p_1^* \omega_M) = \Delta^! p_2^! (1_M) = 1_M.
$$
\end{proof}

One can combine the six-functors to build more general functors between sheaves on topological spaces.
One example is the following technique of convolution, which we will use extensively in the later sections.
Let $X_i$, $i = 1, 2, 3$, be locally compact Hausdorff topological spaces,
and write $X_{ij} = X_i \times X_j$, for $i < j$,  $X_{123} = X_1 \times X_2 \rightarrow X_3 $,
and $\pi_{ij}: X_{123} \rightarrow X_{ij}$ for the corresponding projections.
For $F \in \Sh(X_{12})$, $G \in \Sh(X_{23})$, the \textit{convolution} is defined to be 
$$G \circ_{M_2} F \coloneqq {\pi_{13}}_! (\pi_{23}^* G \otimes \pi_{12}^* F ) \in \Sh(X_{13}).$$
When there is no confusion what $X_2$ is, we will usually suppress the notation and simply write it as $G \circ F$.
This is usually the case when $X_1 = \{*\}$, $X_2 = X$, and $X_3 = Y$ and we think of $X$ as 
the source and $Y$ as the target, $G \in \Sh(X \times Y)$ as a functor sending $F \in \Sh(X)$ to $G \circ (F) \in \Sh(Y)$.
Note that from its expression, this functor is colimit-preserving. 
 
\begin{lemma}[{\cite[Proposition 3.6.2]{Kashiwara-Schapira1}}]\label{convrad}
For fixed $G \in \Sh(X_{23})$, 
the functor $G \circ (-) : \Sh(X_{12}) \rightarrow \Sh(X_{13})$ induced by convolving with $G$
has a right adjoint which we denote it as $\sHom^\circ(G,-): \Sh(X_{13}) \rightarrow \Sh(X_{12})$ 
and is given by
\begin{equation}\label{cvr}
H \mapsto  {\pi_{12}}_* \sHom(\pi_{23}^* G  ,\pi_{13}^! H).
\end{equation}
\end{lemma}

\begin{example}
We note that convolution recovers $*$-pullback and $!$-pushforward.
For example, let $f: X \rightarrow Y$ be a continuous map and denote by $i: \Gamma_f \subseteq X \times Y$ its graph.
Take $X_1 = \{*\}$, $X_2 = X$, and $X_3 = Y$, then for $F \in \Sh(X)$,
$$1_{\Gamma_f} \circ F = {\pi_Y}_! ( 1_{\Gamma_f} \pi_X^* F) = {\pi_Y}_! i_! i^* \pi_X^* F = f_! F.$$
\end{example}

We note that base change implies that convolution satisfies associativity.
\begin{proposition}\label{conv:asso}
Let $F_i \in \Sh(X_{i i +1})$ for $i = 1, 2, 3$. Then
$$ F_3 \circ_{X_3} (F_2 \circ_{X_2} F_1) = (F_3 \circ_{X_3} F_2) \circ_{X_2} F_1.$$
In particular, if $G_1$, $G_2 \in \Sh(X \times X)$, then there is an identification of functors
$$ G_2 \circ (G_1 \circ (-) ) = (G_2 \circ_X G_1) \circ (-).$$
\end{proposition}

We will use a relative version of convolution. 
Let $B$ be a locally compact Hausdorff space viewed as a parameter space. 
Regard $F \in \Sh(X_{12} \times B)$, $G \in \Sh(X_{23} \times B)$ as $B$-family sheaves,
one can similarly define the relative convolution $G \circ|_B F \in \Sh(X_{13} \times B)$
by replacing $\pi_{ij}$ with 
$$\pi_{{ij},B}: X_{123} \times B \rightarrow X_{ij} \times B.$$

In the case of manifolds, convolution satisfies certain compatibility with microsupport.
For $A \subseteq T^* M_{12}$ and $B \subseteq T^* M_{23}$, we set
$$B \circ A = \{ (x,\xi,z,\zeta) \in T^* M_{13} | \exists (y,\eta), (x,\xi,y,\eta) \in A, (y,-\eta, z, \zeta) \in B \}.$$
Note if $A$ and $B$ are Lagrangian correspondences satisfying appropriate transversality condition,
the set $B \circ A$ is the composite Lagrangian correspondence twisted by a minus sign on the second component.
Write $q_{ij}: T^* M_{123} \rightarrow T^* M_{ij}$ to be the projection on the level of cotangent bundles
and $q_{2^a3}$ the composition of $q_{23}$ with the antipodal map on $T^* M_2$.
Then $B \circ A = q_{13} (q_{2^a3}^{-1}  B \cap q_{12}^{-1} A)$ and (4), (6) and (3) of Proposition \ref{mses} implies 
the following corollary.
 
\begin{corollary}[{\cite[(1.12)]{Guillermou-Kashiwara-Schapira}}]\label{msconv}

Assume 
\begin{enumerate}
\item $p_{13}$ is proper on $M_1 \times \supp(G) \cap \supp(F) \times M_3$;
\item $q_{2^a 3}^{-1} \ms(G)  \cap q_{12}^{-1} \ms(F)  \cap 0_{M_1} \times T^* M_2 \times 0_{M_3}
\subseteq 0_{M_{123}}.$ 
\end{enumerate}
$$\ms( G \circ F ) \subseteq \ms(G) \circ \ms(F).$$
\end{corollary}

A similar microsupport estimation holds for the $B$-family case.
One noticeable difference for the microsupport estimation is that instead of $T^* M_{ij}$ and $T^* M_{123}$
one has to consider $T^* M_{ij} \times T^* B$ and $T^* M_{123} \times (T^* B \times_B T^* B)$ instead.
Here $\times_B$ is taken over the diagonal $B \hookrightarrow B \times B$.
Also the projection $r_{ij}: T^* M_{123} \times (T^* B \times_B T^* B) \rightarrow T^* M_{ij} \times T^* B$
for the $B$-component is now given by the first projection (with a minus sign) for $ij =12$, the addition for $ij = 13$,
and the second projection $ij = 23$.
Otherwise the microsupport estimation is similar to the ordinary case.

Finally, we mention the compatibility between microsupport and limits/colimits.
\begin{proposition}[{\cite[Exercise V.7]{Kashiwara-Schapira1},  \cite[2.7]{Jin-Treumann}}]
Let $I$ be a (small) set and $\{F_i\}_{i \in I}$ be a family of sheaves on $M$ index by $I$. 
Then there are microsupport estimations,
$$\ms( \bigoplus_i F_i) \subseteq \overline{ \cup_i \ms(F_i)}, \
\ms(\prod_i F_i) \subseteq \overline{ \cup_i \ms(F_i) }.$$
\end{proposition}

Let $X \subseteq T^* M$ be a conic closed subset, the above microsupport estimation 
combined with the adjoint functor theorem \cite[Corollary 5.5.2.9]{Lurie1} implies that
the inclusion $${\iota_X}_*: \Sh_X(M) \subseteq \Sh(M)$$ admits both a left adjoint $\iota_X^*$ and a right adjoint $\iota_X^!$.
Indeed, to check that ${\iota_X}_*$ preserves limits (resp. colimits), it is enough to check for $\prod_i$ (resp. $\bigoplus$)
since we are in the stable setting.
Now let $F_i \in \Sh_X(M)$ for $i \in I$, the above proposition implies that $\bigoplus_i F_i$ and $\prod_i F_i$, taken in $\Sh(M)$,
both have microsupport contained in $\overline{ \cup_i \ms(F_i)} \subseteq \overline{\cup_i X} = X$ and are objects in $\Sh_X(M)$.
General categorical theory then implies that the category $\Sh_X(M)$ is itself presentable.
We will see in the next two sections that, under some mild regularity condition,
the category $\Sh_\Lambda(M)$ is in fact compactly generated where $\Lambda \subseteq S^* M$ is a singular isotropic.

\subsection{Constructible sheaves}

The theory of constructible sheaves is based on the results of stratified spaces.
Standard references for stratified spaces are \cite{Goresky-MacPherson} and \cite{Mather}.

A \textit{stratification} $\mathscr{S}$ of $X$ is a decomposition of $X$ into to a disjoint union of locally closed subset 
$\{ X_s \}_{ s \in \mathscr{S}}$.
A set $Y \subseteq X$ is said to be \textit{$\mathscr{S}$-constructible} if it is a union of strata in $\mathscr{S}$.
We assume, without further mention, that a stratification should be locally finite and satisfies the \textit{frontier
condition} that $\overline{X_s} \setminus X_s$ is a disjoint union of strata in $\mathscr{S}$.
In this case, there is an ordering which is defined by $s \leq  t$ if and only if $X_t \subseteq \overline{X_s}$.
We always implicitly chose this ordering when regarding $\mathscr{S}$ as a poset.
For example, we will consider its linearization $\mathscr{S} \dMod \coloneqq \PSh(\mathscr{S}^{op},\cV)$.

\begin{definition}
For a given stratification $\mathscr{S}$,
a sheaf $F$ is said to be $\mathscr{S}$-constructible if $F|_{X_s}$ is a local system for all
$s \in \mathscr{S}$.
We denote the subcategory of $\Sh(X)$ consisting of such sheaves by $\Sh_{\mathscr{S}}(X)$.
A sheaf $F$ is said to be constructible if $F$ is $\mathscr{S}$-constructible for some stratification $\mathscr{S}$.
\end{definition}

For $s \in \mathscr{S}$, we denote by $\str(s)$ the smallest $\mathscr{S}$-constructible open set containing $X_s$.
Alternatively, $\str(s) = \coprod_{t \leq s} X_t$.
The poset $\mathscr{S}$ can be then identified with the subposet $\{ \str(s) | s \in \mathscr{S}\}$ of $Op_X$.
Hence, there is a functor $\mathscr{S} \dMod \hookrightarrow \Sh_{\mathscr{S}}(X)$ induced 
by the restricting along the inclusion of poset $\mathscr{S} \hookrightarrow Op_X$.
A stratification is called a \textit{triangulation} if $X = \abs{K}$ is a realization of some simplicial complex $K$
and $\mathscr{S} \coloneqq \{ \abs{\sigma} | \sigma \in K \}$ is given by the simplexes of $K$.
We note that when $\mathscr{S}$ is a triangulation,
the functor $\mathscr{S} \dMod \hookrightarrow \Sh_{\mathscr{S}}(X)$ is an equivalence.
This follows from the fact that each $X_s$ is contractible and the following criterion:

\begin{lemma}[{\cite[Lemma 4.2]{Ganatra-Pardon-Shende3}}]
Let $\Pi$ be a poset with a map to $Op_M$, and let $\cV[\Pi]$ denote its stabilization.
The following are equivalent
\begin{itemize}
\item $\Gamma(U;\cV) = 1_\cV$ for $U \in \Pi$ and $\Gamma(U;\cV) \rightarrow \Gamma(U \setminus V;\cV)$
whenever $U \not \subseteq V$.
\item The composition $\cV[\Pi] \rightarrow \cV[Op_M] \rightarrow \Sh(M)$ is fully faithful 
where the second map is given by $!$-pushforward.
\end{itemize} 
\end{lemma}

Now let $M$ be a $C^\alpha$ manifold for $\alpha \in \NN_{> 0} \cup \{\infty, \omega\}$.
We consider regularity conditions for stratifications.
\begin{definition}
Let $M$ and $N$ be locally closed $C^1$ manifolds of $\mathbb{R}^n$, with $N \subset \overline{M} \setminus M$.
Consider sequences $x_n \in M$ and $y_n \in N$ such that $x_n, y_n \rightarrow y \in N$ and 
$\{T_{x_n} M\}$ converges to $\tau$ (in the corresponding Grassmannian).
Assume also $\{ \RR(x_n - y_n)\}$ converges to $l$.
We say the pair $(N,M)$ satisfies the \textit{Whitney condition} if any such sequence satisfies $\tau \supseteq l$.
In general, we say a pair $(N,M)$ of $C^1$ submanifold manifolds of a $C^\alpha$ manifold $X$
satisfies the Whitney condition if the above condition is satisfied on local charts.
\footnote{The definition of Whitney condition given here is more commonly referred as Whitney condition (b), 
which implies a weaker Whitney condition (a). We will use neither of the conditions directly and thus this simplified choice.}
\end{definition}

\begin{remark}[{\cite[Chapter 4]{Mather}}]
The Whitney condition can also be formulated without working on coordinates.
Recall the normal bundle of the diagonal $\Delta_M \hookrightarrow M \times M$ can be identified 
as the tangent bundle $TM$ of $M$.
The real blow-up $Bl_{\Delta_M} (M \times M)$ can be seen as a disjoint union 
$\mathbb{P}(TM) \amalg (M \times M \setminus \Delta_M)$ of the projective tangent bundle and the off-diagonal.
Then we say  $(N_1,N_2)$ satisfies the Whitney condition if for any sequence 
$(y_n,x_n) \in N_1 \times N_2 \subseteq Bl_{\Delta_M} (M \times M) \setminus \Delta_M$
such that $T_{x_n} N_2 \rightarrow \tau$ and $(y_n,x_n) \rightarrow l \in \mathbb{P}(TM)$, we have $l \subseteq \tau$.
\end{remark}

We say a stratification $\mathscr{S}$ is $C^k$ if each $X_s$ is $C^k$ locally closed manifold.  
A $C^k$ stratification $\mathscr{S}$ is a Whitney stratification if $(X_t,X_s)$ satisfies the Whitney condition for $s \leq t$.
Let $N^* \mathscr{S}$ denote the union $\cup_{s \in \mathscr{S}} N^* X_s$ of the conormals of the strata.
The set $N^* \mathscr{S}$ is a singular conic Lagrangian in $T^* M$ and
the Whitney condition implies a weaker property that $N^* \mathscr{S} \subseteq T^* X$ is closed.
We also use $N^*_\infty \mathscr{S}$ to denote the corresponding singular isotropic at the infinity.
The main advantage of considering Whitney stratifications is the following proposition.

\begin{proposition}[{\cite[Prop. 8.4.1]{Kashiwara-Schapira1}, \cite[Proposition 4.8]{Ganatra-Pardon-Shende3}}]
For a Whitney stratification $\mathscr{S}$ of a $C^1$ manifold $M$, we have $\Sh_{\mathscr{S}}(M) = 
\Sh_{N^*_\infty \mathscr{S}}(M)$ (i.e. having microsupport contained in $N^* \mathscr{S}$ is equivalent to being
$\mathscr{S}$-constructible).
\end{proposition}

\begin{example}
We note that the Whitney condition is crucial for this proposition.
More precisely, we will give an example, on $\RR^3$, 
of a stratification $\sS$ and a sheaf $F$ such that $F \in \Sh_{N^*_\infty \mathscr{S}}(\RR^3)$ but $F \not\in \Sh_{\sS}(\RR^3)$.
Consider the $C^\infty$ map
\begin{align*}
f: \RR^3 &\rightarrow \RR^3 \\
(x,y,z) &\mapsto (z^2 x,z^2 y,z).
\end{align*}
Set $V \coloneqq \{ x^2 + y^2 = z^4\}$, which we note is the image of $C \coloneqq \{x^2 + y^2 = 1\}$ under $f$.
Take the stratification $\sS = \{X_1, X_2, X_3\}$ of $\RR^3$ where $X_1 \coloneqq \{x = z^2, y = 0 \}$,
$X_2 \coloneqq V \setminus X_1$, and $X_3 \coloneqq \RR^3 \setminus V$. 
We note that $\sS$ does not satisfy the Whitney condition at the origin.

\begin{figure}[H]
  \centering
    \includegraphics[scale = 0.12]{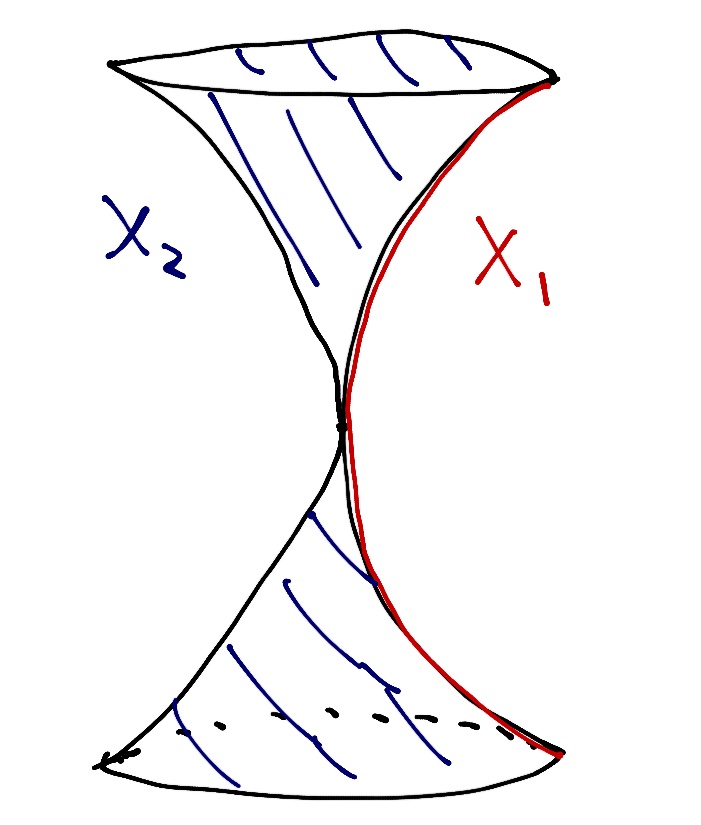}
  \caption{A stratification of $\RR^3$ which provides a counterexample of the above proposition when
  the Whitney condition is not present. The red curve is the single $1$-dimensional stratum $X_1$, 
  the black locally closed surface is the single $2$-dimensional stratum $X_2$, 
  and their complement is the single $3$-dimensional open stratum $X_3$.}
\end{figure}

Now consider the sheaf $F \coloneqq f_* (1_C)$.
We claim that $\ms(F) \subseteq \left(N^* X_1 \cup N^* X_2 \right)$ and we have $F \in \Sh_{N^*_\infty \mathscr{S}}(\RR^3)$.
The claim is clear when $z \neq 0$ since $f|_{ \{z\neq 0\}}$ is a diffeomorphism.
For $z = 0$, since $f$ is proper on $C$, we can apply the microsupport estimation from (3) of Proposition \ref{mses} and conclude that
 its microsupport $\ms(F)$ is bounded by $f_\pi ( (df^*)^{-1} \ms(1_C) )$.
Because $\ms(1_C) = N^* C$, the slice of $(df^*)^{-1} \ms(1_C)$ with $\{ z = 0 \}$ is given by $(\RR^*)^2 \times \{0\}$
at $(x,y,0)$. 
Thus $\ms(F) \cap T^*_{(0,0,0)}\RR^3$ is contained in $N^* X_1$.
However, by base change $F |_{X_1} = (f|_C)_* ( 1_{S^1 \times \{0 \} \cup \{(1,0)\} \times \RR^1} )$ 
is not a locally constant sheaf on $X_1$ so $F \not \in \Sh_{\sS}(\RR^3)$.

\begin{figure}[h]
  \centering
  \includegraphics[scale = 0.12]{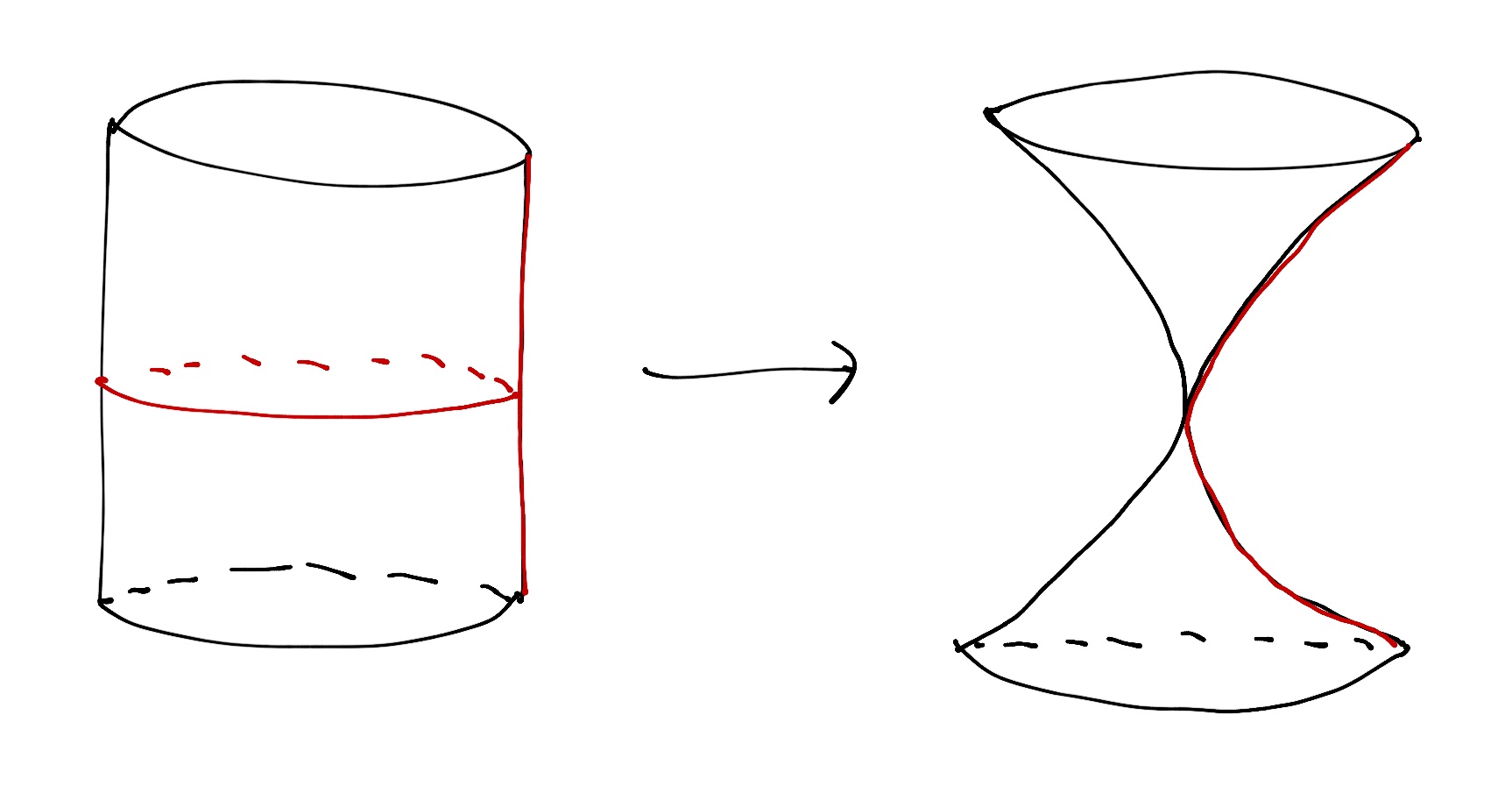}
  \caption{The picture exhibiting the fact that $F$ above is not locally constant on $X_1$.
  Stalks along from $0$ are given simply by $1_\cV = \Gamma(\{*\};\cV)$ but the stalk at $0$ is $\Gamma(S^1;\cV)$.}
\end{figure}

\end{example}

The above proposition is a corollary of the existence of \textit{inward cornerings} defined by applying the 
following lemma which is also proved in the same paper and we will use them for the main theorem of this paper as well.

\begin{lemma}[{\cite[Proposition 2.3]{Ganatra-Pardon-Shende3}}]\label{epnbd}
Fix any $1 \leq p \leq \infty$, and let $\mathscr{S}$ be a $C^p$ Whitney stratification of  $M$.
Fix a relatively compact $\mathscr{S}$-constructible set $Y$.
Let $\mathscr{S}_Y \coloneqq \{s | X_s \subseteq Y\}$ denote the collection of strata consisting of $Y$
and set $N^* \mathscr{S}_Y \coloneqq \cup_{s \in \mathscr{S}_Y} N^* X_s$ 
, which is closed in $T^* X$ by the Whitney condition.
Then there exists a decreasing family $Y^\epsilon$ of open neighborhoods of $Y$ such that as $\epsilon \rightarrow 0$,
\begin{enumerate}
\item $N^* Y^\epsilon$ becomes contained in arbitrary small conic neighborhood of $N^* Y$,
\item $N^* Y^\epsilon \cap N^* \mathscr{S} = \varnothing$ for $\epsilon > 0$.
\end{enumerate}
\end{lemma}

\begin{definition}\label{incor}
For a relative compact $\mathscr{S}$-constructible open set $U$,
an \textit{inward cornering} of $U$ is an open set of the form
$$ U^{-\epsilon} \coloneqq U \setminus \overline{(\partial U)^{\epsilon} }.$$
When $\epsilon > 0$ is small,
the inward cornering $U^{-\epsilon}$ is a codimension $0$ open submanifold 
whose closure $\overline{U^{-\epsilon}}$ is a compact manifold with corners.
The family $U^{-\epsilon}$ depends smoothly on $\epsilon$.
Its outward conormal $N^*_{\infty,out} U^{-\epsilon}$ remains disjoint from $N^*_\infty \mathscr{S}$ as $\epsilon$ changes, 
and converges to $N^*_\infty \mathscr{S}$ uniformly as $\epsilon \rightarrow 0$.  
\end{definition}

Combining with the comment on triangulations, we obtain a simple description of sheaves microsupported in $N^*_\infty \mathscr{S}$ 
for some $C^1$ Whitney triangulation $\mathscr{S}$.

\begin{proposition}[{\cite[Proposition 4.19]{Ganatra-Pardon-Shende3}}]\label{mc=cc}
Let $\mathscr{S}$ be a $C^1$ Whitney triangulation. 
Then, there is an equivalence 
\begin{align*}
\Sh_{N^*_\infty \mathscr{S}}(M) &= \mathscr{S} \dMod \\
1_{X_s} &\leftrightarrow 1_s
\end{align*}
where $1_s$ is the indicator which is defined by 
$$1_s(t) =
\begin{cases}
1, &t \leq s. \\
0, &\text{otherwise}.
\end{cases}
$$
In particular, the category  $\Sh_{N^*_\infty \mathscr{S}}(M)$ is compactly generated whose
compact objects $\Sh_{N^*_\infty \mathscr{S}}(M)^c$ are given by sheaves with compact support and perfect stalks.
\end{proposition}

\subsection{Isotropic microsupport}\label{ims}

We say a subset $\Lambda \subseteq S^* M$ is isotropic if it can be stratified by isotropic submanifolds.
A standard class of isotropic subsets are given by the conormal $N^*_\infty \sS$ of a stratification $\sS$ which we study 
in the last section.
Assume $M$ is real analytic and we recall that a general isotropic which satisfies a decent regularity condition
is bounded by isotropics of this form.

\begin{definition}
A subset $Z$ of $M$ is said to be subanalytic at $x$ if there exists open set $U \ni x$, 
compact manifolds $Y_j^i$ $(i = 1, 2, 1 \leq j \leq N)$ and morphisms
$f_j^i: Y_j^i \rightarrow M$ such that 
$$Z \cap U = U \cap \bigcup_{j=1}^N (f_j^1(Y_j^1) \setminus f_j^2(Y_j^2)).$$
We say $Z$ is subanalytic if $Z$ is subanalytic at $x$ for all $x \in M$.
\end{definition}

\begin{lemma}[{\cite[Corollary 8.3.22]{Kashiwara-Schapira1}}]
Let $\Lambda$ be a closed subanalytic isotropic subset of $S^* M$.
Then there exists a $C^\omega$ Whitney stratification $\mathscr{S}$ such that $\Lambda \subseteq N^* \mathscr{S}$.
\end{lemma}

Combining with the above proposition, we obtain a microlocal criterion for a sheaf $F$ with
subanalytic microsupport being constructible:

\begin{proposition}[{\cite[Theorem 8.4.2]{Kashiwara-Schapira1}}]
Let $F \in \Sh(M)$ and assume $\msif(F)$ is subanalytic. 
Then $F$ is constructible if and only if $\msif(F)$ is a singular isotropic. 
\end{proposition}

Another feature of subanalytic geometry is that relatively compact subanalytic sets form an o-minimal structure.
Thus, one can apply the result of \cite{Czapla} to refine a $C^p$ Whitney stratification to a Whitney triangulation, for $1 \leq p < \infty$.

\begin{lemma}\label{cbs}
Let $\Lambda$ be a subanalytic singular isotropic in $S^* M$.
Then there exists a $C^1$ Whitney triangulation $\mathscr{S}$ such that $\Lambda \subseteq N^* \mathscr{S}$.
\end{lemma}

Combining the above two results, we conclude:
\begin{theorem}\label{extri}
Let $F \in \Sh(M)$ and assume $\msif(F)$ is a subanalytic singular isotropic.
Then $F$ is $\mathscr{S}$-constructible for some $C^1$ Whitney triangulation $\mathscr{S}$. 
\end{theorem}

Collectively, sheaves with the same subanalytic isotropic microsupport form a category with nice finiteness properties.
Let $\Lambda$ be a subanalytic singular isotropic in $S^* M$.

\begin{proposition}
Let $\Lambda$ be a subanalytic singular isotropic in $S^* M$.
The category $\Sh_\Lambda(M)$ is compactly generated.
\end{proposition}

\begin{proof}
Fix a $C^1$ Whitney triangulation $\mathscr{S}$ such that $\Lambda \subseteq N^*_\infty \mathscr{S}$
by Lemma \ref{cbs}.
Recall that $\mathscr{S} \dMod = \Ind (\Perf \mathscr{S})$ is compactly generated.
For $F \in \Sh_\Lambda(M)$, there exists $F_i \in \Sh_{N^*_\infty \mathscr{S} }(M)$ such that 
${\iota_\Lambda}_* F = \varinjlim F_i$.
Thus $$F  = \iota_\Lambda^* {\iota_\Lambda}_* F 
= \iota_\Lambda^* {\iota_\Lambda}_* \varinjlim F_i = \varinjlim \iota_\Lambda^* F_i.$$
Now note that ${\iota_\Lambda}^* F_i$ is compact in $\Sh_\Lambda(M)$ since 
$\iota_\Lambda^* \dashv {\iota_\Lambda}_* \dashv \iota_\Lambda^!$
and the left adjoint of left joint preserves compact objects.
\end{proof}

\begin{corollary}\label{ps=>cp}
When $M$ is compact, if $F \in \Sh_\Lambda(M)$ has perfect stalks, then $F$ is in $\Sh_\Lambda(M)^c$.
\end{corollary}

\begin{proof}
Apply Lemma \ref{cbs} to include $\Sh_\Lambda(M) \subseteq \Sh_\sS(M)$ as in the last proposition and
conclude the statement by Proposition \ref{mc=cc}.
\end{proof}

Now recall from Proposition \ref{shissh}, the presheaf $\Sh$ in $\Cat$ of ($\cV$-valued) sheaves is itself a sheaf.
Since a set can be recovered from its intersections with an open cover, the same argument shows that the assignment 
$U \mapsto \Sh_\Lambda(U)$ forms a sheaf $\Sh_\Lambda$ in $\Cat$ as well.
In fact, since for an open inclusion $j: U \hookrightarrow M$, $j^* = j^!$ is a right adjoint,
$\Sh_\Lambda$ forms a sheaf in $\PrRst = (\PrLst)^{op}$.
By Proposition \ref{clmbrad} and Proposition \ref{cs=st},
passing to left adjoints turns $\Sh_\Lambda$ to a cosheaf in $\PrLcs$, 
and taking compact objects further turns it to a cosheaf in $\st$.

\begin{proposition}\label{shciscsh}
The precosheaf $\Sh_\Lambda^c: \Op_{M} \rightarrow \st$ is a cosheaf.
\end{proposition}

For an inclusion of subanalytic singular isotropics $\Lambda \subseteq \Lambda^\prime$, 
by picking a $C^1$ Whitney triangulation $\mathscr{S}$ such that $\Lambda^\prime \subseteq N^*_\infty \mathscr{S}$, 
a similar consideration as above shows that the inclusion $\Sh_\Lambda(M) \hookrightarrow \Sh_{\Lambda^\prime}(M)$
has both a left and a right adjoint. Thus,

\begin{proposition}
Passing to left adjoint, the inclusion $\Sh_\Lambda(M) \hookrightarrow \Sh_{\Lambda^\prime}(M)$
induces a canonical functor $\Sh_{\Lambda^\prime}(M)^c \twoheadrightarrow \Sh_\Lambda(M)^c$ between compact objects.
\end{proposition}

Let $\Lambda$ be a singular isotropic, $(x,\xi) \in \Lambda$, and consider the restriction of the microstalk functor 
to the compactly generated categories $\mu_{(x,\xi)}:\Sh_\Lambda(M) \rightarrow \cV$.
By applying its left adjoint to the generator $1 \in \cV$, we see that it is
tautologically corepresented by the compact object $\mu_{(x,\xi)}^L (1) \in \Sh_\Lambda(M)^c$.
Furthermore, when there is an inclusion $\Lambda \subseteq \Lambda^\prime$ and $(x,\xi) \in \Lambda^\prime$,
the corepresentative $\mu_{(x,\xi)}^L (1) \in \Sh_{\Lambda^\prime}(M)^c$ 
is sent under $\Sh_{\Lambda^\prime}(M)^c \twoheadrightarrow \Sh_\Lambda(M)^c$ to a similar corepresentative 
in $\Sh_\Lambda(M)^c$ and,
they are tautologically sent to the zero object when  $(x,\xi)$ is a smooth point in $ \Lambda^\prime \setminus \Lambda$.
By Proposition \ref{msbms}, the converse is also true:

\begin{proposition}[Theorem 4.13 of \cite{Ganatra-Pardon-Shende3}]\label{wdstopr}
Let $\Lambda \subseteq \Lambda^\prime$ be subanalytic isotropics and let
$\sD^\mu_{\Lambda^\prime,\Lambda}(T^* M)$ denote the fiber of the canonical functor
$\Sh_{\Lambda^\prime}(M)^c \twoheadrightarrow \Sh_{\Lambda}(M)^c$.
Then $\sD^\mu_{\Lambda^\prime,\Lambda}(T^* M)$ 
is generated by the corepresentatives of the microstalk functors $\mu_{(x,\xi)}$ for smooth Legendrian points 
$(x,\xi) \in \Lambda^\prime \setminus \Lambda$.
\end{proposition}

\section{Isotopies of sheaves}\label{iso_sh}

Let $(X,\omega,Z)$ be a Liouville manifold and $\alpha \coloneqq \iota_Z \omega$ be the Liouville form.
Consider an isotopy of Lagrangian submanifolds conic at infinity $L_t$, $t \in [0,1]$.
One says that the isotopy $L_t$ is \textit{positive} if $\alpha(\partial_t \partial_\infty L_t) \geq 0$.
Standard Floer theory implies that there is a \textit{continuation element} $c(L_t) \in \HF^*(L_0,L_1)$
if $L_0$ and $L_1$ intersect transversally.
For any triple $(K_0,K_1,K_2)$ of transversally intersected Lagrangians,
there exists also a multiplication map $\mu: \HF^* (K_0,K_1) \otimes \HF^* (K_1,K_2) \rightarrow \HF^* (K_0,K_2)$.
Thus, for suitable $K$'s, 
multiplying  $c(L_t)$ induces a map $\HF^*(L_1,K) \rightarrow \HF^*(L_0,K)$ which is usually referred as the 
\textit{continuation map} and is one key ingredient for defining the wrapped Floer category.
See for example \cite[Section 3.3]{Ganatra-Pardon-Shende1} for details.

\subsection{Continuation maps}
We recall here the sheaf-theoretical continuation maps studied in \cite{Guillermou-Kashiwara-Schapira}.
A dual construction can be found in \cite{Tamarkin1} and \cite{Guillermou-Schapira}.
In the sheaf-theoretical setting, objects corresponding to the continuation elements are simply morphisms/maps between sheaves. 
Thus, we simply use the term \textit{continuation maps} to refer both the morphisms between sheaves 
and the induced maps on the $\Hom$'s.
We denote by $(t,\tau)$ the coordinate of $T^* \RR$ 
and by $T^*_\leq \RR = \{ \tau \leq 0 \}$ the set of non-positive covectors.

\begin{lemma}\label{ps=nc}
Let $M$ be a manifold, $[-\infty, \infty]$ be the compactification of $\RR$ at the two infinities, 
$p: M \times \RR \rightarrow M$ be the projection,
$j: M \times \RR \hookrightarrow M \times [-\infty,\infty]$ be the open interior,
and $i_{\pm}: M \times \{ \pm \infty \} \hookrightarrow M \times [-\infty, \infty]$ 
be the closed inclusion at the positive/negative infinity. 
Then for a sheaf $F \in \Sh_{T^* M \times T^*_\leq \RR }(M \times \RR)$,
there are isomorphisms $p_* F = i_-^* j_* F$ and $p_! F[1] = i_+^* j_* F$ identifying
the two pushforwards as nearby cycles at the infinities.
\end{lemma}
 
\begin{proof}

We first prove the case when $\supp (F) \subseteq M \times [-C,C]$ for some $C \in \Rp$.
In this case, $i_-^* j_* F = i_+^* j_* F = 0$ and $p_* F = p_! F$ since $p$ is proper on $\supp(F)$.
Let $x \in M$ be a point. 
Base change implies that $(p_* F)_x = \Gamma( \{ x \} \times \RR ; F |_{ \{ x \} \times \RR})$.
Apply the microsupport estimation $\ms(f^* F) \subseteq f^\# (\ms(F))$  
from (2) of Proposition \ref{ncmses} to the inclusion of the slice at $x$, 
and we obtain $\ms( F |_{ \{ x \} \times \RR}) \subseteq T^*_\leq \RR$
so we reduce to the case $M = \{ * \}$.
In this case, consider the family of open sets $\{ (-\infty,t) \}_{t \in \RR}$.
The noncharacteristic deformation lemma, Lemma \ref{ncdef}, implies that
$\Gamma( \RR;F) \xrightarrow{\sim} \Gamma( (-\infty,t); F)$ for all $t \in \RR$.
Since $\supp(F)$ is compact, the latter vanishes for $t << 0$.

Now for the general case, we first notice there are canonical morphisms $p_* F \rightarrow i_-^* j_* F$
and $i_+^* j_* F \rightarrow p_! F[1]$ functorial on $F$:
Let $j_- : M \times \RR \hookrightarrow M \times [-\infty,\infty)$ denote the open embedding  compactifying the negative end.
For any $G \in \Sh(M \times [-\infty,\infty) )$, there is a fiber sequence
$$ {j_-}_! j_-^* G  \rightarrow G \rightarrow {i_-}_! i_-^* G.$$
Set $G = {j_-}_* F$ and recall that $j_-^*  {j_-}_* = \id$, we obtain the fiber sequence
$$ {j_-}_! F \rightarrow {j_-}_* F \rightarrow {i_-}_! i_-^* {j_-}_* F.$$
Let $p_-: M \times [-\infty, \infty) \rightarrow M$ denote the projection (and similarly for $p_+$).
The canonical morphism $p_* F \rightarrow i_-^* j_* F$ is obtained 
by applying ${p_-}_*$ to the above fiber sequence.
The morphism $i_+^* j_* F \rightarrow p_! F[1]$ can be obtained similarly.

Recall that there is fiber sequence
$$ F_{ M \times (-\infty,0] } \rightarrow F \rightarrow F_{ M \times (0,\infty)}.$$
(5) of Proposition \ref{mses} implies that both $ F_{ M \times (-\infty,0] }$ 
and $F_{ M \times (0,\infty)}$ are contained in $T^* M \times T^*_\leq \RR$.
So it is sufficient to prove the cases when $\supp(F) \subseteq M \times (-\infty,C]$ and
when $\supp(F) \subseteq M \times [-C,\infty) $ for some $C \in \Rp$.

We first prove the cases in which the objects vanish:
Assume $\supp(F) \subseteq [-C,\infty)$. 
We claim $p_* F = 0 = i_-^* j_* F$.
One computes
$$ p_* F = p_* F_{[-C,\infty)} = p_* \lmi{n \rightarrow \infty} F_{ [-C,n] } = \lmi{n \rightarrow \infty} p_* F_{ [-C,n] } = 0$$
by the case when $\supp(F) \subseteq M \times [-C,C]$ for some $C \in \RR$.
Similarly, by considering the colimit $$F_{(-\infty,C]} = \clmi{n \rightarrow \infty} \ F_{(-n,C]},$$
one can conclude $p_! F = 0 = i_+^* j_* F$ when $\supp(F) \subseteq (-\infty,C]$.

Now assume $\supp(F) \subseteq (-\infty,C]$ and we claim $p_* F = i_-^* j_* F$.
Consider again the fiber sequence
$$ {j_-}_! F \rightarrow {j_-}_* F \rightarrow {i_-}_! i_-^* {j_-}_* F.$$
Apply ${p_-}_*$ and notice that ${p_-}_! = {p_-}_*$ for these sheaves because of the compact support assumption.
Thus, we obtain the fiber sequence
$$ p_! F \rightarrow p_* F \rightarrow i_-^* {j_-}_* F.$$
Since $p_! F = 0$ by the previous case, $p_* F = i_-^* {j_-}_* F$.
The other isomorphism can be obtained similarly.
\end{proof}

To define the continuation map, we need a prototype version of Theorem \ref{w=ad}.
For the rest of the section, we use $I$ to denote an open interval and $(t,\tau)$ to denote
the coordinate of its cotangent bundle $T^* I$. 

\begin{proposition}[{\cite[Proposition 4.8]{Guillermou-Kashiwara-Schapira}}]\label{pc}
Let $\iota_*: \Sh_{T^* M \times T^*_\leq I}(M \times I) \hookrightarrow \Sh(M \times I)$
denote the tautological inclusion. 
Then there exist left and right adjoints $\iota^* \dashv \iota_* \dashv \iota^!$ which are given by convolution
$ \iota^* F = 1_{ \{ t^\prime > t \} }[1] \circ F$ and its right adjoint
$\iota^!F = \sHom^\circ ( 1_{ \{ t^\prime > t \} }[1], F)$.
Here we denote by $(t,t^\prime)$ the coordinate of $I^2$.
\end{proposition}

\begin{proof}
Assume for simplicity $I = (0,1)$.
We first show that $ \ms( 1_{ \{t^\prime > t \} } \circ F) \subseteq T^* M \times T^*_\leq I$.
Let $\pi_1, \pi_2: M \times I \times I \rightarrow M \times I$ denote the projection $\pi_1(x, t, t^\prime) = (x,t)$
and $\pi_2(x,t,t^\prime) = (x,t^\prime)$.
Then $ 1_{ \{t^\prime > t \} } \circ F = {\pi_2}_! [(\pi_1^* F)_{M \times \{ t^\prime > t \}  }] $.
In order to estimate the effect of ${\pi_2}_!$ on the microsupport, 
we need the map $\pi_2$ to be proper on the support of the sheaf. 
Thus, let $j:M \times I \times I \hookrightarrow M \times (-1,1) \times I$
denote the open inclusion extending $I \hookrightarrow (-1,1)$, 
$\pi_2^\prime: M \times (-1,1) \times I \rightarrow M \times I$ denote the projection to the first and the third components,
and factorize $1_{ \{t^\prime > t \} } \circ F$ to
${\pi^\prime_2}_* j_! [(\pi_1^* F)_{M \times \{ t^\prime > t \}  }]$.
Before taking ${\pi^\prime_2}_*$, one observes that, 
by (4) and (5) of Proposition \ref{mses} and (1) of Proposition \ref{ncmses}, 
none of the operations introduces non-zero covectors on the second $I$-component to the microsupport 
except when taking $(-)_{M \times \{ t^\prime > t \}  }$, 
covectors of the form $(0,\sigma,-\sigma)$ for $\sigma \in \Rp$ might be added to the cotangent fibers over 
the boundary $\{ t^\prime = t \}$.
Thus 
\begin{align*}
\ms( 1_{ \{t^\prime > t \} } \circ F )
&= \ms({\pi^\prime_2}_* j_! [(\pi_1^* F)_{M \times \{ t^\prime > t \}  }]) \\
&\subseteq  (\pi^\prime_2)_\pi \left(\ms( j_! [(\pi_1^* F)_{M \times \{ t^\prime > t \}  }]) 
\cap T^* M \times 0_{(-1,1)} \times T^* I \right) \\
&\subseteq T^* M \times T^*_\leq I.
 \end{align*}

For the right adjoint $\iota^!$, we note that 
\begin{align*}
\lmi{r \rightarrow \infty} ( 1_{ \{t - r \leq t^\prime \leq t \} } \circ F )
&= \lmi{r \rightarrow \infty} {\pi_2}_! [(\pi_1^* F)_{M \times \{ t - r \leq t^\prime \leq t \}  }]\\
&= \lmi{r \rightarrow \infty} {\pi_2}_* [(\pi_1^* F)_{M \times \{ t - r \leq t^\prime \leq t \}  }] \\
&={\pi_2}_*  \lmi{r \rightarrow \infty}  [(\pi_1^* F)_{M \times \{ t - r \leq t^\prime \leq t \}  }] \\
&={\pi_2}_*  [(\pi_1^* F)_{M \times \{ t^\prime \leq t \}  }].
\end{align*}
Thus one can argue as the left adjoint case.
(Note the last term is different from $1_{ \{t^\prime \leq t \} } \circ F $ in general
since limits do not commute with convolution.)

In sum, we've shown that there are functors 
$$1_{ \{t^\prime > t \} }[1] \circ {(-)}, 
\ \lmi{r \rightarrow \infty} \left( 1_{ \{t - r \leq t^\prime \leq t \}  } \circ {(-)} \right):
\Sh(M \times I) \rightarrow \Sh_{T^* M \times T^*_\leq I}(M \times I).$$
In order to show that these are indeed the desired adjoints, it is sufficient to show, for example, 
that the canonical morphism $1_{\Delta_I} \rightarrow 1_{ \{t^\prime > t \} }[1] $
becomes an isomorphism after convolving with sheaves 
in $\Sh_{T^* M \times T^*_\leq I}(M \times I)$.
Recall that convolving with $1_{\Delta_I}$ is the same as the identity functor.

Consider the fiber sequence 
$$ 1_{ \{t^\prime > t\}} \rightarrow 1_{ \{t^\prime \geq t\}} \rightarrow 1_{\Delta_I}.$$
We have similarly $ 1_{ \{t^\prime \geq t\}} \circ F = {\pi_2}_! (\pi_1^* F)_{M \times \{ t^\prime \geq t \} }$ 
and a similar microsupport estimation implies that, before applying ${\pi_2}_!$,
$\ms\left( (\pi_1^* F)_{M \times \{ t^\prime \geq t \} } \right) \subseteq T^* M \times T^*_\leq I \times T^* I$.
Thus, the last lemma \ref{ps=nc} implies ${\pi_2}_! [ (\pi_1^* F)_{M \times \{ t^\prime \geq t \} } ]$ is the nearby circle
of $ (\pi_1^* F)_{M \times \{ t^\prime \geq t \} }$ at $\infty$ along the first $I$-direction and it is $0$.
Thus $F = 1_{\Delta_I} \circ F \xrightarrow{\sim}  1_{ \{t^\prime > t \} }[1] \circ F$.
A similar argument shows 
$ \lmi{r \rightarrow \infty} ( 1_{ \{t - r \leq t^\prime \leq t \} } \circ F ) \xrightarrow{\sim} F$
for $F$ with the same microsupport condition.

Finally, we notice that $1_{ \{t - r \leq t^\prime \leq t \} } \circ (-)$ 
is the inverse of $1_{ \{t + r > t^\prime > t \} }[1] \circ (-)$.
Thus the functor $1_{ \{t - r \leq t^\prime \leq t \} } \circ (-)$ is equivalent to 
$\sHom^\circ(1_{ \{t + r > t^\prime > t \} }[1], - )$
and
\begin{align*}
&\lmi{r \rightarrow \infty} ( 1_{ \{t - r \leq t^\prime \leq t \} } \circ F )
=  \lmi{r \rightarrow \infty} \sHom^\circ (1_{ \{t + r > t^\prime > t \} }[1], F) \\
&= \sHom^\circ (\clmi{r \rightarrow \infty} \ 1_{ \{t + r > t^\prime > t \} }[1], F) 
= \sHom^\circ ( 1_{ \{ t^\prime > t \} }[1], F).
\end{align*}  

\end{proof}


Now let  $F \in \Sh_{T^* M \times T^*_\leq I}(M \times I)$ and, 
by the preceding lemma, $F \xrightarrow{\sim}  1_{ \{t^\prime > t \} }[1] \circ_I F$.
Let $a \in I$ and let $i_a: M \hookrightarrow M \times I$ denote the slice at $a$.
Applying $i_a^*$ results the isomorphism $i_a^* F \xrightarrow{\sim} 1_{ (-\infty,a)}[1] \circ_I F$.
Recall that for $a \leq b$, there is a canonical morphism $1_{ (-\infty,a)}[1] \rightarrow 1_{ (-\infty,b)}[1]$
induced by the open inclusion $(-\infty,a) \hookrightarrow (-\infty,b)$.

\begin{definition}
For $F \in \Sh_{T^* M \times T^*_\leq I }(M \times I)$ and $a \leq b$
Set $F_x = i_x^* F$ for $x \in I$.
We define the continuation map
$c(F,a,b): F_a \rightarrow F_b$ to be the (homotopically unique)
morphism $c$ that makes the following diagram commute:

$$
\begin{tikzpicture}

\node at (0,2) {$F_a$};
\node at (5,2) {$F_b$};
\node at (0,0) {$1_{ (-\infty,a)}[1] \circ F$};
\node at (5,0) {$1_{ (-\infty,b)}[1] \circ F$};

\draw [->, thick] (0.5,2) -- (4.6,2) node [midway, above] {$c$};
\draw [->, thick] (1.4,0) -- (3.6,0) node [midway, above] {$ $};

\draw [double equal sign distance, thick] (0,1.6) -- (0,0.4) node  [midway, right] {$ $};
\draw [double equal sign distance, thick] (5,1.6) -- (5,0.4) node [midway, right] {$ $};

\end{tikzpicture}
 $$
 
\end{definition}

The continuation maps inherit various properties from $1_{(-\infty,a)}$.
For example, they compose in the sense that $$c(F,a_2,a_3) \circ c(F,a_1,a_2) = c(F,a_1,a_3)$$
since the canonical maps $1_{(-\infty,a_1)} \rightarrow 1_{(-\infty,a_2)} \rightarrow 1_{(-\infty,a_3)}$
compose to $1_{(-\infty,a_1)}  \rightarrow 1_{(-\infty,a_3)}$.
Let $p_{[a,b]}: N \times [a,b] \rightarrow N$ denote the projection.
If $F|_{N \times [a,b]} =  p_{[a,b]}^* G$ is a pullback from $N$ for some $G \in \Sh(N)$,
one can identify $F_a = G = F_b$ through the canonical map $F \rightarrow {i_a}_* i_a^* F$.
In this case, the continuation map $c(F,a,b)$ is equivalent to this identification $F_a = F_b$.
Note because convolution $\circ$ is compatible with colimits, we have the following corollary.
\begin{corollary}
For $F \in \Sh_{T^* N \times T^*_\leq I}(N \times I) $, 
the canonical map $\clmi{r < t} \ F_r \rightarrow F_t$ is an isomorphism.
\end{corollary}
The dual statement for limits is false in general.
\begin{example}
Let $N = \{ * \}$ and take $F = 1_{(-\infty,0]}$. Then $F_t = 1_\cV$ when $t \leq 0$ and $0$ otherwise.
Thus $F_0 = \clmi{r < 0} \ F_r$ but $F_0 \neq \lmi{s > 0} \ F_s$.
\end{example}

However, we will consider the noncharacteristic situation.
\begin{definition}[{\cite{Nadler-Shende}}]
Let $B$ be a manifold. For $F \in \Sh(M \times B)$, 
we say that $F$ is \textit{$B$-noncharacteristic} if the inclusion $i_b: M \times \{b\} \hookrightarrow M \times B$ 
is noncharacteristic for $F$ for all $b \in B$.
Equivalently, $F$ is $B$-noncharacteristic if $\ms(F) \cap ( 0_M \times T^* B) \subseteq 0_{M \times B} $.
\end{definition}

\begin{lemma}\label{ncpro}
Let $F \in \Sh(M \times I)$ be $I$-noncharacteristic. Then, 
\begin{enumerate}
\item The natural morphism $i_t^* F[-1] \rightarrow i_t^! F$ is an isomorphism for $t \in I$.
\item If $G \in \Sh(M \times I)$ such that $\sHom(G,F)$ is $I$-noncharacteristic, 
then $i_t^* \sHom(G,F) = \sHom(i_t^* G, i_t^* F)$.
\item Denote by $q: M \times I \rightarrow I$ the projection. 
If $q$ is proper on $F$, then $q_* F$ is a constant sheaf. 
Moreover, if $F \in \Sh_{T^* M \times T^*_\leq I}(M \times I)$, 
the continuation maps of $F$ are sent to isomorphisms under $q_*$.
\end{enumerate}
\end{lemma}

\begin{proof}

For (i), the equivalence $i_t^* F[1] = i_t^! F$ follows directly from (4) of Proposition \ref{mses}
and the observation that since $q \circ i_t = \id$, one has
$ i_t^! 1_{M \times I} = i_t^! \circ q^! 1_{M}[1] = 1_{M}[1].$

For (ii), apply (i) and use the fact that $i_t^! \sHom(G,F)[-1]  = \sHom(i_t^* G, i_t^! F) [-1]$.

For (iii), (3) of Proposition \ref{mses} implies that 
$\ms(q_* F) \subseteq q_\tau (\ms(F) \cap 0_M \times T^* I) \subseteq 0_{I}$ 
so $q_* F$ is a constant sheaf by Example \ref{loceg} and the fact that $\RR$ is contractible.
For the statement of continuation maps, recall that they are given by
$ 1_{(-\infty,s)} \circ F \rightarrow 1_{(-\infty,t)} \circ F$
for $s \leq t$. Since $q_* F = q_! F$, they are sent to
$\Gamma_c((-\infty,s); q_* F) \rightarrow \Gamma_c((-\infty,t); q_*F)$ by $q_*$
which are isomorphisms.
\end{proof}

\begin{corollary}\label{nclp}
Let $F \in \Sh_{T^* N \times T^*_\leq I}(N \times I)$.
If $F$ is $I$-noncharacteristic, 
then the canonical map $F_t \rightarrow \lmi{s > t} \ F_s$ is an isomorphism.
\end{corollary}

\begin{proof}
We will use $q_{ij}$ to denote the projection from $N \times I \times I$ to the corresponding components
and $p_i$ the projection from $N \times I$.
By the above Lemma \ref{ncpro}, $F_t =  i_t^! F [1]$.
Apply Proposition \ref{pc} and compute by base change and (4) of Proposition \ref{cp6f}, we see that
\begin{align*}
F_t
&= i_t^! F [1] = i_t^! \sHom^\circ ( 1_{ \{ s^\prime > s \} }[1], F) [1] 
= i_t^! {q_{12}}_* \sHom(1_{ N \times \{ s^\prime > s \} }  ,q_{13}^! F) \\
&= {p_1}_* (i_t \times \id_I)^! \sHom(1_{ N \times \{ s^\prime > s \} }  ,q_{13}^! F) 
= {p_1}_* \sHom(1_{N \times (t,\infty)}, F).
\end{align*}
That is, for $t \leq s$, the continuation map also corresponds to the map
$$ {p_1}_* \sHom(1_{N \times (t,\infty)} [1], F) \rightarrow {p_1}_* \sHom(1_{N \times (s,\infty)} [1], F)$$
and, since $*$-push commutes with limits and $\sHom(-,F)$ turns limits to colimits, we have
$$\lmi{s > t} \ F_s 
= \lmi{s > t} \ {p_1}_* \sHom(1_{N \times (s,\infty)} [1], F) = {p_1}_* \sHom( 1_{N \times (t,\infty)} [1], F) = F_t.
$$
\end{proof}


We conclude this section with a homotopical invariant property of the continuation map in the following setting.
Let $I$ and $J$ be open intervals and let $(t,\tau)$ and $(s,\sigma)$ be the corresponding coordinates
of their cotangent bundles.
Let $G \in \Sh(M \times I \times J)$ be a sheaf such that $\ms(G) \subseteq \{ \tau \leq 0 \}$.
For any $x \in I$, we use $G_{t = x} \coloneqq G|_{M \times \{x\} \times J }$ to denote the restriction 
and similarly for $G_{s = y}$, $y \in J$.
Note by (2) of Proposition \ref{ncmses}, the same condition $\ms(G_{s = y}) \subseteq \{ \tau \leq 0 \}$ holds.
Assume further that there exists $a \leq b$ in $I$ such that $\ms(G_{t = a})$, $\ms(G_{t = b}) \subseteq T^* M \times 0_J$.
By Lemma \ref{constbms}, this implies that there exist $F_a$, $F_b \in \Sh(M)$ such that 
$G_{t = a} = p_s^* F_a$ and $G_{t = b} = p_s^* F_b$
where we use $p_s: M \times J \rightarrow M$ to denote the projection.
Note that, for each $y \in J$, the restriction $G_{s = y}$ induces a continuation map 
$c(G,y,a,b): F_a \rightarrow F_b$.

$$
\begin{tikzpicture}

\draw [thick] (0,0) rectangle (5,3);

\draw [thick] (1,0) -- (1,3);
\draw [thick] (4,0) -- (4,3);

\node at (0.6,2.5) {$F_a$};
\node at (3.6,2.5) {$F_b$};

\draw [->, thick] (1.1,1.5) -- (3.9,1.5) node [midway, above] {$c(G,y,a,b)$};
\draw [->, thick] (1.1,0.5) -- (3.9,0.5) node [midway, above] {$c(G,y^\prime,a,b)$};

\node at (1,3.25) {$t = a$};
\node at (4,3.25) {$t = b$};
\node at (-0.7,0.5) {$s = y^\prime$};
\node at (-0.7,1.5) {$s = y$};

\end{tikzpicture}
 $$

\begin{proposition}\label{homotopy_independence_of_continuation_maps}
The morphism $c(G,y,a,b)$ is independent of $y \in J$.
\end{proposition}

\begin{proof}
Since $\ms(G) \subseteq \{ \tau \leq 0 \}$, a family version of Proposition \ref{pc} implies
$$G = 1_{\Delta_I \times J} \circ|_J G 
\xrightarrow{\sim} 1_{ \{s^\prime > s \}  \times J} [1] \circ|_J G$$ 
is an isomorphism where $\circ|_J$ is the $J$-parametrized convolution.
In particular, there is an isomorphism $G_{t = a} \xrightarrow{\sim} 1_{ (-\infty,a) \times J} [1] \circ|_J G$
and thus a ($J$-parametrized) continuation map
$$c_J(G,a,b): G_{t = a} \rightarrow G_{t = b}.$$

For $y \in J$, let $i_y: M \rightarrow M \times J$ denote the inclusion of the slice at $y$.
By Proposition \ref{cp6f}, there is equivalence $i_y^* (K \circ|_J G) = K|_{s = y} \circ G_{s = y}$ for $K \in \Sh(I \times J)$.
This  implies that $c_J(G,a,b)$ restricts to $i_y^* c_J(G,a,b) = c(G_{s = y},a,b) \eqqcolon c(G,y,a,b)$.
Hence, the ${i_y}^* \dashv {i_y}_*$ adjunction induces a commuting diagram,

$$
\begin{tikzpicture}

\node at (0,2) {$ G_{t=a}$};
\node at (5,2) {$ G_{t=b}$};
\node at (0,0) {$ {i_y}_* {i_y}^* G_{t =a}$};
\node at (5,0) {$ {i_y}_* {i_y}^* G_{t =b} $};

\draw [->, thick] (0.6,2) -- (4.4,2) node [midway, above] {$c_J(G,a,b)$};
\draw [->, thick] (1,0) -- (4,0) node [midway, above] {${i_y}_* c(G,y,a,b)$};


\draw [->, thick] (0,1.7) -- (0,0.3) node [midway, left] {$$};
\draw [->, thick] (5,1.7) -- (5,0.3) node [midway, left] {$$};

\end{tikzpicture}
$$

which is equivalent to 

$$
\begin{tikzpicture}

\node at (0,2) {$ p_s^* F_a$};
\node at (5,2) {$ p_s^* F_b$};
\node at (0,0) {${i_y}_*  F_a$};
\node at (5,0) {${i_y}_*  F_b $};

\draw [->, thick] (0.6,2) -- (4.4,2) node [midway, above] {$c_J(G,a,b)$};
\draw [->, thick] (0.6,0) -- (4.4,0) node [midway, above] {${i_y}_* c(G,y,a,b)$};


\draw [->, thick] (0,1.7) -- (0,0.3) node [midway, left] {$$};
\draw [->, thick] (5,1.7) -- (5,0.3) node [midway, left] {$$};

\end{tikzpicture}.
$$

Since $J$ is contractible, the horizontal arrows become isomorphism after applying ${p_s}_*$.

$$
\begin{tikzpicture}

\node at (0,2) {$ F_a$};
\node at (5,2) {$ F_b$};
\node at (0,0) {$ F_a$};
\node at (5,0) {$ F_b $};

\draw [->, thick] (0.4,2) -- (4.6,2) node [midway, above] {${p_s}_* c_J(G,a,b)$};
\draw [->, thick] (0.4,0) -- (4.6,0) node [midway, above] {$c(G,y,a,b)$};


\draw [double equal sign distance, thick] (0,1.7) -- (0,0.3) node [midway, left] {$ $};
\draw [double equal sign distance, thick] (5,1.7) -- (5,0.3) node [midway, left] {$ $};

\end{tikzpicture}.
$$
That is, the continuation map $c(G,y,a,b)$ is equivalent to ${p_s}_* c_J(G,a,b)$ for all $y \in J$.

\end{proof}

\begin{remark}
One can see from the proof that the continuation maps enjoy higher homotopical independence.
\end{remark}

\subsection{Sheaf-theoretical wrappings}\label{Sheaf-theoretical wrappings}

We specialize to the cases of $I$-family sheaf 
which come from the \textit{Guillermou-Kashiwara-Schapira sheaf quantization} in this section.
Recall that when $M$ is a smooth manifold, its cotangent bundle admits a canonical symplectic structure $(T^* M, d \alpha)$.
The Liouville form $\alpha$ is compatible with the $\Rp$-action which freely acts on $\dT^* M$.
Thus, there is an induced contact structure on the cosphere bundle $S^* M$.
It can be realized as a contact hypersurface of $\dT^* M$ by picking a Riemannian metric.
There is a dictionary between homogeneous symplectic geometry of $\dT^* M$ and contact geometry of $S^* M$.
Thus, we will use them interchangeably when one language is more convenient. 
See subsection \ref{hsnct} of the Preliminary for a more detailed discussion.

\begin{definition}
Let $M$, $B$ be manifolds and $I$ be an open interval containing $0$. 
We say a $C^\infty$ map $\Phi: S^* M \times I \times B \rightarrow S^* M$ is a \textit{$B$-family of contact isotopies}
if for each $(t,b) \in I \times B$, 
the map $\varphi_{t,b} \coloneqq \Phi(-,t,b)$ is a contactomorphism
and $\varphi_{0,b} = \id_{S^* M}$ for all $b \in B$.
\end{definition}

As remarked above, a $B$-family of contact isotopies $\Phi$ corresponds to a $B$-family of homogeneous symplectic isotopies
(of degree $1$), which we abuse the notation and denote it by $\Phi$ as well.
For fixed $b \in B$, we let $V_{\Phi_b}$ denote the vector field generated by $\varphi_{t,b}$.
Since $\varphi_{t,b}$ is homogeneous, $V_{\Phi_b}$ is a Hamiltonian vector field 
with $\alpha(V_{\Phi_b})$ being its Hamiltonian.
The latter is the function which evaluates to $\alpha_{\varphi_{t,b}(x,\xi)} ( \frac{\partial}{\partial t} \varphi_{(t,b)}(x,\xi) )$ 
at $\varphi_{t,b}(x,\xi)$. 

\begin{proposition}
For each $B$-family of homogeneous symplectic isotopies $\Phi$,
there is a unique conic Lagrangian submanifold $\Lambda_\Phi$ in $\dT^* (M \times M) \times T^* I \times T^* B$ 
which is determined by the equation $T^*_{t,b} (I \times B) \circ \Lambda_\Phi$ = $\Lambda_{\varphi_{t,b}}$ where the later is 
$\{ \left(x,-\xi,\phi_{t,b}(x,\xi) \right) | (x,\xi) \in \dT^* M \}$, the twisted graph of $\varphi_{t,b}$.
More precisely, it is given by the formula 
\begin{equation}\label{asso_lag}
\Lambda_\Phi = \left\{ \left(x, -\xi, \varphi_{t,b}(x,\xi), t, - \alpha(V_{\Phi_b})(\varphi_{t,b}(x,\xi)),
b, - \alpha_{\varphi_{t,b}(x,\xi)} \circ d (\Phi \circ i_{x,\xi,t})_b(\cdot) \right) \right\}
\end{equation}
where the parameters run through $(x,\xi) \in \dT^* M$, $t \in I$, $b \in B$, 
and the map $i_{x,\xi,t}$ is the inclusion of B as the $(x,\xi,t)$-slice.
We use the same notation $\Lambda_\Phi$ to denote its projection to $S^* (M \times M \times I \times B)$
which is a Legendrian submanifold.
\end{proposition}

The following theorem of Guillermou, Kashiwara, and Schapira
is a categorification of the more classical statements of quantization which usually have operators as the quantized objects.
The proof given there is the non-family case.
Since the (global) existence is proved by using the uniqueness property to glue local existence 
and the local picture depends smoothly on the family $J^n$,  
the same proof holds for the family version with minor modification \cite[Remark 3.9.]{Guillermou-Kashiwara-Schapira}.

\begin{theorem}[{\cite[Proposition 3.2]{Guillermou-Kashiwara-Schapira}}]
Let $M$ be a manifold. 
For a $J^n$-family contact isotopies $\Phi:  S^* M \times I \times J^n \rightarrow S^* M$ where $J$ is an open interval,
there exists a unique sheaf kernel
$K(\Phi) \in \Sh(M \times M \times I \times J^n)$ such that 
$\msif(K(\Phi)) \subseteq \Lambda_\Phi$ 
and $K(\Phi) |_{t = 0} = 1_{\Delta_M \times J^n}$.
Moreover, $\msif(K(\Phi)) = \Lambda_\Phi$ is simple along $\Lambda_\Phi$,
both projections $\supp(K) \rightarrow M \times I \times J^n$ are proper, 
and the composition is compatible with convolution in the sense that 
\begin{enumerate}
\item $K(\Psi \circ \Phi) = K(\Psi) \circ|_{I \times J^n} K(\Phi)$,
\item $K(\Phi^{-1}) \circ|_{I \times J^n} K(\Phi) = K(\Phi) \circ|_{I \times J^n} K(\Phi^{-1}) = 1_{\Delta_M \times I \times J^n}$.
\end{enumerate}
Here $\Phi^{-1}$ is the $J^n$-family of isotopies given by $\Phi^{-1}(-,t,b) \coloneqq \phi_{t,b}^{-1} $.
\end{theorem} 

\begin{remark}
The equality $\msif(K(\Phi)) = \Lambda_\Phi$ as well as a few other properties of $K(\Phi)$ followed by the uniqueness
is explained in \cite{Guillermou1}.
\end{remark}

We refer the above process of obtaining the sheaf kernel $K(\Phi)$ from a contact isotopy $\Phi$
or a $J^n$-family of contact isotopies $\Phi$, for $n > 1$,
as the \textit{Guillermou-Kashiwara-Schapira sheaf quantization} or GKS sheaf quantization in short.

\begin{example}\label{skreeb}
The construction in \cite[Example 3.10, Example 3.11]{Guillermou-Kashiwara-Schapira} works more generally:
Consider a manifold $M$ and take a Riemannian metric $g$ so that the corresponding injectivity radius is positive.
Denote by $H$ the homogeneous Hamiltonian $H(x,\xi) \coloneqq \sqrt{g_x(\xi,\xi)}$, $(x,\xi) \in \dT^* M$
and $\Phi$ the corresponding positive isotopy.
For small $-1<< s <0$, denote by $Z_s \coloneqq \{ (x,y) \in M \times M | d(x,y) \leq \abs{s}\}$ 
the closed subset of the pairs of points with distance less than $s$ where $d(x,y)$ is the metric induced on $M$ by $g$.
Then the slice $K(\Phi) |_s$ is given by the sheaf $1_{Z_s}$ and the continuation map from time-$s$ to time-$0$ is given by
the canonical map $1_{Z_s} \rightarrow 1_\Delta$. 
To get a description of the continuation maps for the positive time, we note that 
$\sHom(1_\Delta, p_1^* \omega_M) = \Delta_* \Delta^! p_2^! 1_M = 1_\Delta$ by (i) of Lemma \ref{omegap}.
Since $H(x,\xi)$ is time-independent, we conclude by the uniqueness statement that, for small $0< t <<1$,
the time-$t$ continuation map is given by
$1_\Delta = \sHom(1_\Delta, p_1^* \omega_M) \rightarrow \sHom(1_{Z_{-t}}, p_1^* \omega_M) = K(\Phi)|_t$. 
\end{example}

A corollary of the GKS sheaf quantization construction 
is that contact isotopies act on sheaves and the action is compatible with the microsupport:
\begin{corollary}[{\cite[(4.4)]{Guillermou-Kashiwara-Schapira}}]
Let $\Phi:S^* M \times I \rightarrow S^* M$ be a contact isotopy.
Then the convolution 
\begin{align*}
K(\Phi)|_t \circ (-): \Sh(M) &\rightarrow \Sh(M) \\
F &\mapsto K(\Phi)|_t \circ F
\end{align*} is an equivalence whose inverse is given by $K(\Phi^{-1})|_t \circ (-)$.
For a sheaf $F \in \Sh(M)$,
there is an equality $\msnz(K(\Phi) \circ F) = \Lambda_\Phi \circ \dot{\ms}(F)$. 
In particular, if we set $F_t \coloneqq (K(\Phi) \circ_M F) |_{M \times \{t\}}$, then $\msif(F_t) = \phi_t \msif(F)$ for $t \in I$.
Furthermore, if $F$ has compact support, then so does $F_t$ for all $t \in I$.
\end{corollary}

We will consider the notion of \textit{wrapping} for sheaves.
Recall that in contact geometry, 
a wrapping is usually referring to a one-parameter deformation of Legendrians $L_t$ in a contact manifold $Y$.
The wrapping is positive (resp. negative) if $\alpha(\partial_t L_t) \geq 0$ (resp. $\alpha(\partial_t L_t) \leq 0$) 
for some compatible contact form $\alpha$.
An exercise is that such a deformation $L_t$ can always be extended to a contact isotopy $\Phi$ on $Y$.
Since deformations of singular isotropics are not yet available at this moment,
we first consider globally defined contact isotopies on $S^* M$, 
and then use them to deform sheaves through GKS sheaf quantization. 

\begin{remark}
The term wrapping comes from the example on $S^* (\RR/\ZZ) \cong \RR \times \{\pm \infty\}$ with the isotopy given by
$\phi_t(x,\infty) = (x + t,\infty)$ and $\phi_t(x,-\infty) = (x - t, -\infty)$.
In this paper, we will use the term \textit{positive/negative wrapping} to mean either a positive/negative isotopy, 
a family of sheaves induced by such an isotopy, 
or the corresponding family of singular isotropics of those sheaves by taking $\msif(-)$.
Since these two notions are dual to each other, we will mainly work with positive isotopies and simply refer them as
\textit{wrappings} when the context is clear. 
\end{remark}

We will consider the totality of all such wrappings in the next section.
For now, we consider wrapping a single sheaf and develop a perturbation trick which we will use later.
Let $\Phi: S^* M \times I \rightarrow S^* M$ be a contact isotopy
and denote its GKS sheaf quantization by $K(\Phi)$.
For $F \in \Sh(M)$, the convolution $K(\Phi) \circ F$ is an object in $\Sh(M \times I)$,
which we think of it as an isotopy of $F$ and denote it by $F^\Phi$ for simplicity.
By abusing the notation, we write $F_t = i_t^* ( F^\Phi) \in \Sh(M)$ where $i_t: M \times \{t\} \hookrightarrow M \times I$
is the slice at $t$ so that $F_0 = F$.
When $\Phi$ is positive, the expression of $\Lambda_\Phi$ implies that $ \ms(K(\Phi) ) \subseteq \{ \tau \leq 0\}$
and there is continuation map $K(\Phi)_s \rightarrow K(\Phi)_t$ for $s \leq t$ and it induces
continuation maps $F_s \rightarrow F_t$ for $s \leq t$ on $F$.

\begin{proposition}[Perturbation trick, {\cite[Lemma 2.10]{Zhou2}}\footnote{The author thanks Wenyuan Li for informing him, during revisiting the paper, that this proposition has been previously obtained by Peng Zhou, using the microlocal Morse lemma, in the reference cited.}] \label{pert} 
Let $F$ and $G \in \Sh(M)$ be sheaves such that $\supp(F) \cap \supp(G)$ is compact.
Let $\Phi: S^* M \times I \rightarrow S^* M$ be a positive isotopy such that
$\varphi_t(\msif(F) )  \cap \msif(G) = \varnothing$ for $t > 0$.
Then the continuation map $F \rightarrow F_{t_0}$ induces an isomorphism
$$\Hom(G,F) \xrightarrow{\sim} \Hom(G,F_{t_0})$$
for $t_0 > 0$.
\end{proposition}

\begin{remark}
We note that when $G$ is cohomologically constructible with perfect stalks,
the object $\Hom(G,F_{t_0})$ from the last proposition is the same as 
$\Gamma(M; \sHom(G,1_M) \otimes F_{t_0})$ by (7) of Proposition \ref{mses}.
This can been seen as a sheaf-theoretic analogue of the procedure of making
Lagrangians intersect transversally. 
\end{remark}

Before we start the proof, we recall some constancy results from general sheaf theory.
Let $f: X \rightarrow Y$ be a continuous map between locally compact Hausdorff space.
Denote by $\Sh_f(X)$ the subcategory of $\Sh(X)$ consists of objects $F$ satisfying the condition
$F |_{f^{-1}(y)} \in \Loc( f^{-1}(y))$ for all $y \in Y$.
We note that $f^* G \in \Sh_f(X)$ for $G \in \Sh(Y)$.

\begin{proposition}[{\cite[Proposition 2.7.8]{Kashiwara-Schapira1}}]\label{cfl} 
Assume there is an increase sequence of closed subsets $\{X_n\}$ such that $X_n \subseteq \Int (X_{n+1})$,
$X = \cup_n X_n$, and $f_n \coloneqq f|_{X_n}$ is proper with contractible fibers.
Then, the adjunction $f_* : \Sh_f(X) \leftrightarrows \Sh(Y): f^*$ is an equivalence of categories.
\end{proposition}

\begin{corollary}\label{constc}
Let $M$, $B$ be manifolds and assume $B$ is contractible.
Let $p: M \times B \rightarrow M$ denote the projection.
Then $F \in \Sh(M \times B)$ is of the form $p^* G$ for some $G \in \Sh(M)$
if and only if $F |_{\{x\} \times B}$ is locally constant for all $x \in M$.
In this case, $G = p_* F$.
\end{corollary}


By Example \ref{loceg}, Corollary \ref{constc}, and (2) of Proposition \ref{ncmses}, we conclude:
\begin{lemma}\label{constbms}
Let $B$ be a contractible manifold and $p: M \times B \rightarrow M$ be the projection.
A sheaf $F \in Sh(M)$ satisfies $p^* p_* F \xrightarrow{\sim} F$ if and only if $\ms(F) \subseteq T^* M \times 0_B$.
\end{lemma}

\begin{proof}[Proof of Proposition \ref{pert}]
We note that the microsupport estimation $\msnz( F^\Phi) \subseteq \Lambda_\Phi \circ \msnz(F)$
implies that the sheaf $F^\Phi$ is $I$-noncharacteristic.
Thus by Corollary \ref{nclp}, the canonical map $F \rightarrow \lmi{t > 0} \, F_t$ is an isomorphism.
Since taking global sections $\Gamma(M;-)$ preserves limits,
the continuation map $\Hom(G,F) \xrightarrow{\sim} \Hom(G,F_{t_0})$ can be obtained from the composition of
$$\sHom(G,F_0) = \sHom(G, \lmi{t > 0} F_t) = \lmi{t > 0} \sHom(G,F_t) \rightarrow \sHom(G,F_{t_0})$$
where the last morphism is the projection from the limit diagram.
We claim that this diagram is actually a constant diagram after applying $\Gamma(M;-)$.

Consider the sheaf $\sHom(p^* G, F^\Phi) \in \Sh(M \times I)$.
Apply (7) of Proposition \ref{mses} and the microsupport estimation
$$ \ms(\sHom(p^* G, F^\Phi) ) 
\subseteq \left( \ms(F^\Phi) + (- \ms(G) \times 0_I ) \right)$$
implies that $\sHom(p^* G, K(\Phi) \circ F) \in \Sh_{T^* M \times T^*_\leq I}(M \times I)$ 
is also $I$-noncharacteristic.
A similar computation as in Corollary \ref{nclp} implies that 
continuation maps obtained from applying $\sHom(G,-)$ to the continuation maps of $F$
is the same as those associated to $\sHom(p^* G, F^\Phi)$.
Denote by $q: M \times I \rightarrow I$ the projection to $I$.
The assumption that $\supp(G) \cap \supp(F)$ is compact implies that
$q$ is proper on $\supp\left( \sHom(p^* G, F^\Phi) \right)$ so we may apply the microsupport
estimation to conclude, by Lemma \ref{constbms} above, that
$q_*  \sHom(p^* G, F^\Phi)$ is a constant sheaf.
Finally, properness and $I$-noncharacteristic assumptions that 
$$\left(q_* \sHom(p^* G, F^\Phi) \right)_t = \Hom(G,F_t)$$ which concludes the proof.

\end{proof}

\subsection{The category of positive wrappings}

We will define the category of positive wrappings whose morphisms will be given by concatenation. 
In order to define concatenation easily, we assume that the isotopies are constant near the end points,
and the interval $I$ will be a closed interval for this section.
This requirement doesn't lose much information since 
for any positive contact isotopy $\Phi: S^* M \times [0,1] \rightarrow S^* M$,
one can always make it have constant ends through a homotopy of positive isotopies.
For example, pick a non-decreasing $C^\infty$ function $\rho$ on $\mathbb{R}$ such that $\rho |_{(-\infty,1/3]} \equiv 0$
and $\rho |_{[2/3,\infty)} \equiv 1$,
then an example of such a modification is given by $\tilde{\Phi}(x,\xi,t,s) = \Phi(x,\xi, (1-s) t + s \rho(t))$.
By Proposition \ref{homotopy_independence_of_continuation_maps}, they induce equivalent continuation maps and two such identifications can itself be identified
by a similar consideration and so on.
Thus, when we mention isotopies obtained through natural constructions 
such as by integrating from a time-independent vector field,
we will implicitly assume such a deformation procedure.

\begin{definition}
Let $I = [t_0,t_1]$, $J = [s_0, s_1]$ be two closed intervals.
We use $I \# J$ to denote the concatenated interval $(I \amalg J) / \{ t_1 \sim s_0 \}$.
For isotopies $\Phi: S^* M \times I \rightarrow S^* M$, $\Psi: S^* M \times J \rightarrow S^* M$,
the \textit{concatenation} isotopy $\Psi \# \Phi: S^* M \times (I \# J) \rightarrow S^* M$ is the isotopy which is given by
$$ (\Psi \# \Phi)(x,\xi,t) =
\begin{cases}
\Phi(x,\xi,t), &t \in I, \\
\Psi( \Phi(x,\xi,t_1),t), &t \in I^\prime.
\end{cases}
$$
If $I = J$, one can also define the pointwise composition $\Psi \circ \Phi: S^* M \times I \rightarrow S^* M $
by $$ (\Psi \circ \Phi) (x,\xi,t) = \Psi( \Phi(x,\xi,t),t) ).$$
\end{definition}

Note that, up to a scaling, $\Psi \circ \Phi$ and $\Psi \# \Phi$ are homotopic to each othert and,
if both $\Phi$ and $\Psi$ are positive, they are homotopic through positive isotopies.

\begin{definition}
Let $\Omega \subseteq S^* M$ be an open subset. 
We say a contactomorphism $\varphi: S^* M \rightarrow S^* M$ is \textit{compactly supported} on $\Omega$
if $\varphi$ equals $\id_{S^* M}$ outside a compact set $C$ in $\Omega$.
Similarly, a contact isotopy $\Phi: S^* M \times I \rightarrow S^* M$ is \textit{compactly supported} on $\Omega$ 
if $\varphi_t = \id$ outside a fixed compact set $C$ in $\Omega$ for all $t \in I$.
\end{definition}

\begin{definition}
We define \textit{the category $W(\Omega)$ of compactly supported positive wrappings} on $\Omega$ as follows:
An object of $W(\Omega)$ is a pair $(\varphi, [\Phi])$ such that $\varphi$ is a compactly supported 
contactomorphism and $[\Phi]$ is a homotopy class of compactly supported isotopies, 
defined on a closed interval $I$, having $\varphi$ as its end point and realizing it as Hamiltonian.
Note that the degenerate case $I = \{ * \}$ is allowed.  
We will often simply write $(\varphi,[\Phi])$ by $\Phi$ without emphasizing that it is a homotopy class 
through the paper when it is irrelevant. 
A $1$-morphism $\Psi: [\Phi_0] \rightarrow [\Phi_1]$ is a positive isotopy $\Psi$
such that $[\Phi_1] = [\Psi \# \Phi_0]$. 
Composition of $1$-morphisms is given by concatenation.
For $\Psi_0, \Psi_1: [\Phi_0] \rightarrow [\Phi_1]$, a $2$-morphism is a positive family of isotopies 
$\Theta: S^* M \times I \times J \rightarrow S^* M$ which is constant near the end points on the $J$-direction
such that $\Theta(-,t,s_i) = \Psi_i(-,t)$ and $\Theta(-,t_i,s) = \Phi_i(-)$, $i = 0, 1$.
Here $t_i$ and $s_i$ are the end points of $I$ and $J$.
An $n$-morphism will be a homotopy between $n-1$ morphisms with the 
obvious boundary conditions and similar constancy requirements.
\end{definition}


We will show that the category $W(\Omega)$ is filtered so colimits over it are, by definition, filtered colimits.
Recall that a $1$-category $\sC$ is filtered if, 
\begin{enumerate}
\item $\sC$ is non-empty,
\item for any $X, Y \in \sC$, there is $Z \in \sC$ with morphisms $X \rightarrow Z$ and $Y \rightarrow Z$, and, 
\item for any more morphism $f, g: X \rightarrow Y$, there exist $h: Y \rightarrow Z$ such that $h \circ f = h \circ g$.
\end{enumerate}

This is the same as the condition that
for any (ordered) $n$-simplex $K$, $n \in [-1,1]$ and any functor $F: K \rightarrow \sC$, 
there is an extension $\hat{F}$ on $K^\triangleright$
where $K^\triangleright$ is the $n+1$ simplex obtained by adding a final cone point to $K$.
For example, we can realize a pair of morphisms $f, g: X \rightarrow Y$ 
as a hollowed triangle consisting of vertices $X, X, Y$ and edges $\id_X, f, g$ without the presence of the face. 
A final cone point $Z$ provides a morphism $h:Y \rightarrow Z$ for the edge between $Y$ and $Z$.
The existence of the three new faces and the fact that the only $2$-morphism 
in a $1$-category is the strict equality implies  $h \circ f = h \circ g$.
Now, we provide the definition for a category to be filtered in the $\infty$-categorical setting.

\begin{definition}
A category $\sC$ is filtered if for any simplex $K$ and any functor $F: K \rightarrow \sC$,
there is an extension $\hat{F}: K^\triangleright \rightarrow \sC$.
\end{definition}

To illustrate what additional conditions are imposed on $n$-morphisms, for $n \geq 2$, 
we consider the following case of a $2$-simplex.

\begin{example}
Consider the case when $K = S^2$ is the $2$-sphere, or more precisely,
when $K = \Delta^2$ is the standard $2$-simplex such that the base face has three vertices being a fixed object $X$,
three edges being $\id_X$, and the face being the trivial identification.
This is essentially the situation that there are objects $X$, $Y$, a $1$-morphism $f: X \rightarrow Y$,
with a non-trivial $2$-automorphism $T$ on $f$.
The condition of $\sC$ being filtered means that there exist $g: Y \rightarrow Z$ such that
the auto equivalence $g \circ T$ on $g \circ f$ is trivial, that is, $g \circ T = \id_{g \circ f}$.
\end{example}

The following proposition is a version of \cite[Lemma 3.27]{Ganatra-Pardon-Shende1}.

\begin{proposition}\label{wc:fl}
The category $W(\Omega)$ is filtered. 
\end{proposition}

\begin{proof}
Similarly to the situation in classical algebraic topology, it is sufficient to check the case when $K = S^n$, 
the $n$-sphere for $n \in \ZZ_{\geq 0}$. 

When $n = 0$, we are given two homotopy classes of contact isotopies $\Phi_0$ and $\Phi_1$ with the same end point $\varphi$, 
and the goal is to find another contact isotopy $\Phi$ and two positive contact isotopies $\Psi_0$ and $\Psi_1$ such that
$[\Phi] = [\Psi_0 \# \Phi_0] = [\Psi_1 \# \Phi_1]$.
We first notice that, up to a rescaling, $[\Phi_0^{-1} \# \Phi_0] = [\Phi_1^{-1} \# \Phi_1] = [\id_{S^* M}]$.
So it is sufficient to modify $\Phi_0^{-1}$ and $\Phi_1^{-1}$ by composing some $\Phi^\prime$ so that 
$\Phi^\prime \circ \Phi_0^{-1}$ and $\Phi^\prime \circ \Phi_1^{-1}$ are positive.
Let $H_0$, $H_1$ denote their Hamiltonians.
Since $\Phi_0^{-1}$ and $\Phi_1^{-1}$ are compactly supported, there exists a compact set $C \subseteq \Omega$
such that $H_0$ and $H_1$ are zero outside $C$.
Pick a positive real number $r$ such that $r > \max(\abs{H_0},\abs{H_1})$,
relative compact open sets $U$, $V$ in $S^* M$ such that $C \subset U \subseteq \overline{U} \subseteq V \subseteq \Omega$,
and a bump function $\rho$ such that $\rho|_U \equiv 1$ and $\rho \equiv 0$ outside $V$.
The contact isotopy $\Phi^\prime$ generated by $r \rho$ will satisfy the requirement by the Leibniz rule.

When $n > 0$, we are given a family of morphism $\Psi_\theta: \Phi_0 \rightarrow \Phi_1$ 
parametrized by $\theta \in S^{n-1}$
such that $[\Phi_1] = [\Psi_\theta \# \Phi_0]$, and we have to show that, by possibly further concatenation, 
this family can be made to be null-homotopy through positive isotopies.
By precomposing  $\Phi_0^{-1}$, we may assume there is an $S^{n-1}$-family of positive isotopies $\Psi_\theta$ 
and a fixed (not necessarily positive) isotopy $\Phi$, such that, for each $\theta \in S^{n-1}$, 
there exists a homotopy $\Sigma_\theta: S^* M \times I \times [0,1] \rightarrow S^* M$ 
connecting $\Phi$ to $\Psi_\theta$.
We can extend this map to a $D^n$-family of isotopy $\Sigma: S^* M \times I \times D^n \rightarrow S^* M$
by $\Sigma(x,\xi,t,r \theta) = \Sigma_\theta(x,\xi,t,r)$ where we write elements in $D^n$ 
as $r \theta$ by $r \in [0,1]$ and $\theta \in S^{n-1}$. 
Now the same compactness argument as before shows that there is a positive isotopy $\Phi^\prime$
such that $\Phi^\prime \circ \Sigma$ is positive. 
\end{proof}

Let $F: \sC \rightarrow \sD$ be a functor.
For any diagram $p: \sD \rightarrow \mathscr{E}$, the colimits $\clmi{C} (p \circ F)$ and  $\clmi{D}\, F$ exist if either one exists.
Thus, it is well-defined to write the canonical map $\clmi{C} (p \circ F) \rightarrow \clmi{D}\, F$.

\begin{definition}
A functor $F: \sC \rightarrow \sD$ is cofinal if, for any diagram $p: \sD \rightarrow \mathscr{E}$,
the canonical map $\clmi{C} (p \circ F) \rightarrow \clmi{D} \, F$ is an isomorphism.
\end{definition}

In the $1$-categorical setting, a more classical convention is that a functor is cofinal if and only if,
\begin{enumerate}
\item for any $d \in \sD$, there exists $c \in \sC$ and a morphism $d \rightarrow F(c)$, 
\item for any morphism $f, g: d \rightarrow F(c)$, there exists $h: c \rightarrow c^\prime$ such that $F(h) \circ f = F(h) \circ g$.
\end{enumerate}
An equivalent way of saying it is that the fiber product $\sC \times_{ \sD} d/\sD$ 
is non-empty and connected for all $ d \in \sD$.
Here, $d/\sD$ is the over category whose objects are morphisms of the form $d \rightarrow d^\prime$ and 
a morphism $h: (f: d \rightarrow d^\prime) \rightarrow (g: d \rightarrow d^{\prime \prime})$ 
is given by a morphism $h: d^\prime \rightarrow d^{\prime \prime}$ such that $h \circ f = g$,
the fiber product is taken over the canonical projection $d/\sD \rightarrow \sD$ by 
$(d \rightarrow d^\prime) \mapsto d^\prime$ and $F$.
Recall a $1$-category is said to be connected if the associated $1$-groupoid (by formally inverting morphisms) is connected.
The equivalence of these definitions is the Quillen's theorem A.
In the $\infty$-categorical setting it states:

\begin{theorem}[Quillen's Theorem A]
A functor $F: \sC \rightarrow \sD$ is cofinal if and only if 
the fiber product $\sC \times_{ \sD} d/\sD$ is contractible for any $d \in \sD$.
\end{theorem}

Now consider the following construction:
For $n = 1, 2, \cdots$, take a family of open set $\Omega_n \subseteq S^* M$ such that 
$\Omega_n \subseteq \overline{\Omega_n} \subseteq \Omega_{n +1}$,
$\cup_{n \in \ZZ \geq 0 } \, \Omega_n =  S^* M$, and $\Omega_n \subseteq S^* M$ is relative compact.
For $n> 0$, pick bump function $\rho_n$ such that $\rho_n \leq \rho_{n+1}$,
$\rho_n |_{\Omega_n} \equiv n$, and vanishes outside $\Omega_{n+1}$.
Let $\Phi_n: S^* M \times [0,n] \rightarrow S^* M$ be the isotopy generated by $\rho_n$. 
Since $\rho_1 \leq \rho_2 \leq \cdots \rho_n \leq \cdots$, 
there exists positive isotopy $\Psi_n: S^* M \times [n,n+1] \rightarrow S^* M$ such that 
the $\Phi_n$'s and $\Psi_n$'s form a sequence $\id \xrightarrow{\Psi_0} \Phi_1 \xrightarrow{\Psi_1} \cdots$ in $W(\Omega)$.
That is, the above data organizes to a functor $\Phi : \ZZ_{\geq 0} \rightarrow W(\Omega)$.

\begin{lemma}\label{cfnc}
The functor $\Phi: \ZZ_{\geq 0} \rightarrow W(\Omega)$ is cofinal.
\end{lemma}

\begin{proof}
By Quillen's Theorem A, we need to show that $\ZZ_{\geq 0} \times_{W(\Omega)} \left( \Phi/W(\Omega) \right)$ 
is contractible.
Let $\Phi \in W(\Omega)$ and let $H$ denote its Hamiltonian.
Since $\Phi$ is compactly supported there exist a compact set $C \subseteq S^* M$
such that $H$ vanishes outside $C$.
Pick $n$ large such that $C \subseteq \Omega_n$ and $\max(H) \leq n$.
Then the factorization $\Phi_n = ( \Phi_n \circ \Phi^{-1}) \circ \Phi$ provides an morphism $ \Phi \rightarrow \Phi_n $
since composition is homotopic to concatenation.
Thus, the fiber product $\ZZ_{\geq 0} \times_{W(\Omega)} \left( \Phi/W(\Omega) \right)$
is equivalent to $\{ n \in \ZZ_{\geq 0 } | \rho_n \geq H \}$. 
Since $\rho_{n+1} \geq \rho_n$ by our construction, 
the latter is equivalent to the poset of integers larger than $\min \{n: \rho_n \geq H \}$ and is contractible.
\end{proof}

Now we quantize the above construction.
For $(\varphi,[\Phi]) \in W(\Omega)$ where $\Phi$ is defined on $M \times [t_0,t_1]$,
set $w(\Phi) \coloneqq K(\Phi) |_{t = t_1}$ to be the restriction of the GKS sheaf quantization $K(\Phi_i)$ at the end point.
We note that since we require the end point $\varphi$ to be fixed, the sheaf $w(\Phi) \in \Sh(M \times M)$ depends only
on the homotopy class $[\Phi]$ by formula \ref{asso_lag} and Lemma \ref{constbms}.
For a morphism $\Psi: \Phi_0 \rightarrow \Phi_1$, since $\Psi$ is positive, 
formula \ref{asso_lag} again implies there is a continuation map $c(\Psi):w(\Phi_0) \rightarrow w(\Phi_1)$.
Similarly, for a pair of morphisms $\Psi_0, \Psi_1: \Phi_0 \rightarrow \Phi_1$, 
if there is a homotopy $\Theta$ between $ \Psi_0$ and $\Psi_1$,
Proposition \ref{homotopy_independence_of_continuation_maps} implies that $K(\Theta)$ provides an identification
between the continuation maps $c(\Psi_0), c(\Psi_1): w(\Phi_0) \rightarrow w(\Phi_1)$.

\begin{definition}
Organizing the above construction, we obtain a functor $w: W(\Omega) \rightarrow \Sh(M \times M)$ 
sending an object $\Phi$ to the corresponding sheaf kernel $w(\Phi)$, 
a $1$-morphism $\Psi: \Phi_0 \rightarrow \Phi_1$ to the continuation map $c(\Psi): w(\Phi_1) \rightarrow w(\Phi_0)$
, and higher morphisms to higher equivalences of continuation maps.
We will refer this functor as the \textit{wrapping kernel functor}.
\end{definition}

For a sheaf $F \in \Sh(M)$ and a contact isotopy $\Phi: S^* M \times [0,1] \rightarrow S^* M$, 
convolving with $K(\Phi)$ produces a sheaf $K(\Phi) \circ F$ on $\Sh(M \times [0,1])$
and we use the notation $F_t = (K(\Phi) \circ F)|_t$ for $t \in [0,1]$ as before. 
Base change \ref{bsch} and compatibility of six-functor formalism \ref{cp6f} implies that 
there is an identification 
$$(K(\Phi) \circ F) |_{t = 1} = K(\Phi) |_{t = 1} \circ F = w(\Phi) \circ F$$
 functorial on $\Phi$ and $F$.
When $\Psi: \Phi_0 \rightarrow \Phi_1$ is a positive family of isotopies, we use
$c(\Psi,F): w(\Phi_0) \circ F \rightarrow w(\Phi_1) \circ F$ to denote the induced continuation map.
To simplify the notation, we sometimes use $F^{w(\Phi)}$ to denote $w(\Phi) \circ F$. 
When there is no need to specify the isotopy, we simply write it as $F^w$.
Similarly, when $\Psi$ is unspecified, we simply write $c: F^w \rightarrow F^{w^\prime}$ for the continuation map. 
We prove a locality property which we will use later.

\begin{proposition}\label{loclem}
Let $\Phi_0$, $\Phi_1: S^* M \times I \rightarrow S^* M$ be contact isotopies.
If $\Phi_0 = \Phi_1 $ for on $\msif(F) \times I$,
then $w(\Phi_0) \circ F = w(\Phi_1) \circ F$.
Similarly, let $\Psi_0, \Psi_1$ be positive contact isotopies.
If $\Psi_0= \Psi_1 $ on an open neighborhood $\Omega_0$ of $\msif(F)$, then $c(\Psi_0, F) = c(\Psi_1,F)$.
\end{proposition}

\begin{proof}
We abuse the notation and use $\Phi_i$ to denote the corresponding homogeneous symplectic isotopies.
By convolving with $\Phi_1^{-1}$, it is enough to assume $\Phi_1 = \id_{S^* M}$ and show that $F^{w(\Phi_0)} = F$.
We have $\ms(K(\Phi_0) \circ F) \subseteq \ms(K(\Phi_0) ) \circ \ms(F) \subseteq \ms(F) \times 0_I$.
Hence by Proposition \ref{constbms}, $F_t$ is constant along constant on $t$.

Similarly, let $H_i \geq 0$ denote the Hamiltonian of $\Psi_i$, $i = 0, 1$.
Let $\Psi_s$ be the homotopy of isotopies between $\Psi_0$ and $\Psi_1$
generated by the Hamiltonian $\widehat{H}(x,\xi,t,s) = (1 - s) H_0(x,\xi,t) + s H_1(x,\xi,t)$.
Since $\Psi_0= \Psi_1 $ on $\Omega_0$, $\Psi_s |_{\Omega_0}$ is constant on $s$.
Thus Proposition \ref{homotopy_independence_of_continuation_maps} applies to $K(\Psi_s) \circ F$ and we conclude $c(\Psi_0,F) = c(\Psi_1,F)$.
\end{proof}

The above construction defines a functor from $\Sh(M)$ to $[W(\Omega),\Sh(M)]$ 
since the expression $w(\Phi) \circ F$ is functorial on $F$.
Further composing with the functor of taking limits and colimits defines functors $\wrap^{\pm}(\Omega): \Sh(M) \rightarrow \Sh(M)$.
Since the subcategory of $W(\Omega)$ consists of objects $(\varphi,[\Phi])$ 
such that there exists a positive isotopy $\Phi$ representing $[\Phi]$ is cofinal, we can informally write the formula by
$$ \wrap^+(\Omega)F = \clmi{F \rightarrow F^w} \ F^w, \ \wrap^-(\Omega)G = \lmi{G^{w^-} \rightarrow G} G^{w^-}.$$
With this definition, we can generalize Proposition \ref{pc} to Theorem \ref{w=ad}.


\begin{proof}[\hypertarget{proof_w=ad}{Proof} of Theorem \ref{w=ad}]
Set $\Omega = S^* M \setminus X$ and let $F \in \Sh(M)$. 
We first show that for any $(x,\xi) \in \Omega$, $(x,\xi) \not\in \msif(\wrap^\pm(F))$, i.e., 
for any function $f$ defined near $x$ such that $f(x) = 0 $ and $df_x \in \Rp \xi$, 
the restriction map $(\wrap^\pm(F))_x \rightarrow \Gamma_{ \{f < 0 \} } (\wrap^\pm(F))_x $ is an isomorphism.
Since the situation is local and $df_x \neq 0$, by changing a coordinate, we may assume
$f = x_1$ the first coordinate function near $x = 0$.
Pick a family of open balls $U_i$ centered at $x$ such that $U_i \supseteq \overline{U_{i+1}} \supseteq U_{i+1}$ and $\cap_i U_i = \{x\}$.
The stalk $ \Gamma_{ \{f < 0 \} } \left( \wrap^\pm(F) \right)_x $ can be computed by 
the colimit $\clmi{i} \, \Gamma\left(U_i \cap \{x_1 < 0 \}; \wrap^\pm(F)\right)$.

We first prove the negative case. 
For each $i$, we take a small positive wrapping $\Phi_i$ supported in $\Omega$ 
such that $w(\Phi_i) \circ 1_{U_i \cap \{x_1 < 0 \} }= 1_{\tilde{U}_i}$ 
with $0 \in \tilde{U}_i$ and $\tilde{U}_i$ shrinks to $x$ as $i \rightarrow \infty$.
For example, take $U_i \times C_i$ in $\Omega$ containing $(x,\xi)$ where
$\{C_i\}$ is a family of small balls on the fiber direction with a condition similar to the $\{ U_i \}$.
For each $i$, pick a bump function $\rho_i$ on $S^* M$ supported on $U_i \times C_i$ and equals $1$ near $(x,\xi)$.
Take $H_i$ to be the Hamiltonian associated to the Reeb flow with shrinking speed and modify it to $\rho_i H_i$.
Finally, take $\Phi_i$ to be the isotopy associated to $\rho_i H_i$.

We compute,
\begin{align*}
\Gamma\left(U_i \cap \{x_1 < 0 \}; \wrap^-(F)\right)
&= \lmi{ W(\Omega)}  \Hom \left( 1_{U_i \cap \{x_1 < 0 \} }, w(\Phi)  \circ F \right) \\
&= \lmi{  W(\Omega)}  \Hom \left( w(\Phi_i) \circ  1_{U_i \cap \{x_1 < 0 \} },
w(\Phi_i)  \circ w(\Phi) \circ F \right) \\
&= \lmi{ W(\Omega)}  \Hom \left(  1_{\tilde{U}_i },  w(\Phi_i \circ|_I \Phi) \circ F \right) \\ 
&= \Gamma \left( \tilde{U}_i ; \wrap^-(F) \right).
\end{align*}
Here we use the fact that $w(\Phi_i) \circ$ is an equivalence for the second equation.
For the last equation, we use the fact that negative wrappings of the form $\Phi_i \circ \Phi$ are initial in $W(\Omega)$.
Take $i \rightarrow \infty$ and we conclude $(\wrap^-(F))_x \xrightarrow{\sim} \Gamma_{ \{f < 0 \} } (\wrap^-(F))_x $
and $\Gamma_{ \{f \geq 0 \} } (\wrap^-(F))_x  = 0$.

Now we turn to the positive case.
We take the same family of $U_i$, $\Phi_i$ and $\tilde{U}_i$, and compute,
\begin{align*}
\Gamma\left(U_i \cap \{x_1 < 0 \}; \wrap^+(F)\right)
&=  \Hom \left( 1_{U_i \cap \{x_1 < 0 \} }, \clmi{ W(\Omega)} \, ( w(\Phi)  \circ F) \right) \\
&=  \Hom \left( w(\Phi_i) \circ 1_{U_i \cap \{x_1 < 0 \} },
w (\Phi_i) \circ ( \clmi{  W(\Omega)} \  w(\Phi) \circ F ) \right) \\
&=  \Hom \left(1_{\tilde{U}_i \cap \{x_1 < 0 \} },
\clmi{ W(\Omega)} (w(\Phi_i \circ|_I \Phi) \circ F ) \right) \\
&=  \Gamma\left(\tilde{U}_i; \wrap^+(F)\right). 
\end{align*}
We use the fact that $w(\Phi_i) \circ $ is a left adjoint so it commutes with colimits for the third equation.
Take $i \rightarrow \infty$ and we get  $(\wrap^+(F))_x \xrightarrow{\sim} \Gamma_{ \{f < 0 \} } (\wrap^+(F))_x $
and $\Gamma_{ \{f \geq 0 \} } (\wrap^+(F))_x  = 0$.

From the above computation, we see that $\wrap^\pm(\Omega) : \Sh(M) \rightarrow \Sh (M)$ factorizes to $\Sh_X (M)$.
Finally, we show that $\wrap^+(\Omega) \dashv \iota_* \dashv \wrap^-(\Omega)$. 
Take $G \in \Sh_X(M)$ and $F \in \Sh(M)$.
We compute,
\begin{align*}
\Hom \left(G,\wrap^-(F) \right)
&= \Hom \left(G, \lmi{  W(\Omega)} w(\Phi) \circ F \right) \\
&=  \lmi{ W(\Omega)} \Hom \left(G, w(\Phi) \circ F \right) \\
&= \lmi{W(\Omega)} \Hom \left(w(\Phi^{-1}) \circ G, F \right) \\
&= \lmi{ W(\Omega)} \Hom \left( G, F \right) = \Hom \left( \iota_* G, F\right).
\end{align*}
The second to last equality is implied by Proposition \ref{loclem} since $\Phi$ is compactly supported away from 
$\Lambda \supseteq \msif(G)$.
A similar computation shows that $$\Hom \left( \wrap^+(F), G \right) = \Hom \left(F, \iota_* G \right).$$
\end{proof}

\begin{remark}
We mention that, aside from the prototype cases \cite{Tamarkin1, Guillermou-Kashiwara-Schapira, Guillermou-Schapira} 
mentioned in the introduction,
special cases for such geometric descriptions can be found in, for example, \cite{Kuwagaki1} 
in the setting of toric homological
mirror symmetry which are defined by using the group structure of the torus 
and are crucial for matching the data with the coherent side.  
\end{remark}

\begin{notation}
Recall that we use $\iota_X^*: \Sh(M) \rightarrow \Sh_X(M)$ to denote the left adjoint of the inclusion
${\iota_X}_*: \Sh_X(M) \hookrightarrow \Sh(M)$ for a conic closed subset $X \subseteq T^* M$.
By the above Theorem $\ref{w=ad}$, when $\Lambda \subseteq S^* M$ is a singular isotropic,
we will use the notation $\wrap_\Lambda^+(M) = \iota_\Lambda^*: \Sh(M) \rightarrow \Sh_\Lambda(M)$
when emphasizing that $\iota_\Lambda^*$ is given by wrappings.
When there is no ambiguity for the ambient manifold $M$, we simply write it as $\wrap_\Lambda^+$.
We will use a similar notation for the right adjoints.
\end{notation}

\section{The category of wrapped sheaves}

In this section, we make our main construction of this dissertation and
mimic the definition of $\WF(T^* M, \Lambda)$, the (partially) wrapped Fukaya category associated to the pair $(T^* M, \Lambda)$,
and define the category of wrapped sheaves $\wsh_\Lambda(M)$ using the techniques developed in Section \ref{iso_sh}.
We note that although this and the next section are logically independent from \cite{Ganatra-Pardon-Shende1,
Ganatra-Pardon-Shende2}, many of the proofs are adaptations from the wrapped Fukaya setting to the wrapped sheaves setting.

\subsection{Definition}

Let $M$ be a real analytic manifold and $\Lambda \subseteq S^* M$ be a closed subset.
Let $\widetilde{\wsh}_\Lambda (M)$ be the small subcategory of $\Sh(M)$ 
generated under finite colimits and retracts by sheaves of the form $F^{w(\Phi)}$ 
where $F$ is a sheaf with compact support such that $\msif(F)$ is a subanalytic singular isotropic 
and $\msif(F) \cap \Lambda = \varnothing$,
and $\Phi$ is a contact isotopy compactly supported away from $\Lambda$.
To encode the effect of wrappings,
we take $$\sC_\Lambda(M) \coloneqq \langle \cof(c(\Psi,F)) | \Psi \in \Mor(W(S^* M \setminus \Lambda)), 
F \in \widetilde{\wsh}_\Lambda (M) \rangle$$
to be the subcategory of $\widetilde{\wsh}_\Lambda (M)$ generated by the cofibers of the continuation maps.

\begin{definition}
We define the category of \textit{wrapped sheaves} associated to $(M,\Lambda)$ to be the quotient category
\begin{align*}
\wsh_\Lambda(M) 
&\coloneqq \widetilde{\wsh}_\Lambda(M)/ \sC_\Lambda(M) \\
&\coloneqq \cof \left( \sC_\Lambda(M) \hookrightarrow \widetilde{\wsh}_\Lambda(M) \right)
\end{align*}
where the $\cof$ is taken in $\st$. 
\end{definition}

\begin{remark}\label{iso=sm}
Localization identifies sheaves which are isotopic to each other:
Let $F \xrightarrow{c} F^w \rightarrow \cof(c)$ be a fiber sequence in $\widetilde{\wsh}_\Lambda(M)$ 
induced by a continuation map.
Since a quotient map is exact and $\cof(c) = 0$ in $\wsh_\Lambda(M)$, the fiber sequence
becomes $F \xrightarrow{c} F^w \rightarrow 0$ and hence $c: F \rightarrow F^w$ is an isomorphism in $\wsh_\Lambda(M)$. 
Now let $\Phi$ be any isotopy compactly supported away from $\Lambda$.
By Proposition \ref{wc:fl}, $\Phi$ can be modified to be positive by a further wrapping.
That is, there exists $\Psi: \id \rightarrow \Phi^\prime$ and $\Psi^\prime: \Phi \rightarrow \Phi^\prime$ in $W(S^* M \setminus \Lambda)$ 
and thus there are continuation maps $F \xrightarrow{c(\Psi,F)} F^{\Phi}$ and $F \xrightarrow{c(\Psi^\prime,F)} F^{\Phi^\prime}$.
As a result, the two objects $F$ and $F^\Phi$ are isomorphic in $\wsh_\Lambda(M)$.
\end{remark}

To simplify the notation, we will use  $\Hom_w$ to denote the Hom spaces of the localized category when
the context is clear.
As the localization is essentially surjective, 
we usually implicitly assume a preimage $F \in \widetilde{\wsh}_\Lambda(M)$ for objects in $\wsh_\Lambda(M)$.
By Proposition \ref{hombc}, there are identifications
$$ \Hom_w (X,Y) 
= \clmi{Y \xrightarrow{\alpha} Y^\prime } \Hom(X,Y^\prime)
= \clmi{X^\prime \xrightarrow{\beta} X } \Hom(X^\prime,Y)$$
where $\alpha$ and $\beta$ run through morphisms whose cofibers $\cof(\alpha)$, $\cof(\beta)$
 are in $\sC_\Lambda(M)$.
We will show that it is enough to take the colimit over $W(S^* M \setminus \Lambda)$ in our case, which
is the same colimit over all continuation maps $c: F \rightarrow F^w$ by cofinality.
We begin with the case of Homing out of objects in $\sC_\Lambda(M)$ and in this case 
$\Hom_w$ vanishes. 

\begin{lemma}\label{homvn}
Let $G \in \sC_{\Lambda}(M)$ and $F \in \widetilde{\wsh}_\Lambda(M)$, we have
$$ \clmi{c: F \rightarrow F^w}  \Hom(G,F^w) = 0  = \Hom_w(G,F)$$
where $F \xrightarrow{c} F^w$ runs through all continuation maps.
\end{lemma}

\begin{proof}
We first consider the case when $G \in \sC_{\Lambda}(M)$ 
is built from iterated cones and shifts of $\cof(c)$ for some continuation map $c$. 
For such a $G$, we may assume $G$ fits into a cofiber sequence $H \xrightarrow{c} H^{w(\Phi)} \rightarrow G$ by induction. 
Apply $\Hom(-,F^{w(\Phi^\prime)})$ and we obtain the cofiber sequence 
$$\Hom(G,F^{w(\Phi^\prime})) \rightarrow \Hom(H^{w(\Phi)},F^{w(\Phi^\prime)}) \rightarrow \Hom(H,F^{w(\Phi^\prime)}).$$
Then one compute 
\begin{align*}
\Hom(H^{w(\Phi)},F^{w(\Phi^\prime)}) 
&= \Hom(w(\Phi) \circ H,w(\Phi^\prime) \circ F) \\
&= \Hom( H,w(\Phi^{-1})\circ w(\Phi^\prime) \circ F) \\
&= \Hom( H,w(\Phi^{-1} \circ \Phi^\prime) \circ F).
\end{align*}

Take $\clmi{F \rightarrow F^{w^\prime} }$ over isotopies $\Phi^\prime$ of the form $\Phi \circ  \Psi$, which is cofinal, 
and we obtain the fiber sequence
$$\clmi{F \rightarrow F^{w^\prime} } \Hom(G,F^{w^\prime}) \rightarrow \clmi{F \rightarrow F^{w^\prime} } \Hom(H,F^{w^\prime}) 
\xrightarrow{\sim} \clmi{F\rightarrow F^{w^\prime} } \Hom(H,F^{w^\prime}).$$
This implies that $\clmi{F \rightarrow F^{w^\prime} } \Hom(G,F^{w^\prime}) = 0$.

Now let $G^\prime$ be a retract of the same $G$ as above.
Taking $\clmi{ F \rightarrow F^w} \Hom(-,F^w)$ makes $\clmi{ F \rightarrow F^w} \Hom(G^\prime,F^w)$
a retract of $\clmi{F \rightarrow F^{w^\prime} } \Hom(G,F^{w^\prime}) = 0$.
Since the only retract of a zero object is a zero object, 
we conclude that $\clmi{ F \rightarrow F^w} \Hom(G^\prime,F^w) = 0$.
\end{proof}

\begin{proposition}
For $F, G \in \widetilde{\wsh}_\Lambda(M)$, we have
$$\Hom_w(G,F) = \clmi{ F \rightarrow F^w} \,  \Hom(G,F^w) $$
where $F \xrightarrow{c} F^w$ runs through all continuation maps.
\end{proposition}

\begin{proof}
Consider any morphism $\alpha: G^{\prime} \rightarrow G$ such that $\fib(\alpha) \in \sC_\Lambda(M)$
which we will denote it as $G^{\prime} \xrightarrow{qis.} G$ when the exact morphism $\alpha$ is not relevant.
Now take a continuation map $c:F \rightarrow F^w$ and apply
$\Hom(-,F^w)$ to the fiber sequence $\fib(\alpha) \rightarrow G^{\prime} \rightarrow G$, 
we obtain the fiber sequence $$ \Hom(G,F^w) \rightarrow  \Hom(G^{\prime},F^w) \rightarrow \Hom( \fib(\alpha),F^w).$$
Now recall that $\Hom_w$ can be computed by either varying the first or the second factor. 
As a result, we can first take colimit over such $\alpha: G^{\prime} \rightarrow G$ and we obtain the fiber sequence
$$ \Hom(G,F^w) \rightarrow  \Hom_{w}(G,F^w) 
\rightarrow \clmi{G^{\prime} \xrightarrow{qis.} G } \Hom(\fib(\alpha),F^w).$$

Then we take colimit over $F \rightarrow F^w$ and get
$$\clmi{F \rightarrow F^w} \, \Hom(G,F^w) \rightarrow \clmi{F \rightarrow F^w} \, \Hom_{w}(G,F^w) 
\rightarrow \clmi{F \rightarrow F^w} \clmi{ G^{\prime} \xrightarrow{qis.} G } \, \Hom(\fib(\alpha),F^w).$$

Since colimits commute with each other, the above Lemma \ref{homvn} implies 
$$\clmi{F \rightarrow F^w} \clmi{G^{\prime} \xrightarrow{qis.} G  } \Hom(\fib(\alpha),G^{\prime})
= \clmi{ G^{\prime} \xrightarrow{qis.} G  } \, \clmi{F \rightarrow F^w} \, \Hom(\fib(\alpha),G^{\prime}) = 0$$ 

or, equivalently $$\clmi{F \rightarrow F^w} \, \Hom(G,F^w) \xrightarrow{\sim} 
\clmi{F \rightarrow F^w} \, \Hom_w(G,F^w) \xrightarrow{\sim}  \Hom_{w}(G,F).$$
Here we use that fact that $F \rightarrow F^w$ is an isomorphism in $\wsh_\Lambda(M)$.
\end{proof}


We note that the above construction is covariant on the open sets of $M$.
For an open set $U \subseteq M$, we set $\Lambda|_U = \Lambda \cap S^* U$.
We abuse the notation and use $\widetilde{\wsh}_\Lambda(U)$ to denote the category $\widetilde{\wsh}_{\Lambda|_U}(U)$.
When there is an inclusion of open sets $U \subseteq V$, 
objects in $\widetilde{\wsh}_\Lambda(U)$ can be naturally regarded as objects in $\widetilde{\wsh}_\Lambda(V)$
since we require them to have compact support in $U$.
We define $\sC_\Lambda(U)$ similarly and note that 
$\sC_\Lambda(U) \subseteq \widetilde{\wsh}_{\Lambda}(U) \cap \sC_\Lambda(V)$.
Thus there is a canonical map $\wsh_{\Lambda} (U) \rightarrow \wsh_{\Lambda}(V)$.
\begin{definition}
The above construction defines a covariant functor $\wsh_\Lambda: \Op_M \rightarrow \st$, 
i.e., a precosheaf with coefficient in idempotent complete small stable categories.
We refer it as the precosheaf of \textit{wrapped sheaves} associated to $\Lambda$.
\end{definition}

Note also that this construction is contravariant on the closed set $\Lambda$.
That is, if $\Lambda \subseteq \Lambda^\prime$ is an inclusion of closed subset of $S^* M$,
there is a canonical map $\wsh_{\Lambda^\prime}(U) \rightarrow \wsh_\Lambda(U)$ 
for $U \subseteq M$ by a similar consideration.
In other words, there is a morphism $\wsh_{\Lambda^\prime} \rightarrow \wsh_\Lambda$ between precosheaves.

\begin{conjecture}
Consider the case when $\Lambda$ is a subanalytic singular isotropic.
Inspired by homological mirror symmetry, Nadler defines in \cite{Nadler-pants} a conic cosheaf
$\operatorname{\msh}_\Lambda^w: \Op_{T^* M} \rightarrow \st$
through a purely categorical construction and terms it as the cosheaf of \textit{wrapped microlocal sheaves}.
One main property of $\operatorname{\msh}_\Lambda^w$ is that its restriction to the zero section, the `wrapped sheaves', 
is the cosheaf $\Sh_\Lambda^c$ discussed in Proposition \ref{shciscsh}.
We will reserve the term `wrapped sheaves' for the geometrically constructed category $\wsh_\Lambda$ through this paper.
Corollary \ref{maincor} of the main theorem asserts that these two cosheaves are the same after all.
As a result, we expect to extend the construction $\wsh_\Lambda$ to the cotangent bundle as well.
\end{conjecture}

\subsection{Generation}

We find a set of generators of $\wsh_\Lambda(M)$ when $\Lambda$ is a singular isotropic. 
We first prove a special case of the K\"unneth formula which we will refer to as the stabilization lemma.
We will be going back and forth between the homogeneous symplectic notation and the contact notation,
so for a conic set $X$ in $T^* M$, we will use the notation $X^\infty \coloneqq (X \setminus 0_M)/\Rp$
to denote its boundary at infinity.
The corresponding lemma in the Fukaya setting can be found in \cite[(1.7), (1.8)]{Ganatra-Pardon-Shende2}.
Fix $n \in \mathbb{R}^n$.
Let $M$ be a real analytic manifold and $\Lambda \subseteq S^* M$ be a closed subset. 
We set 
$$\Lambda^{st} = \left( ( \Rp \Lambda \cup 0_M) \times 0_{\mathbb{R}^n} \right)^\infty
\subseteq S^* (M \times \mathbb{R}^n).$$

Pick a small ball $B \subseteq \mathbb{R}^n$ centered at $0$.
For $F \in \Sh(M)$, by (6) of Proposition \ref{mses}, there is microsupport estimation
$$ \ms(F \boxtimes 1_B) \subseteq \ms(F) \times N^*_{out}(B).$$
As a result, exterior tensoring with $1_B$ induces a functor $$- \boxtimes 1_B:
\widetilde{\wsh}_\Lambda (M) \rightarrow \widetilde{\wsh}_{ \Lambda^{st} } \left(  M \times \mathbb{R}^n \right).$$
We claim that this functor induces a fully faithful functor on the quotient.
We first recall a lemma.

\begin{lemma}
Let $\sC$, $\sD$ be stable categories, 
$S$ and $T$ be sets of morphisms in $\sC$ and $\sD$
which are closed under composition and contain identities.
Set $\sC_0 \coloneqq \langle \cof(s) | s \in S \rangle$ and $\sD_0 \coloneqq \langle \cof(t) | t \in T \rangle$.
Let $F: \sC \rightarrow \sD$ be a functor such that for all $X_0 \xrightarrow{s} X_1 \in S$, 
there exists $F(X_1) \xrightarrow{t} Y \in T$ such that
$t \circ F(s) \in T$. 
Then $F|_{\sC_0}$ factors through $\sD_0$ and $F$ descends to a functor
$\bar{F}: (\sC/ \sC_0) \rightarrow (\sD/\sD_0)$ filling the commutative square
with the quotient functors:

$$
\begin{tikzpicture}

\node at (0,2) {$\sC_0$};
\node at (4,2) {$\sC$};
\node at (8,2) {$\sC/ \sC_0$};
\node at (0,0) {$\sD_0$};
\node at (4,0) {$\sD$};
\node at (8,0) {$\sD/\sD_0$};

\draw [right hook-latex, thick] (0.5,2) -- (3.6,2) node [midway, above] {$ $};
\draw [->>, thick] (4.5,2) -- (7.4,2) node [midway, above] {$ $};
\draw [right hook-latex, thick] (0.5,0) -- (3.6,0) node [midway, above] {$ $};
\draw [->>, thick] (4.5,0) -- (7.4,0) node [midway, above] {$ $};


\draw [->, thick] (0,1.7) -- (0,0.3) node [midway, left] {$F|_{\sC_0}$};
\draw [->, thick] (4,1.7) -- (4,0.3) node [midway, left] {$F$};
\draw [->, thick] (8,1.7) -- (8,0.3) node [midway, left] {$\bar{F}$}; 

\end{tikzpicture}.
$$
\end{lemma}

\begin{proof}
Let $X_0 \xrightarrow{s} X_1$ and $F(X_1) \xrightarrow{t} Y \in T$ be as above.
By the Lemma \ref{stsl}, there exists a fiber sequence $F(\cof(s)) \rightarrow \cof(t \circ F(s)) \rightarrow \cof(t)$. 
\end{proof}

\begin{proposition}[The stabilization lemma]\label{stb}
The functor $- \boxtimes 1_B$ defined above
descends to a fully faithful functor on the quotient
$$\wsh_\Lambda(M) \hookrightarrow \wsh_{\Lambda^{st}}( M \times \mathbb{R}^n )$$
which we will refer to as \textit{the stabilization functor}.
\end{proposition}

\begin{proof}
Let $\Phi: \dT^* M \times I \rightarrow \dT^* M$ be a positive homogenous symplectic isotopy
and $H = \alpha(\Phi_* \partial_t)$ be its corresponding Hamiltonian.
We note that since $H$, in general, cannot be defined smoothly on the entire $T^* M$, it is not always possible to extend it smoothly
to $\dT^* (M \times \mathbb{R}^n)$ by setting the dependence on the second component to be constant.
Nevertheless, $H^2$ is defined smoothly on $T^* M$ since $H$ is homogeneous of degree 1 by Lemma \ref{homogeneous_hamiltonian_extends}.\footnote{The author learns this trick from \cite[section 2]{Nadler-Shende}.}
Pick a bump function $\rho$ on $S^* (M \times \mathbb{R}^n)$
such that $ (\supp(H) \times \overline{T^* B})^\infty \subseteq  \Int( \supp(\rho) )$ so that 
$H^2 + \rho \abs{\xi}^2 > 0$ when $H > 0$.
Here we use the same notation $\rho$ to denote its pullback on $\dT^* M$.
Then set $\tilde{H} \coloneqq \sqrt{ H^2 + \rho \abs{\xi}^2}$ and we denote its corresponding homogeneous 
isotopy on $\dT^* (M \times \mathbb{R}^n)$ by $\tilde{\Phi}$.

The above construction implies that 
$$
\ms\left( (K(\Phi^{-1}) \boxtimes_I 1_{\Delta_{\mathbb{R}^n} \times I} ) \circ_I K( \tilde{\Phi} ) \right)
\subseteq \ms(K(\Phi^{-1}) \boxtimes_I 1_{\Delta_{\mathbb{R}^n} \times I}) \circ_I \ms( K( \tilde{\Phi} ) ) 
\subseteq \{ \tau \leq 0\}
$$
since $H(x,\xi) \leq \tilde{H}(x,\xi,t,\tau)$.
Thus there is a continuation map 
$$1_{\Delta_{M \times \mathbb{R}^n} } \rightarrow 
(w(\Phi^{-1}) \boxtimes 1_{\Delta_{\mathbb{R}^n} } ) \circ w( \tilde{\Phi} )$$
or equivalently $$w(\Phi) \boxtimes 1_{\Delta_{\mathbb{R}^n}} \rightarrow w(\tilde{\Phi})$$ which precomposes with 
$1_{\Delta_{M \times \mathbb{R}^n} } \rightarrow w(\Phi) \boxtimes_I 1_{\Delta_{\mathbb{R}^n} } $
to $1_{\Delta_{M \times \mathbb{R}^n} } \rightarrow w(\tilde{\Phi})$.
Thus, the last lemma implies that the functor 
$$\wsh_\Lambda(M) \rightarrow \wsh_{\Lambda^{st}}(M \times \mathbb{R}^n)$$ 
is well-defined.

A similar argument implies that for any positive isotopy $\Psi$ on $\dT^* (M \times \mathbb{R}^n)$
with Hamiltonian $\tilde{H}$, there exists $H$ on $\dT^* M$ and $\rho$ on $\dT^* \mathbb{R}^n$ 
such that $\tilde{H} \leq \sqrt{ H^2 + \rho^2}$.
Thus wrappings coming from the product is cofinal. 
Since $B$ is contractible, there is an isomorphism 
$\Hom(F \boxtimes 1_B, G \boxtimes 1_{\tilde{B}}) = \Hom(F,G)$
for any larger ball $\tilde{B}$ in $\mathbb{R}^n$.
This implies

\begin{align*} 
\Hom_{\wsh_{\Lambda^{st}}( M \times \mathbb{R}^n )}(G \boxtimes 1_B, F \boxtimes 1_B) 
&= \clmi{\tilde{\Phi} \in W (S^* (M \times \mathbb{R}^n )  \setminus \Lambda^{st}  )} 
\Hom \left(G \boxtimes 1_B, (F \boxtimes 1_B)^{\tilde{\Phi}} \right) \\
&= \clmi{\tilde{\Phi} \in W (S^* (M \times \mathbb{R}^n )  \setminus \Lambda^{st}  ) } 
 \Hom \left(G \boxtimes 1_B, F^\Phi \boxtimes 1_{\tilde{B} }\right) \\
&= \clmi{\Phi \in W(S^* M \setminus \Lambda)} \Hom \left(G , F^\Phi \right) \\
&= \Hom_{\wsh_\Lambda(M)}(G,F).
\end{align*}

Thus, the stabilization functor
$\wsh_\Lambda(M) \hookrightarrow \wsh_{\Lambda^{st}}( M \times \mathbb{R}^n )$ is fully faithful.
\end{proof}


Now we show that the category $\wsh_\varnothing(M)$ is generated by one object.
More precisely, we say an open set $B \subseteq M$ is a ball if $B$ is relative compact, contractible and $\overline{B}$ is a closed disk.
Now let $B$ be a ball such that there exists an open chart $U$ containing $\overline{B}$.
Since all such balls are smoothly isotopic to each other inside $M$, 
lifting such isotopies implies that the object $1_B$,
where $B$ is a such ball, is independent of the choice of the exact ball.
In order to show that $1_B$ is a generator, we need a class of auxiliary objects.

\begin{definition}
We say that an open set $B \subseteq M$ is a stable ball if it is relative compact, contractible, 
and $\overline{B}$ has a smooth boundary in $M$.
\end{definition}

One can check that a stable ball is a ball up to a stabilization by the famous corollary of the h-cobordism theorem.
The following statements are Theorem 5.12 and Corollary 5.13 in \cite{Ganatra-Pardon-Shende3}.

\begin{theorem}
A stable ball of dimension $\geq 6$ with simply connected boundary is a ball.
\end{theorem}

\begin{corollary}
Let $M$ be a stable ball. Then $B \times I^k$ is a ball provided $\dim B + k \geq 6$ and $k \geq 1$.
\end{corollary}

\begin{proof}
This is implied by a combination of the van Kampen Theorem and the Poinc\'are duality for manifolds with boundary
$$H^{ \dim N - k }(N,\partial N) = H_k (N).$$
\end{proof}

\begin{lemma}\label{ballg}
Assume $M$ is connected.
The category $\wsh_\varnothing(M)$ is generated under finite colimits and retractions by $1_B$ for any small ball $B$. 
\end{lemma}

\begin{proof}
Let $F \in \wsh_\varnothing(M)$ be an object.
By Remark \ref{iso=sm}, we may assume $F$ is a sheaf with compact support, subanalytic isotropic microsupport, and perfect stalks.
By Proposition \ref{extri}, there is a Whitney triangulation $\mathscr{T}$ such that $F$ is $\mathscr{T}$-constructible.
Since $\Sh_{N^*_\infty \mathscr{T}}(M)^c = \Perf \mathscr{T}$ is generated under finite colimits and retractions 
by $1_{\str(t)}$ for $t \in \mathscr{T}$, we may assume $F = 1_{\str(t)}$.
We claimed that the object $1_{\str(t)}$ is isomorphic to $1_B$ for some small ball $B \subseteq \str(t)$.
Note the open set $\str(t)$ is relatively compact and contractible, however, $\overline{\str(r)}$ might not be a manifold with boundary
and modification needs to be made.

Apply the inward cornering construction in Definition \ref{incor} to $U = \str(t)$, 
we obtain a family of $\str(t)^{-\epsilon}$ depending smoothly on $\epsilon$.
When $\epsilon$ is small, the object $1_{\str(t)^{-\epsilon}}$ in $\wsh_\Lambda(M)$ 
is independent of $\epsilon$ so we abuse the notation and simply denote it by $1_{\str(t)^-}$.
As there is no stop restriction, the canonical map $1_{\str(t)^{-} } \rightarrow 1_{\str(t)}$
becomes an isomorphism in $\wsh_\varnothing(M)$ through the positive wrapping obtained by taking $\epsilon \rightarrow 0$.
The closure $\overline{\str(t)^{-}}$ is a manifold with corners when $\epsilon$ is small, i.e., the boundary $\overline{\str(t)^{-}}$
can be modified by the boundary of the inclusion $ [0,\infty)^k \times \mathbb{R}^{n - k} \subseteq \mathbb{R}^n$ 
for some $k \geq 1$.
By Example \ref{conems}, $\msif( 1_{(0,\infty)^k \times \mathbb{R}^{n - k}})$ is smooth   
and $1_{(0,\infty)^k \times \mathbb{R}^{n - k}}$ can be wrapped to some $1_{V_\delta}$ 
where $V_\delta \coloneqq \{ x \in \mathbb{R}^n | d(x,[0,\infty)^k \times \mathbb{R}^{n - k}) < \delta \}$ by the Reeb flow.
One can see from the local model that the boundary of $V_\delta$ is smooth for small $\delta$.
Thus, we may further replace $1_{\str(t)^{-} }$ 
by some $1_{U}$ such that $U$ is relative compact, contractible and has smooth boundary, i.e., a stable ball.

Finally pick a ball $B \subseteq U$ and consider the canonical morphism $1_B \rightarrow 1_U$ induced by the inclusion.
Apply the stabilization lemma for $n$ large and we see from the last corollary that we may assume $U$ to be a ball as well.
In this case, the canonical map $1_B \rightarrow 1_U$ coincide with the continuation map 
obtained by the standard Reeb flow and is an isomorphism.
Thus, the original map is an isomorphism in $\wsh_\varnothing(T^* M)$ and we see that $1_B$ generates.
\end{proof}

We assume for the rest of this section that $\Lambda$ is a singular isotropic.
To study generation for the general case, we need the following lemma to perform general position argument.

\begin{definition}[{\cite[Definition 1.6]{Ganatra-Pardon-Shende2}}]
Let $Y^{2n-1}$ be a contact manifold.
We say a set $\mathfrak{f}$ is mostly Legendrian if there is a decomposition
$\mathfrak{f} = \mathfrak{f}^{subcrit} \cup \mathfrak{f}^{crit} \subseteq Y$ for which $\mathfrak{f}^{subcrit}$
is closed and is contained in the smooth image of a second countable manifold of dimension $< n-1$,
and $\mathfrak{f}^{crit}$ is a Legendrian submanifold.
\end{definition}

We note that a closed singular isotropic $\Lambda$ is in particular mostly Legendrian.

\begin{lemma}[{\cite[Lemma 2.2 and Lemma 2.3]{Ganatra-Pardon-Shende2}}]\label{gpa}
Let $Y^{2n-1}$ be a contact manifold and $\mathfrak{f}$ be mostly Legendrian.
\begin{enumerate}
\item Let $\Lambda \subseteq Y$ be a compact Lagrangian. 
Then $\Lambda$ admits cofinal wrappings $\Lambda \rightsquigarrow \Lambda^w$ 
with $\Lambda^w$ disjoint from $\mathfrak{f}$.
\item Let $\Lambda_1$, $\Lambda_2 \subseteq Y$ be compact Legendrians disjoint from $\mathfrak{f}$.
Consider the space of positive Legendrian isotopies $\Lambda_1 \rightsquigarrow \Lambda_2$.
Then the subspace of isotopies which
\begin{enumerate}[label*=\arabic*.]
\item remains disjoint from $\mathfrak{f}^{subcrit}$ and
\item intersect $\mathfrak{f}^{crit}$ only finitely many times, each time passing transversally at a single point,
\end{enumerate}
is open and dense.
\end{enumerate}
\end{lemma}

For inclusion of singular isotropics $\Lambda \subseteq \Lambda^\prime \subseteq S^* M$, 
general position argument implies that the induced map
$\wsh_{\Lambda^\prime}(M) \rightarrow \wsh_{\Lambda}(M)$ is always essentially surjective.
In order to study the fiber of this map, we consider the following objects:

Let $\Lambda$ be a subanalytic singular isotropic and let $(x,\xi) \in \Rp \Lambda $ be a smooth point.
Consider a proper analytic $\Lambda$-Morse function $f: M \rightarrow \mathbb{R}$.
We assume there exists $\epsilon > 0$
such that $x$ is the only $\Lambda$-critical point over $f^{-1} ( [-\epsilon, \epsilon] )$
with critical value $0$, $d f_x = \xi$ and $f^{-1}(-\infty, \epsilon)$ is relatively compact.
By our assumption, both $1_{f^{-1}(-\infty,\pm \epsilon)}$ are objects of $\widetilde{\wsh}_\Lambda(M)$.

\begin{definition}
A sheaf-theoretical linking disk at $(x,\xi)$ (with respect to $\Lambda$) is an object $D_{(x,\xi)}$ of the form 
$$\cof ( 1_{f^{-1}(-\infty,- \epsilon)} \rightarrow 1_{f^{-1}(-\infty,\epsilon)} ) \in \wsh_\Lambda(M)$$
where the arrow is induced by the inclusion of opens $f^{-1}(-\infty,- \epsilon) \subseteq f^{-1}(-\infty,\epsilon)$
given by a function $f$ with the above properties.
Note that by scaling $f$ with $r \in \Rp$, we see that the object $D_{(x,\xi)}$ depends only on $(x,\xi)$'s image in $S^* M$.
Thus, we also use the same notation $D_{(x,\xi)}$ for $(x,\xi) \in S^* M$. 
\end{definition}

\begin{remark}\label{ldlp}
Since $d f \neq 0$ over $f^{-1} ( [-\epsilon, \epsilon] )$, the fibers $f^{-1}(t)$ for $t \in [-\epsilon,\epsilon]$ are smooth submanifolds.
Thus, the canonical map $ 1_{f^{-1}(-\infty,- \epsilon)} \rightarrow 1_{f^{-1}(-\infty,\epsilon)}$
is also given by the continuation map of the wrapping $\{ N^*_{\infty,out} f^{-1}(-\infty,t) \}_{t \in [-\epsilon,\epsilon]}$ which passes through
$\Lambda$ transversely exactly once at $(x,[\xi])$.
Extend this wrapping to a global one $\Psi$. 
Since there is no other intersection with $\Lambda$, we can decompose $\Psi$ to $\Psi_+ \# \Psi_0 \# \Psi_-$ so that $\Psi_\pm$
do not intersect $\Lambda$ and $\Psi_0$ only moves points near $(x,\xi)$.
This way, we can see $D_{(x,\xi)}$ can be presented as a cofiber induced by an expanding open half-plane.
$$\includegraphics[scale = 0.2]{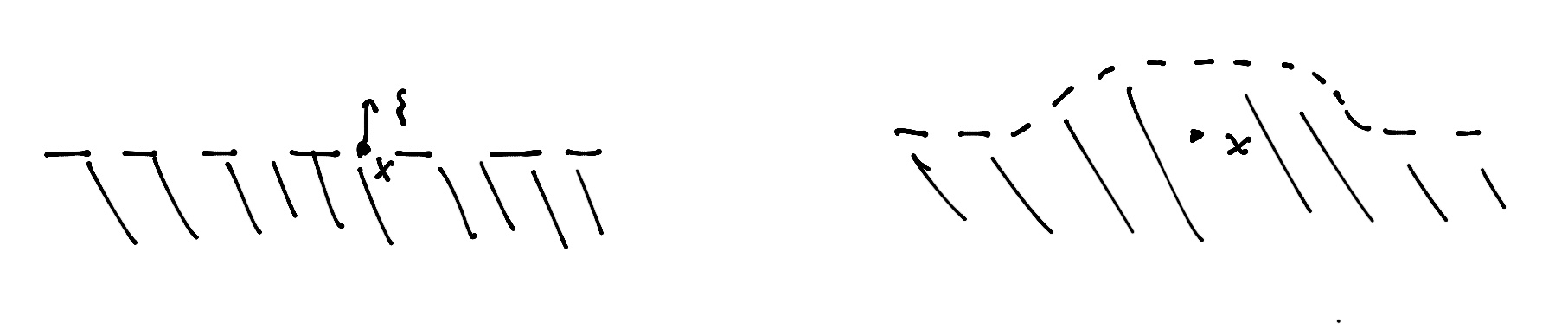}$$
Thus, $D_{(x,\xi)}$ can also be presented as a cofiber induced by inclusions of small balls.
$$\includegraphics[scale = 0.2]{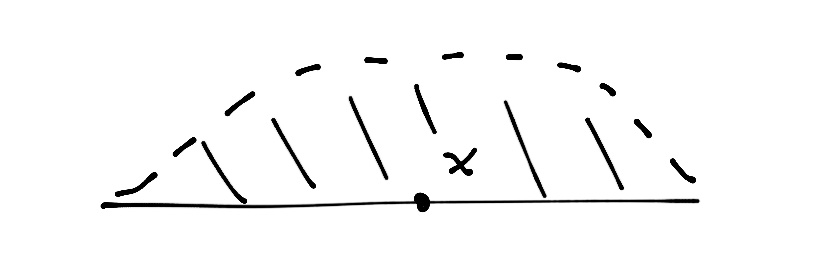}$$
\end{remark}

\begin{proposition}\label{wgstopr}
Let $\Lambda \subseteq \Lambda^\prime$ be subanalytic singular isotropics and
let $\sD^w_{\Lambda^\prime,\Lambda}( M)$ denote the fiber of the canonical map 
$\wsh_{\Lambda^\prime}(M) \rightarrow \wsh_{\Lambda}(M)$.
Then $\sD^w_{\Lambda^\prime,\Lambda}(M)$ 
is generated by the sheaf-theoretical linking disk $D_{(x,\xi)}$ for smooth Legendrian points 
$(x,\xi) \in \Lambda^\prime \setminus \Lambda$.
\end{proposition}

\begin{proof}
Let $F \in \wsh_\Lambda(M)$. 
We assume that $\msif(F)$ is a subanalytic isotropic and pick a Whitney triangulation $\mathscr{T}$ 
such that $\msif(F) \subseteq N^*_\infty \mathscr{T}$.
Fixed a particular way to construct $F$ out of sheaves of the form $M_{\str(t)}$ for some $M \in \cV_0$
by taking finite steps of cofibers and use $\{F_i\}_{ \{i \in A \} }$ to denote those $M_{\str(t)}$'s which show up in these steps.
Note that it is possible that their microsupport $\msif(F_i)$ intersect $\Lambda$.
However, we see from the proof of Lemma \ref{ballg}, 
the microsupport $\msif(F_i)$ of these $F_i$'s are smooth Legendrians in $S^* M$.
Thus, we can apply Lemma \ref{gpa} and assume $\msif(F_i) \cap \Lambda^\prime = \varnothing$ for $i \in A$. 
By the microsupport triangular inequality (1) of Proposition \ref{mses},
the sheaves appear in the the cofiber sequences which build $F$ from these $F_i$'s
do not intersect $\Lambda^\prime$ as well. 
In particular, $\msif(F) \cap \Lambda^\prime = \varnothing$.
Similarly, an application of Lemma \ref{gpa} implies that we can assume the existence of
a cofinal wrapping sequence $F \rightarrow F^{w_1} \rightarrow F^{w_2} \rightarrow \cdots$ 
such that $\msif(F^{w_n}) \cap \Lambda^\prime = \varnothing$.
This implies that the canonical map 
$$\widetilde{\wsh}_{\Lambda^\prime}(M)/
\left( \sC_\Lambda(M) \cap \widetilde{\wsh}_{\Lambda^\prime}(M) \right)
\rightarrow \wsh_\Lambda(M) \coloneqq \widetilde{\wsh}_\Lambda(M) / \sC_\Lambda(M)$$
is an equivalence.
Thus, we can apply Lemma \ref{prlsn} to the diagram

$$
\begin{tikzpicture}

\node at (0,2) {$\sC_{\Lambda^\prime}(M)$};
\node at (6,2) {$\widetilde{\wsh}_{\Lambda^\prime}(M)$};
\node at (12,2) {$\wsh_{\Lambda^\prime}(M)$};
\node at (0,0) {$\sC_\Lambda(M) \cap \widetilde{\wsh}_{\Lambda^\prime}(M)$};
\node at (6,0) {$\widetilde{\wsh}_{\Lambda^\prime}(M)$};
\node at (12,0) {$\wsh_\Lambda(M)$};

\draw [right hook-latex, thick] (1,2) -- (4.7,2) node [midway, above] {$i$};
\draw [->>, thick] (7.3,2) -- (10.8,2) node [midway, above] {$p$};
\draw [right hook-latex, thick] (1.9,0) -- (4.7,0) node [midway, above] {$j$};
\draw [->>, thick] (7.3,0) -- (10.8,0) node [midway, above]  {$q$};

\draw [right hook-latex, thick] (0,1.7) -- (0,0.3) node [midway,right] {$ $};
\draw [double equal sign distance, thick] (6,1.7) -- (6,0.3) node [midway, right] {$ $};
\draw [->>, thick] (12,1.7) -- (12,0.3) node [midway,right] {$ $};

\end{tikzpicture}
$$
which implies that  $\sD^w_{\Lambda^\prime,\Lambda}(M) 
= \left( \sC_\Lambda(M) \cap \widetilde{\wsh}_{\Lambda^\prime}(M) \right) /\sC_{\Lambda^\prime}(M)$.
Now $\sC_\Lambda(M)$ is the category generated by the cofibers $\cof( c(\Psi,F) )$ of continuation maps
whose wrapping $\psi_t( \msif(F))$ avoids $\Lambda$, and 
$ \sC_\Lambda(M) \cap \widetilde{\wsh}_{\Lambda^\prime}(M)$ 
is generated by a similar construction except we now only requires the end points to avoid $\Lambda$.
The claim is that the quotient is generated by the sheaf-theoretic linking disks 
$D_{(x,\xi)}$ for smooth Legendrian points $(x,\xi) \in \Lambda^\prime \setminus \Lambda$.

Lemma \ref{stsl} implies that if $H_1 \rightarrow H_2 \rightarrow H_3$ is a cofiber sequence,
then $\cof(c(\Psi,H_1)) \rightarrow \cof( c(\Psi,H_2)) \rightarrow \cof( (c(\Psi,H_3))$ is also a fiber sequence.
Thus, it is enough to assume $F$ has smooth Legendrian microsupport by the discussion at the beginning of the proof.
Let $\Psi$ be a positive isotopy such that  $\psi_t( \msif(F))$ does not touch $\Lambda$ and, by general position argument, 
we may assume $\psi_t(\msif(F))$ touches $\Lambda^\prime$ for finitely many times
and transversally through one point $p \in \Lambda^\prime \setminus \Lambda$ each time.
Decomposing $\Psi$ to $\Psi = \Psi_k \# \cdots \# \Psi_1$ so that passing happens only once during the duration of each $\Psi_i$.
Since $c(\Psi) = c(\Psi_k) \circ \cdots \circ c(\Psi_1)$, 
it is sufficient to prove the case when only one such passing at $(x,\xi)$ appears by induction with Lemma \ref{stsl}.

Let $q \in \msif(F)$ be the point so that the path $\psi_t(q)$ pass $(x,\xi) \in \Lambda^\prime$
and $U$ small open ball near $q$ in $S^* M$. 
We again decompose $\Psi$ to $\Psi = \Psi_+ \# \Psi_0 \# \Psi_-$ such that 
there's no passing happing during $\Psi_\pm$ and $\Psi_0$ only moves points in $U$.
In this case, $c(\Psi_\pm,F)$ are isomorphisms and we can further assume $\Psi$ only moves points in $  U$.
Now set $F_t = (K(\Psi) \circ F )|_{M  \times \{t\}} $ so $\msif(F_t) = \psi_t \left( \msif(F) \right)$.
We use one last general position argument to assume the front projection $\pi_\infty: \msif(F_t) \rightarrow \pi_\infty( \msif(F_t) )$ 
is finite near $\psi_t(U)$ so $\pi(U \cap \msif(F_t)) \subseteq \pi(U)$ is a hyperplane.
Thus, we reduce to the local picture defining $D_{(x,\xi)}$ discussed in the last Remark \ref{ldlp}.
\end{proof}

We combine the above two results to deduce a generation result for a special case.
Let $\mathscr{S}$ be a Whitney triangulation.
For each stratum $s \in \mathscr{S}$, we pick a small ball $B_s$ which centered at $X_s$ and contained in $\str(s)$
such that $N^*_{\infty,out} B_s \cap N^*_\infty \mathscr{S} = \varnothing$
and consider $1_{B_s} \in \wsh_{N^*_{\infty,out} \mathscr{S}}(M)$.
This is possible because of the Whitney condition.
Again different choices of such small balls induce the same objects in $\wsh_{N^*_\infty \mathscr{S}}(M)$
since they are isotopic to each other in the base by isotopies respecting the stratification
and the lifting isotopies on the microsupport won't touch $N_\infty^*(\mathscr{S})$.
(See \cite{Mather}.)

\begin{proposition}\label{strbg}
The set $\{ 1_{B_s} \}_{ s \in \mathscr{S} }$ generates 
$\wsh_{N^*_\infty \mathscr{S}  }(M)$ under finite colimits and retractions.
\end{proposition}

\begin{proof}
Set $\mathscr{S}_{\leq k} = \{s \in \mathscr{S} | \dim X_s \leq k \}$.
We claim, when $k < n -1$,
$\{ 1_{B_s} \}_{ s \in \mathscr{S}_{\leq k} }$ plus $1_B$ for any small ball $B$ whose closure
$\overline{B}$ is disjoint from any stratum of dimension $\leq k$ generates $\wsh_{N_\infty^* \mathscr{S}_{\leq k} }(M)$.
To see this, we note that the case $k = -1$ is Lemma \ref{ballg}.
Assume the case for $k < n - 2 $ and consider the projection $\wsh_{N_\infty^* \mathscr{S}_{\leq k} } (M) 
\twoheadrightarrow \wsh_{N_\infty^* \mathscr{S}_{\leq k-1}} (M)$. 
We note that $M \setminus \cup_{s \in \mathscr{S}_{k + 1} } X_s$ is path connected by standard Hausdorff dimension theory
since these $X_s$'s have codimension $\geq 2$.
Thus, $1_B$ is independent of the choice of $B$.
By Proposition \ref{wgstopr}, the fiber of the above projection is generated under finite colimits and retracts 
by sheaf-theoretic linking disks $D_{(x,\xi)}$ 
for $(x,\xi) \in N^*_\infty \mathscr{S}_{\leq k+1 } \setminus N^*_\infty \mathscr{S}_{\leq k} $. 
But $D_{(x,\xi)}$ can be written as the cofiber $\cof(1_B \rightarrow 1_{B_s} )$ 
by the local picture mentioned in Remark \ref{ldlp}.
Finally, apply a similar argument to the projection $ \wsh_{N_\infty^* \mathscr{S} } (M) 
= \wsh_{N_\infty^* \mathscr{S}_{\leq n-1} } (M) \rightarrow \wsh_{N_\infty^* \mathscr{S}_{\leq n-2}} (M)$
implies the proposition.
\end{proof}

Recall in the proof of Lemma \ref{ballg}, we show that $1_{B_s} \rightarrow 1_{\str(s)}$ is an isomorphism in
$\wsh_\varnothing(M)$.
The later object is, however, not an object in $\wsh_{N^*_\infty \mathscr{S}}(M)$.
Instead, we consider the object $1_{\str(s)^-}$ where $\str(s)^-$ is a small inward cornering of $\str(s)$ of Definition \ref{incor}.
Choose $B_s$ small so that $B_s \subseteq \str(s)^-$.
We claim that the canonical map $1_{B_s} \rightarrow 1_{\str(s)^-}$ is an isomorphism.
By the Yoneda embedding, it is an isomorphism if the corresponding morphism 
$\Hom_w(-,1_{B_s}) \rightarrow \Hom_w(-,1_{\str(s)^-})$ is an isomorphism
as presheaves on $\wsh_{N^*_\infty \mathscr{S}} (M)$.
The following statements are in directly parallel with Proposition 5.18, Lemma 5.21, and Proposition 5.24 in \cite{Ganatra-Pardon-Shende3}.

\begin{lemma}
For a $\mathscr{S}$-constructible relatively compact open set $U$, we have
$$\Hom_w(1_{B_s},1_{U^{-}}) = 
\begin{cases}
1 &\str(s) \subseteq U \\
0 &\text{otherwise}
\end{cases}
$$
\end{lemma}

\begin{proof}
The construction of $U^-$ tautologically provides a cofinal sequence $1_{U^{-\epsilon} }$ so
$$ \Hom_w(1_{B_s},1_{U^- }) =
 \clmi{\epsilon \rightarrow 0}   \Hom(1_{B_s}, 1_{U^{-\epsilon} }).$$
First consider the case when $\str(s) \subseteq U$.
Since $B_s \subseteq \str(s)$ has a non-zero distance from $\partial U$, it is contained in $U^{-\epsilon}$ for $\epsilon << 1$
and the left hand side is $1$. 
When $\str(s) \cap U = \varnothing$, the $\Hom$ is clearly $0$ so we assume $s$ is a stratum on the boundary of $U$.
In this case, one can conclude the result by refining the wrapping to a family $U^{-\epsilon,-\delta}$
where we add the center of the ball $B_s$ to the stratification and $\delta$ denotes the parameter
which corresponds to this new stratum. See \cite[Proposition 2.10]{Ganatra-Pardon-Shende3}.
\end{proof}

\begin{proposition}\label{b=str}
The canonical map $1_{B_s} \rightarrow 1_{\str(s)^-}$ is an isomorphism in $\wsh_{N^*_\infty \mathscr{S} }(M)$.
\end{proposition}

\begin{proof}
We proceed by induction on the codimension of $s$. 
When $s$ has codimension zero, we may replace $M$ by $\str(s)$ and it becomes Lemma \ref{ballg}.

Now the previous lemma and the proposition \ref{strbg} implies
$$\Hom(1_{B_s},1_{{\str(t)}^- }) = \Hom(1_{\str(s)^-},1_{{\str(t)}^- }) = 0$$
for $t$ of strictly smaller codimension than $s$.
By induction, $1_{B_t} \xrightarrow{\sim} 1_{\str(t)^-}$ for such $t$'s.
The later generates a subcategory which contains the fiber of the projection 
$\wsh_{N_\infty^* \mathscr{S}}(T^* M) \rightarrow \wsh_{N_\infty^* \mathscr{S}_{\leq \dim s} }(T^* M)$.
This implies that it is enough to show the isomorphism in the category $\wsh_{N_\infty^* \mathscr{S}_{\leq \dim s} }(T^* M)$.
This is a special case of the following lemma applying to the case $Y = \str(s)^-$, $X = s \cap \str(s)^-$,
and $Z = t \cap \str(s)^-$.
\end{proof}

\begin{lemma}
Let $X^m \subseteq Y^n$ be an inclusion of stable balls, with $\partial X \subseteq \partial Y$. Assume there exists
another stable ball (with corners) $Z^{m+1} \subseteq Y^n$ such that $\partial Z$ is the union of $X$ with a smooth
submanifold of $\partial Y$. 
Then the canonical map $1_{B_\epsilon(x)} \rightarrow 1_Y$ 
is an isomorphism in $\wsh_{N_\infty^* X}(Y)$ for any $x \in X$.
\end{lemma}

\begin{proof}
Reduce to the case of balls by stabilization.
In this case, $Y$ is a unit ball, $X$ is the intersection of $Y$ with a linear subspace, 
and $Z$ is the intersection of $Y$ with a closed half-plane with the boundary being the linear subspace.
The positive isotopy which expands $1_{B_\epsilon(x)}$ to $1_Y$ is disjoint from $N_\infty^* X$.
\end{proof}

\begin{corollary}
The set $\{ 1_{\str(s)^-} \}_{ s \in \mathscr{S} }$ generates 
$\wsh_{N^*_\infty \mathscr{S}  }(M)$ under finite colimits and retractions.
\end{corollary}

\section{The comparison morphism}
Let $M$ be a real analytic manifold and $\Lambda \subseteq S^* M$ a subanalytic singular isotropic.
We define in this section, by an abuse of notation, a \textit{comparison} functor 
$$\wrap_\Lambda^+(M): \wsh_\Lambda(M) \rightarrow \Sh_\Lambda(M)^c$$
and show that it is an equivalence of category.
Since such functors combine to a \textit{comparison} morphism 
$\wrap_\Lambda^+ : \wsh_\Lambda \rightarrow \Sh_\Lambda^c$ between
precosheaves, the last statement will imply that $\wsh_\Lambda$ is in particular a cosheaf
for this case.

\subsection{Definition}
Let $\Lambda \subseteq S^* M$ be a closed singular isotropic subset.
Recall from Proposition \ref{w=ad} that
the inclusion $\Sh_\Lambda(M) \hookrightarrow \Sh(M)$ has a left adjoint given by the positive infinite wrapping functor 
$$\wrap_\Lambda^+(M): \Sh(M) \rightarrow \Sh_\Lambda(M).$$
Geometrically, it takes a sheaf $F$ to the limiting object over increasingly positive wrappings.
Since a continuation map $ c: F \rightarrow F^w$
tautologically becomes an isomorphism after applying $\wrap_\Lambda^+(M)$,
the functor $\wrap_\Lambda^+(M)$ vanishes on $\sC_\Lambda(M)$.
\hypertarget{infwrapfun}{We} abuse the notation and denote the resulting functor on the quotient category also by
$$ \wrap_\Lambda^+(M): \wsh_\Lambda(M) \rightarrow \Sh_\Lambda(M).$$

\begin{remark}
We note that in general when $\Lambda$ is not a singular isotropic, the category on the right hand side is much larger.
For example, when $\Lambda = S^* M$ and $M$ is non-compact, 
$\wrap_{S^* M}^+(M)$ is the trivial inclusion $\{0\} \hookrightarrow \Sh(M)$.
\end{remark}
We first notice that in this case the restriction of $\wrap_\Lambda^+(M)$ on $\wsh_\Lambda(M)$ takes image in 
the subcategory consisting of compact objects.

\begin{lemma}
Let $\Lambda$ be a subanalytic singular isotropic.
For $F \in \wsh_\Lambda(M)$, the sheaf $\wrap_\Lambda^+(M) (F)$ is a compact object.
\end{lemma}

\begin{proof}
Let $F \in \wsh_\Lambda(M)$ and $\varinjlim F_i$ be a filtered colimit in $\Sh_\Lambda(M)$. We compute, 
\begin{align*}
\Hom(\wrap_\Lambda^+(M) (F), \varinjlim F_i) 
&= \clmi{\Phi \in W(S^* M \setminus \Lambda)} \Hom(w(\Phi) \circ F, \varinjlim F_i) \\
&= \clmi{\Phi \in W(S^* M \setminus \Lambda)} \Hom(F, w(\Phi^{-1}) \circ \varinjlim F_i) \\
&= \clmi{\Phi \in W(S^* M \setminus \Lambda)} \Hom\left(F,  \varinjlim  (w(\Phi^{-1}) \circ F_i) \right) \\
&= \clmi{\Phi \in W(S^* M \setminus \Lambda)} \Hom (F, \varinjlim F_i)  = \Hom(F,\varinjlim F_i).
\end{align*}

For the last equality, we use the fact that $\Phi$ is supported away from 
$\Lambda \supseteq \msif(F_i)$ so $w(\Phi^{-1}) \circ F_i = F_i$ by Lemma \ref{loclem}.
Now pick a Whitney triangulation $\mathscr{S}$ such that $F$ is $\mathscr{S}$-constructible
and $\Lambda \subseteq N^* \mathscr{S}$.
In this case, the $\Hom$ can be computed in $\Sh_{N^*_\infty \mathscr{S}}(M) = \mathscr{S} \dMod $.
Since $\Sh_{N^*_\infty \mathscr{S}}(M)^c$ consists exactly objects with compact support and perfect stalks, 
$F$ is compact in $\Sh_{N^*_\infty \mathscr{S}}(M)$.
Thus $\Hom(F,\varinjlim F_i) = \varinjlim \Hom(F,\ F_i) $ and a backward computation as above implies that
$$\Hom(\wrap_\Lambda^+(M) (F), \varinjlim F_i)  =  \varinjlim \Hom(\wrap_\Lambda^+(M) (F), F_i) $$
so $\wrap_\Lambda^+(M) (F) \in \Sh_\Lambda(M)^c$ is compact.
\end{proof}

We note that this map is compatible with the precosheaf structure on both side.

\begin{lemma}
Let $j:U \subseteq M$ be an open set. 
The restriction $j^*: \Sh_\Lambda(M) \rightarrow \Sh_{\Lambda|_U}(U)$ has left and right adjoints which are given by
$\wrap_\Lambda^+\circ j_!$ and $\wrap_\Lambda^- \circ j_*$.
Hence, taking left adjoint induces a functor $\wrap_\Lambda^+  \circ j_!
:\Sh_{\Lambda|_U}(U)^c \rightarrow \Sh_\Lambda(M)^c$ between compact objects.
\end{lemma}

\begin{proof}
We use the fact that the left adjoint of a left adjoint preserves compact objects.
\end{proof}

Note when $\Omega \subseteq \Omega^\prime$, there is equivalence $\wrap^+(\Omega^\prime) \circ \wrap^+(\Omega) 
= \wrap^+(\Omega) \circ \wrap^+(\Omega^\prime) = \wrap^+(\Omega^\prime)$.
Thus, by the above lemma, there is commuting diagram for an inclusion of opens $j: U  \hookrightarrow V$:
$$
\begin{tikzpicture}

\node at (0,2) {$\wsh_\Lambda(U)$};
\node at (5.5,2) {$\Sh_{\Lambda|_U}(U)^c$};
\node at (0,0) {$\wsh_\Lambda(V)$};
\node at (5.5,0) {$\Sh_{\Lambda|_V}(V)^c$};

\draw [->, thick] (0.9,2) -- (4.4,2) node [midway, above] {$\wrap_{\Lambda |_U}^+(U)$};
\draw [->, thick] (0.9,0) -- (4.4,0) node [midway, above] {$\wrap_{\Lambda |_V}^+(V)$};


\draw [->, thick] (0,1.7) -- (0,0.3) node [midway, left] {$j_!$};
\draw [->, thick] (5.5,1.7) -- (5.5,0.3) node [midway, right] {$\wrap_{\Lambda|_V}^+(V) \circ j_!$};
\end{tikzpicture}
$$

\begin{definition}\label{infwrapfun}
We will refer the morphism $\wrap_\Lambda^+: \wsh_\Lambda \rightarrow \Sh_\Lambda^c$ between precosheaves
defined by the above diagram as the \textit{comparison} morphism.
\end{definition}

Similarly, when $\Lambda \subseteq \Lambda^\prime$, recall the left adjoint of the inclusion 
$\Sh_\Lambda(M) \hookrightarrow \Sh_{\Lambda^\prime}(M)$ is given by
$\wrap_\Lambda^+(M)$ and thus there is a commuting diagram:

$$
\begin{tikzpicture}

\node at (0,2) {$\wsh_{\Lambda^\prime}(M)$};
\node at (5.5,2) {$\Sh_{\Lambda^\prime}(M)^c$};
\node at (0,0) {$\wsh_\Lambda(M)$};
\node at (5.5,0) {$\Sh_{\Lambda}(M)^c$};

\draw [->, thick] (1,2) -- (4.5,2) node [midway, above] {$\wrap_{\Lambda^\prime}^+(M) $};
\draw [->, thick] (1,0) -- (4.5,0) node [midway, above] {$\wrap_{\Lambda^\prime}^+(M) $};


\draw [->, thick] (0,1.7) -- (0,0.3) node [midway, left] {$ $};
\draw [->, thick] (5.5,1.7) -- (5.5,0.3) node [midway, right] {$\wrap_\Lambda^+(M)$};

\end{tikzpicture}
$$

One can see this is compatible with the corestrictions on both side.
Thus, there is a commuting diagram in precosheaves with coefficient in $\st_w$:

$$
\begin{tikzpicture}

\node at (0,2) {$\wsh_{\Lambda^\prime}$};
\node at (5,2) {$\Sh_{\Lambda^\prime}^c$};
\node at (0,0) {$\wsh_\Lambda$};
\node at (5,0) {$\Sh_{\Lambda}^c$};

\draw [->, thick] (0.7,2) -- (4.4,2) node [midway, above] {$\wrap_\Lambda^+ $};
\draw [->, thick] (0.7,0) -- (4.4,0) node [midway, above] {$\wrap_{\Lambda^\prime}^+ $};


\draw [->, thick] (0,1.7) -- (0,0.3) node [midway, left] {$ $};
\draw [->, thick] (5,1.7) -- (5,0.3) node [midway, left] {$ $};

\end{tikzpicture}.
$$

The main \hyperlink{main}{theorem} of this paper, Theorem \ref{main}, is that the comparison functor
$$\wrap_\Lambda^+(M): \wsh_\Lambda(M) \rightarrow \Sh_\Lambda(M)$$ is an equivalence.
As a corollary, the comparison morphism $$\wrap_\Lambda^+: \wsh_\Lambda \rightarrow \Sh_\Lambda^c$$
is an isomorphism so, in particular, $\wsh_\Lambda$ is a cosheaf since $\Sh_\Lambda^c$ is.

\subsection{Sufficient condition for fully faithfulness}

For the rest of the section, we work with a fixed pair $(M,\Lambda)$ such that $\Lambda \subseteq S^* M$
is a subanalytic singular isotropic.
We would like to study the effect of $\wrap_\Lambda^+$ on the $\Hom$.
Since $\wrap_\Lambda^+$ is defined by a colimit,
the canonical $\Hom_w(G,F) \rightarrow \Hom(\wrap_\Lambda^+ G, \wrap_\Lambda^+ F)$ 
can be obtained from the following few steps. 
By definition of colimits, there is a canonical map \textit{colimiting continuation map} $F^w \rightarrow \wrap_\Lambda^+ F$ 
for any wrapping $w$.
This induces, for any other wrapping $w^\prime$, a map between the $\Hom$'s
$\Hom(G^{w^\prime},F^w) \rightarrow \Hom(G^{w^\prime},\wrap_\Lambda^+ F)$.
Since convolving with $w(\Phi)$ is an auto-equivalence on $\Sh(M)$,
there is a canonical map 
$$\Hom(G,F^w) = \Hom(G^{w^\prime},(F^w)^{w^\prime}) \rightarrow \Hom(G^{w^\prime},\wrap_\Lambda^+(F)).$$
Take limit over $w^\prime$ and then colimit over $w$, we obtain the map between $\Hom$'s
$$\Hom_w(G,F) = \clmi{w} \Hom(G,F^w) \rightarrow \lmi{w^\prime} \Hom(G^{w^\prime},\wrap_\Lambda^+(F))
= \Hom(\wrap_\Lambda^+(G),\wrap_\Lambda^+(F)).$$
In short, we have the following lemma.

\begin{lemma}\label{samehom}
Running $F$ and $G$ through a set of generators of $\wsh_\Lambda(M)$.
If the limiting continuation map $F \rightarrow \wrap_\Lambda^+ F$ becomes an isomorphism after applying $\Hom(G,-)$ 
for all such $G$, then the canonical map on $\Hom$
$\Hom_w(G,F) \rightarrow \Hom(\wrap_\Lambda^+(G),\wrap_\Lambda^+(F))$
is an isomorphism for any $F, G \in \wsh_\Lambda(M)$.
\end{lemma}

Pick any cofinal functor $\Psi: \ZZ_{\geq 0} \rightarrow W(S^* M \setminus \Lambda)$ 
which corresponds to a sequence of wrappings $\id \xrightarrow{\Psi_0} \Phi_1 \xrightarrow{\Psi_1} \Phi_2 \rightarrow \cdots$.
For convenience, we scale it so that $\Psi_i$ has domain $S^* M \times [i,i+1]$.
This sequence of positive family of isotopies patches to a positive isotopy
$\Psi: S^* M \times [0,\infty) \rightarrow S^* M$
whose restriction on $S^* M \times [i,i+1]$ is $\Psi_i \# \Phi_i$.
Note that $\Psi$ has a non-compact support by the cofinal criterion Lemma \ref{cfnc}.
By the GKS sheaf quantization, there is a sheaf kernel $K(\Psi)$ on $M \times M \times [0,\infty)$
such that $K(\Psi)|_{M \times M \times [i,i+1]} = K(\Psi_i)$.
For $F \in \widetilde{\wsh}_\Lambda(M)$,  we write $F^\Psi = K(\Psi) \circ F$ and
let $F^{w_n}$ denote $F^{w(\Phi_n)} = F^\Psi |_{M \times \{n\}}$ the resulting sheaves under the wrapping $\Psi$.
it is enough to study the morphism $$\Hom(G,F^{w_n}) \rightarrow \Hom(G,\wrap_\Lambda^+(F) )$$
which is induced from the sequence of wrappings $$F \rightarrow F^{w_1} \rightarrow \cdots \rightarrow
F^{w_n} \rightarrow \cdots \rightarrow \wrap_\Lambda^+(F).$$

\begin{definition}
Let $X$ be a topological space, $j:X \times \RR \hookrightarrow X \times  (-\infty,\infty]$
and $i:X \times \{\infty\} \rightarrow X \times (-\infty,\infty]$ be the inclusions as open and closed subset.
We call the composition $\psi = i^* \circ j_*: \Sh(X \times (-\infty,\infty) ) \rightarrow \Sh(X)$ the \textit{nearby cycle} functor.
\end{definition}

\begin{lemma}\label{w=ps}
The colimit $\wrap_\Lambda^+ F$ can be computed as the nearby cycle at infinity of the sheaf $F^\Psi$.
That is, $\wrap_\Lambda^+ F = \psi F^\Psi$.
\end{lemma}

\begin{proof}
Since $\Psi_n$ are cofinal, $\wrap_\Lambda^+(F) = \clmi{n \in \mathbb{N}} \ F^{w_n}$.
By the construction above, for each $n>1$, $F^{w_n}$ is given by 
$1_{\{n\}} \circ F^{\Psi} \xrightarrow{\sim}  1_{(0,n)}[1] \circ F^\Psi$
and the continuation map between them is induced by $1_{(0,n)} \rightarrow 1_{(0,m)}$ for $m \geq n$.
Since convolution commutes with colimit, $\wrap_\Lambda^+(F) 
= (\clmi{n \in \mathbb{N}}\ 1_{(0,n)}[1]) \circ F^\Psi = 1_{(0,\infty)}[1] \circ F^\Psi = p_! F^\Psi[1]$
where $p: M \times [0,\infty) \rightarrow M$ is the projection.
The latter is the same as $\psi F^\Psi = i^* j_* F^\Psi$ by Lemma \ref{ps=nc}
because $\Psi$ is a positive isotopy and so $ \ms(F^\Psi) \subseteq \{ \tau \leq 0 \}$.
\end{proof}

To study the (co)limiting continuation map $\Hom(G,F) \rightarrow \Hom(G,\psi F^\Psi)$,
we use a similar trick as in Proposition \ref{pert} and consider the object $\sHom(p^* G, F^\Psi)$
where we use $p: M \times \RR \rightarrow M$ to denote the projection. 
This time, we have to study its behavior near the infinity.

\begin{lemma}\label{homtinf}
Let $F, G \in \Sh(M \times \RR)$ be sheaves on $M \times \RR$ such that $F$ and $\sHom(G,F)$ are $\RR$-noncharacteristic,
and $q$ is proper on $\supp(G)$ and $\supp(\psi G)$ is compact.
Then $\Hom(i_s^* G, i_s^* F)$ is contant on $s \in \RR$ and equals to 
$\Gamma \left(M; \psi \sHom(G,F) \right)$
\end{lemma}

\begin{proof}
The statement over $\RR$ follows from Lemma \ref{ncpro}.
Since we assume $\supp(\psi(G))$ is proper, we can apply base change on the lager space $M \times (-\infty,\infty]$
to obtain the statement for $\Gamma \left(M; \psi \sHom(G,F) \right)$ since $*$-push of $(a,b) \hookrightarrow [a,b]$
sends constant sheaves to constant sheaves.
\end{proof}

Since $\infty$ is a boundary point, we cannot conclude the equivalence using transversality. 
Such situation is considered by Nadler and Shende in \cite{Nadler-Shende} 
where they developed the theory of nearby cycle to study the canonical map 
$$\Gamma \left(M; \psi \sHom(G,F) \right) \rightarrow \Hom(\psi G, \psi F)$$ 
which we recall now.

\begin{definition}[{\cite[Definition 2.2]{Nadler-Shende}}]
A closed subset $X \subseteq S^* M$ is \textit{positively displaceable from legendrians} (pdfl) 
if given any Legendrian submanifold $L$ (compact in a neighborhood of $X$), 
there is a 1-parameter positive family of Legendrians $L_s$, $s\in (-\epsilon,\epsilon)$
(constant outside a compact set), such that $L_s$ is disjoint from $X$ except at $s = 0$.
\end{definition}

\begin{definition}[{\cite[Definition 2.7]{Nadler-Shende}}]
Fix a co-oriented contact manifold $(\mathcal{V},\xi)$ and positive contact isotopy $\eta_s$.
For any subset $Y \subseteq \mathcal{V}$ we write $Y[s] \coloneqq \eta_s(Y)$.
Given $Y, Y^\prime \subseteq \mathcal{V}$ we define the \textit{chord length spectrum}
of the pair to be the set lengths of Reeb trajectories from $Y$ to $Y^\prime$:
$$ cls(Y \rightarrow Y^\prime) = \{ s \in \mathbb{R} | Y[s] \cap Y^\prime \neq \varnothing \}$$
we term $cls(Y) \coloneqq cl(Y \rightarrow Y)$ the chord length spectrum of $Y$.
\end{definition}

\begin{definition}[{\cite[Definition 2.9]{Nadler-Shende}}]
Given a parameterized family of pairs $(Y_b,Y^\prime_b)$ in $S^* M$ over $b \in B$
we say it is gapped if there is some interval $(0,\epsilon)$ uniformly avoided by all $cl(Y_b \rightarrow Y^\prime_b)$. 
In case $Y = Y^\prime$, we simply say $Y$ is gapped.
\end{definition}

\begin{definition}[{\cite[Definition 3.17]{Nadler-Shende}}]
Given a subset $ X \subseteq T^*(M \times J)$, we define its \textit{nearby subset} by
$$\psi (X) \coloneqq \overline{\Pi(X)} \cap T^* (M \times (-1,\infty] ) |_{M \times \{ \infty\}}.$$
\end{definition}

The main theorem for the nearby cycles in \cite{Nadler-Shende} is the following:
\begin{theorem}[{\cite[Theorem 4.2]{Nadler-Shende}}]\label{ncl}
Let $F$, $G$ be sheaves on $M \times J$. Assume

\begin{enumerate}
\item $\ms(F)$ and $\ms(G)$ are $J$-noncharacteristic;
\item $\psi(\ms(F))$ and $\psi(\ms(G))$ are pdfl;
\item The family of pairs in $S^* M$ determined by $(\ms_\pi(F),\ms_\pi(G))$ is gapped for some
fixed contact form on $S^* M$.
\end{enumerate}
Then 
$$ \Gamma(M; \psi \sHom(F,G)) \rightarrow \Hom(\psi(F),\psi(G))$$
is an isomorphism. 
\end{theorem}

Now we apply the theory of nearby cycle to the infinite wrapping functor.

\begin{lemma}
Let $F \rightarrow F^{w_1} \rightarrow F^{w_2} \rightarrow \cdots$ be a sequence in $ \widetilde{\wsh_\Lambda}(M)$ 
as in Lemma \ref{w=ps}.
If for any conic open neighborhood $\mathcal{U}$ of $\Lambda$, 
there exists $n$ such that $\ms(F^{w_n}) \subseteq \mathcal{U}$.
Then the sequence is cofinal.
\end{lemma}

\begin{proof}
Pick a decreasing conic open neighborhood $\mathcal{U}_1 \supseteq \mathcal{U}_2 \supseteq \mathcal{U}_3 
\supseteq \cdots \supseteq \Lambda$ of $\Lambda$ such that $\overline{\mathcal{U}_n} \subseteq \mathcal{U}_{n+1}$
and $\cap_n \, \mathcal{U}_n = \Lambda$.
By taking a subsequence of the $F^{w_n}$'s, we may assume $\ms(F^{w_n}) \subseteq \mathcal{U}_n$ for $n \geq 1$.
Again by taking a subsequence of both the $F^{w_n}$'s and $\mathcal{U}_n$'s, 
we may further assume $\mathcal{U}_{n + 1} \cap \ms(F^{w_n}) =\varnothing$.
By the locality property \ref{loclem}, the continuation map $c: F^{w_n} \rightarrow F^{w_{n+1}}$ depends only
on the value of $\Psi$ on $\mathcal{U}_n$.
Thus, we may modify $\Psi$ on $T^* M \setminus \overline{\mathcal{U}_{n-1}}$ to satisfy the condition in Lemma \ref{cfnc}.
Since taking subsequence won't change the colimit, the original sequence is cofinal.
\end{proof}

\begin{theorem}\label{ffc}
Let $F \in \widetilde{\wsh_\Lambda}(M)$. 
Assume there is a sequence of wrappings $$\Phi_0 \xrightarrow{\Psi_0} \Phi_1 \rightarrow \cdots$$ 
which glues to a (non-compactly supported) 
positive contact isotopy $\Psi: S^* M \times [0,\infty) \rightarrow S^* M$ such that, for any neighborhood $\mathcal{U}$ of $\Lambda$,
$\psi_s(\msif(F)) \subseteq \mathcal{U}$ for $s>>0$.
Then for $G \in \widetilde{\wsh_\Lambda}(M)$ the canonical map,
$$\Hom(G,F^{w_n}) \rightarrow \Hom(G,\psi F^\Psi)$$
is an isomorphism for $n >>0$.
Thus, the canonical map
$$\Hom_w(G,F) \rightarrow \Hom(\wrap_\Lambda^+ G,\wrap_\Lambda^+F)$$
is an isomorphism.
\end{theorem}

\begin{proof}
By Lemma \ref{cfnc}, the sequence $F^{w_n}$ is cofinal and $\wrap_\Lambda^+ F$ 
is computed by $\clmi{n \in \mathbb{N}} \, F^{w_n}$. 
Thus, the first statement implies that $\wrap_\Lambda^+$ induces isomorphisms on the $\Hom$
by Lemma \ref{samehom}.

Now note that $\msif(G)$ in $S^* M$ is compact since $\supp(G)$ is compact 
and the front projection $\pi_\infty:S^* M \rightarrow M$ is proper.
Since a manifold is in particular a regular topological space, there exist open sets $\mathcal{U}$ and $\mathcal{V}$
containing $\Lambda$ and $\msif(G)$ such that $\mathcal{U} \cap \mathcal{V} = \varnothing$.
By restricting to $n >> 0$, we may assume $\psi_s(\msif(F)) \subseteq \mathcal{U}$ is thus disjoint from $\msif(G)$,
which implies that $\sHom(p^* G,F^\Psi)$ is $J$-noncharacteristic. 
Lemma \ref{homtinf} then implies that $\Hom(G,F^{w_n}) = \Gamma( M; \psi \sHom(p^* G, F^\Psi) )$.

So it is sufficient to check the conditions of Theorem \ref{ncl} hold for the pair $p^* G$ and $F^\Psi$.
The set $\ms(p^* G) = \ms(G) \times 0_J$ is tautologically $J$-noncharacteristic.
For $F^\Psi$, we recall that  $\dot{\ms}(F^\Psi) =  \Lambda_\Psi \circ F 
= \{ (\Psi(x,\xi,t) | (x,\xi) \in \dot{\ms}(F), t \in [0,\infty) \}$ which implies that $F^\Psi$ is $J$-noncharacteristic.

By picking a shrinking neighborhood $\mathcal{V}_n$ of $\Lambda$,
we see that the nearby set $\psi(\msif(F^\Psi))$ is contained in $\Lambda$.
One can pick a Whitney triangulation $\mathscr{S}$ such that $\Lambda \subseteq N^*_\infty \mathscr{S}$.
Similarly, up to an isotopy, there exists a Whitney triangulation $\mathscr{T}$ such that $\msif(G) \subseteq N^*_\infty \mathscr{T}$. 
The singular isotropics $N^*_\infty \mathscr{S}$ and $N^*_\infty \mathscr{T}$ are pdfl by Lemma \ref{gpa}.

Now the same argument showing $\sHom(p^* G,F^\Psi)$ is $J$-noncharacteristic implies that
there is an $N \in \mathbb{N}$ such that $\psi_s (\msif(F)) \subseteq \mathcal{V}$ for $s \geq N$.
Thus, when restricting to $M \times [N,\infty)$, $p^* G$ and $F^\Psi$ are microlocally disjoint and 
the gapped condition is tautologically satisfied.
\end{proof}

\begin{corollary}
If there exists a generating set $S$ in $\widetilde{\wsh_\Lambda}(M)$
such that each $F \in S$ admits a sequence $F^{w_n}$ constructed in the above manner
such that $\msif(F^{w_n})$ is contained in arbitrary small neighborhood of $\Lambda$
for $n$ large, then the comparison functor 
$$\wrap_\Lambda^+: \wsh_\Lambda(M) \hookrightarrow \Sh_\Lambda(M)^c$$
is fully faithful.
\end{corollary}

\subsection{Proof of the main theorem} \label{proof_of_the_main_theorem}

We first consider the special case when $\Lambda = N^*_\infty \mathscr{S}$ for some Whitney triangulation $\mathscr{S}$.

\begin{theorem}
The comparison functor $\wrap_{N^*_\infty \mathscr{S} }^+ : \wsh_{N^*_\infty \mathscr{S}}(M) 
\rightarrow \Sh_{N^*_\infty \mathscr{S}}(M)^c$ is an equivalence.
\end{theorem}

\begin{proof}
By Proposition \ref{ballg} and Proposition \ref{b=str}, $\wsh_\mathscr{S}(M)$ has $\{ 1_{\str(s)^-} \}$ as a set of generators.
Recall that $1_{\str(s)^-}$ is defined to be an unspecified inward cornering $\str(s)^{-\epsilon}$ for small enough $\epsilon$.
As mentioned in Definition \ref{incor}, the construction of $\str(s)^{-\epsilon}$ is made so that 
$N^*_{\infty,out} \str(s)^{-\epsilon}$ is disjoint from $N^*_\infty \mathscr{S}$ and 
is contained in arbitrary small neighborhood of $N^*_\infty \mathscr{S}$ as $\epsilon \rightarrow 0$.
Since $\ms( 1_{\str(s)^{-\epsilon} } ) = N^*_{\infty,out} \left( \str(s)^{-\epsilon} \right)$, Theorem \ref{ffc} applies
and $\wrap_{N^*_\infty \mathscr{S} }^+$ is fully faithful.
We see from the same generators $1_{\str(s)^-}$ that $\wrap_{N^*_\infty \mathscr{S} }^+$ is essential surjective
since $\{ 1_{\str(s)} \}$ form a set of generators of $\Sh_{N^*_\infty \mathscr{S}}(M)^c$ by Proposition \ref{mc=cc}.
\end{proof}

To prove the general case, we have to match a special class of objects on both sides.
The following lemma is a sheaf-theoretic variant of the proof \cite[Theorem 5.36]{Ganatra-Pardon-Shende3}.

\begin{lemma}
Let $\Lambda$ be a subanalytic singular isotropic and $(x,\xi) \in \Lambda$ be a smooth point.
For any $F \in \Sh(M)$ such that $\ms(F)$ is contained in $\Lambda$ near $(x,\xi)$, there is an equivalence
$$\Hom(\wrap_\Lambda^+ D_{(x,\xi)},F) = \mu_{(x,\xi)} F.$$
That is, the object $\wrap_\Lambda^+  D_{(x,\xi)}$ co-represents $\mu_{(x,\xi)}$.
\end{lemma}

\begin{proof}
Recall $D_{(x,\xi)}$ is defined to be the cofiber of the canonical map 
$ 1_{f^{-1}(-\infty,- \epsilon)} \rightarrow 1_{f^{-1}(-\infty,\epsilon)}$ 
where $f$ is a proper analytic function defined near $x$ satisfying the following conditions:
There exists an $\epsilon > 0$, 
so that $f$ has only one $\Lambda$-critical point $x$ over $f^{-1}[-\epsilon,\epsilon]$
with $f(x) = 0$, $d f_x = \xi$ and $f^{-1}(-\infty,\epsilon)$ is relatively compact.
Proposition \ref{mstk} implies the function $f$ defines a microstalk functor $\mu_{(x,\xi)}$,
by $$\mu_{(x,\xi)}(F) \coloneqq \Gamma_{ \{f \geq 0\}}(F)_x.$$
As in Remark \ref{ldlp}, the local picture of the fibers $f^{-1}(\{t\})$ for $t \in [-\epsilon,\epsilon]$ is a hyperplane near $x$.
Let $\Psi$ denote any global extension of the wrapping $N^*_{\infty,out} f^{-1}(-\infty,t)$ for $t \in [-\epsilon,\epsilon]$.
If we modify $\Psi$ to $\Psi_0$ by multiplying a bump function supported near $(x,\xi)$ on its Hamiltonian,
the resulting wrapping $(1_{ f^{-1}(-\infty,\epsilon) })^{\Psi_0}$ will appear as expanding
 $f^{-1}(-\infty,\epsilon)$ to some large open set where the expansion happens only near $x$ in the $\xi$ codirection.
The cofiber $\cof\left( c(\Psi_0,1_{ f^{-1}(-\infty,\epsilon) }) \right)$ can be seen as the cofiber induced by
some small open neighborhood $U$ and its open subset $U \cap 1_{ f^{-1}(-\infty,\epsilon)}$.
$$ \includegraphics[scale = 0.2]{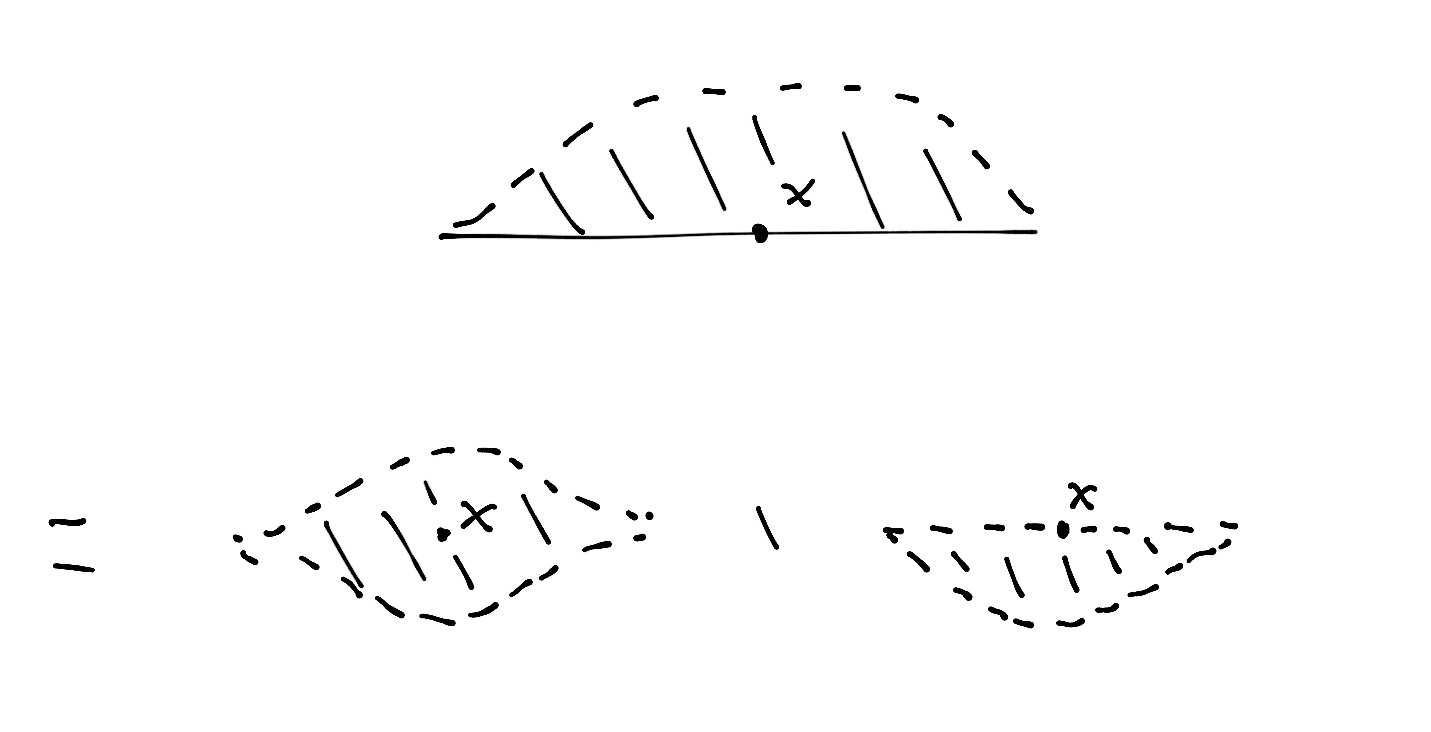}$$
Since the graph $\Gamma_{d f}$ does not intersection $\Lambda$ except at $(x,\xi)$, 
there is an isomoprhism $\Gamma(  f^{-1}(-\infty,\epsilon) ;F) =\Gamma( f^{-1}(-\infty,0) ;F)$ 
by the non-characteristic deformation lemma \ref{ncdef}. 
Thus, $\Gamma_{ \{f \geq 0\}}(F)_x$ can be computed as a colimit 
\begin{align*}
\Gamma_{ \{f \geq 0\}}(F)_x 
&= \clmi{\Psi_0} \Hom(\cof \left( 1_{ f^{-1}(-\infty,0) } 
\rightarrow (1_{ f^{-1}(-\infty,0) })^{w(\Phi_0)} \right) ,F) \\
&= \clmi{\Psi_0} \Hom(\cof \left( 1_{ f^{-1}(-\infty,\epsilon) } 
\rightarrow (1_{ f^{-1}(-\infty,0) })^{w(\Phi_0)} \right) ,F)
\end{align*}
by picking $\Psi_0$ so that the corresponding $U$ as above forms a neighborhood basis of $x$.
Similarly, $\Hom((1_{ f^{-1}(-\infty,0) })^w ,F)$ can be replaced by $\Gamma(  f^{-1}(-\infty,\epsilon) ;F)$
such that the maps are compatible with inclusions of the corresponding open sets.
That is, we are taking colimit over a constant functor and thus
\begin{align*}
\Gamma_{ \{f \geq 0\}}(F)_x 
&= \clmi{\Psi_0} \Hom(\cof \left( 1_{ f^{-1}(-\infty,\epsilon) } 
\rightarrow (1_{ f^{-1}(-\infty,0) })^{w(\Phi_0)} \right) ,F)  \\
&= \clmi{\Psi_0} \Hom( \cof(1_{f^{-1}(-\infty,- \epsilon)} \rightarrow 1_{f^{-1}(-\infty,\epsilon)}), F) \\
&= \Hom( \cof(1_{f^{-1}(-\infty,- \epsilon)} \rightarrow 1_{f^{-1}(-\infty,\epsilon)}), F) = \Hom( D_{(x,\xi)},F)
\end{align*}

Finally, we recall that $\wrap_\Lambda^+$ is defined by the restriction of the left adjoint of 
the tautological inclusion $\Sh_\Lambda(M) \hookrightarrow \Sh(M)$ on $\widetilde{\wsh}_\Lambda(M)$,
and we conclude that $\Gamma_{ \{f \geq 0\}}(F)_x = \Hom(\wrap_\Lambda^+ D_{(x,\xi)},F)$.
\end{proof}

\begin{remark}
Corepresentatives of the microstalk functors $\mu_{(x,\xi)}: \Sh_\Lambda(M) \rightarrow \cV$
are frequently considered since they often provide a preferred set of generators.
For example, Zhou in \cite{Zhou1} finds an explicit description of corepresentatives in the case of FLTZ skeleton first
considered in \cite{Fang-Liu-Treumann-Zaslow}, 
and uses it to match them with certain line bundles on the coherent side, 
which gives an explicit description to the equivalence proved in \cite{Kuwagaki1} through descent argument.
A common recipe for finding such a description is to first find a sheaf $F$ 
which is constructed locally near $x$, and is thus not necessarily in $\Sh_\Lambda(M)$,
but still satisfies the identification $\Hom(F,-) = \mu_{(x,\xi)}$ on $\Sh_\Lambda(M)$. 
Then one constructs a one-parameter family of sheaves $F_t$, $t \in [0,1]$, such that 
$F_0 = F$, $F_1 \in \Sh_\Lambda(M)$, and $\Hom(F_t,-)$ remains constant as $t$ varies.
This lemma can be seen as an abstraction for such a construction when subanalytic structure is presented. 
\end{remark}


\begin{proof}[\hypertarget{proof_main}{Proof} of Theorem \ref{main}]
Pick a Whitney triangulation $\mathscr{S}$ such that $\Lambda \subseteq N^*_\infty \mathscr{S}$.
We use $D^w_{N^*_\infty \mathscr{S},\Lambda}(M)$ denote to the subcategory in 
$\wsh_{N^*_\infty \mathscr{S}}(M)$ generated by the sheaf-theoretical
linking discs $D_{(x,\xi)}$ at Legendrian points of $N^*_\infty \mathscr{S} \setminus \Lambda$ and, similarly,
$D^\mu_{N^*_\infty \mathscr{S},\Lambda}(M)$ the subcategory in $\Sh_{N^*_\infty \mathscr{S}}(M)$ 
generated by the corresponding microstalk representatives.
By Proposition \ref{wdstopr} and Proposition \ref{wgstopr},
they are the fiber of the projections $\wsh_{N^*_\infty \mathscr{S}}(M) \rightarrow \wsh_\Lambda(M)$
and $\Sh_{N^*_\infty \mathscr{S}}(M)^c \rightarrow \Sh_\Lambda(M)^c$ respectively.
Thus, there is a commuting diagram 

$$
\begin{tikzpicture}

\node at (0,2) {$D^w_{N^*_\infty \mathscr{S},\Lambda}(M)$};
\node at (6,2) {$\wsh_{N^*_\infty \mathscr{S}}(M)$};
\node at (12,2) {$\wsh_\Lambda(M)$};
\node at (0,0) {$D^\mu_{N^*_\infty \mathscr{S},\Lambda}(M)$};
\node at (6,0) {$\Sh_{N^*_\infty \mathscr{S}}(M)^c$};
\node at (12,0) {$\Sh_{\Lambda}(M)^c$};

\draw [->, thick] (1.4,2) -- (4.6,2) node [midway, above] {$ $};
\draw [->, thick] (7.5,2) -- (10.8,2) node [midway, above] {$ $};
\draw [->, thick] (1.4,0) -- (4.9,0) node [midway, above] {$ $};
\draw [->, thick] (7.1,0) -- (11.1,0) node [midway, above] {$ $};


\draw [right hook-latex, thick] (0,1.7) -- (0,0.3) node [midway, left] {$ $};
\draw [double equal sign distance, thick] (6,1.7) -- (6,0.3) node [midway, left] {$\wrap_{N^*_\infty \mathscr{S}}^+$};
\draw [->, thick] (12,1.7) -- (12,0.3) node [midway, left] {$\wrap_\Lambda^+$}; 

\end{tikzpicture}.
$$
The last lemma implies that the equivalence
$\wrap_{N^*_\infty \mathscr{S}}^+ : \wsh_\mathscr{S}(T^* M) \xrightarrow{\sim} \Sh_\mathscr{S}(M)^c$ 
restricts to $$\wrap_{N^*_\infty \mathscr{S}}^+: D^w_{N^*_\infty \mathscr{S},\Lambda}(M)
 \xrightarrow{\sim} D^\mu_{N^*_\infty \mathscr{S},\Lambda}(M).$$
Hence, Lemma \ref{prlsn} implies that $\wrap_\Lambda^+$ is an equivalence as well.
\end{proof}

\appendix

\section{Appendix}

We collect some background knowledge from symplectic geometry and higher categorical theory 
used in the main body of the paper in this appendix.

\subsection{Homogenous symplectic geometry and contact geometry}\label{hsnct}

We recall some facts from homogenous symplectic geometry and contact geometry,
and explain how they are interchangeable with each other.
We assume that the contact manifolds in this sections are co-orientable.
Let $(X, d \alpha)$ be a Liouville manifold and let $Z$ denote its Liouville vector field.
We define a homogeneous symplectic manifold to be a Liouville manifold such that the Liouville flow
induces a proper and free $\Rp$-action. 
In this case, the quotient $X/\Rp$ is a manifold.

\begin{definition}
A subset $Y \subseteq (X,d \alpha)$ is \textit{conic} if it is preserved under the $\Rp$-action.
\end{definition}

\begin{proposition}
A coisotropic submanifold $Y$ is conic if and only if $\alpha |_{TY^{d \alpha}} = 0$
where we use $TY^{d \alpha}$ to denote the symplectic orthogonal complement of $TY$.
\end{proposition}

\begin{proof}
Y is coisotropic iff $TY^{d \alpha} \subseteq TY$. 
$Y$ is conic if and only if $Z(y) \in T_y Y$ for $y \in Y$ which implies $\alpha_y(w) = d \alpha_y (Z(y),w) =0$
for all $w \in TY^{d \alpha}$.
Note this direction always holds.
On the other hand, the same equation implies that $Z(y) \in TY$ if $\alpha |_{TY^{d \alpha}} = 0$.
Since $Y$ is coisotropic, $Z(y)$ is in particular in $T_y Y$.
\end{proof}

\begin{corollary}\label{cl}
A Lagrangian submanifold $L \subseteq (X,d \alpha)$ is conic if and only if $\alpha |_L = 0$.
\end{corollary}

\begin{example}
The Liouville vector field $Z$ of the cotangent bundle $T^* M$ can be written locally by
$Z = \sum \xi_i \partial_{\xi_i}$ where the $\xi_i$'s are the dual coordinates of local coordinates $x_i$ of $M$.
The Liouville flow is given by $\Phi^Z_s(x,\xi) = (x, e^s \xi)$
and the $\Rp$-action is simply the multiplication, $r \cdot (x,\xi) = (x, r \xi)$ for  $r \in \Rp$.
\end{example}

\begin{proposition}\label{cq}
The one form $\alpha$ descends to a contact form $\overline{\alpha}$ on $X/\Rp$.
\end{proposition}

\begin{proof}
Let $X$ denote the Liouville vector field associated to $\alpha$ (which is non-vanishing since $\Rp$ acts freely).
By definition, $\alpha(Z) = \omega(Z,Z) = 0$ so it defines a section on  
$(T X/ \langle Z \rangle)^* = T^* (X/\Rp)$.
This is a contact form since on $X$, $\iota_Z \omega \wedge d (\iota_Z \omega)^{n-1} 
= \iota_Z \omega \wedge (\Lie{Z} \omega)^{n-1} = \iota_Z \omega \wedge \omega^{n-1} 
= \frac{1}{n} \iota_Z \omega^{n-1}$ and $T (X/\Rp)$ can be identified as vectors transversal to $Z$.
\end{proof}

\begin{example}
The example we will study in this paper is the cotangent bundle away from the zero section $\dT^* M$ 
for some smooth manifold $M$.
Pick a metric $g$ and restrict the projection $p: \dT^* M \rightarrow S^* M$ to $\{ (x,\xi) | g_x(\xi,\xi) =1 \}$
and denote it as $p_g$. The map $p_g$ is a diffeomorphism because its domain is transversal to 
the $\Rp$-action and $p_g$ is clearly one-to-one. 
Its inverse $s:S^* M \rightarrow \dT^* M$ provides $S^* M$ a global contact form $s^* \alpha_{can}$.

Note any such section gives the same contact structure but there might not be any contactomorphism sending 
one contact form to another. A more intrinsic description of the contact structure is $\eta_{[x,\xi]} = \ker \xi$.
\end{example}

\begin{lemma}
A homogeneous symplectomorphism $\psi: (X, d \alpha) \rightarrow (Q, d \beta)$ preserves 
the Liouville form, i.e., $\psi^*  \beta = \alpha$.
\end{lemma}

\begin{proof}
let $\psi$ be homogeneous and $\psi^* d \beta = d \alpha$.
We denote the Liouville vector fields by $Z$ and $Y$ and the corresponding flow by $\phi^Z_t$ and $\phi^Y_t$, $t \in \mathbb{R}$. 
Since $\psi$ is a homogeneous symplectomorphism, we have $\psi (\phi^Z_t(x)) = \phi^Y_t( \psi(x))$ for all $x \in X$.
Differentiate the equation and evaluate at $0$, we obtain that $d \psi_{x} (Z(x) ) = Y( \psi(x))$, i.e., $Y = \psi_* Z$.
So for any differential $p$-form $\nu$ on $Q$, 
\begin{align*}
\left( \psi^* (\iota_Y \nu) \right)(v_1,\cdots,v_{p-1})
&=(\iota_Y \nu) (\psi_* v_1, \cdots, \psi_* v_{p-1}) \\
&= \nu(Y, \psi_* v_1, \cdots, \psi_* v_{p-1}) \\
&= (\psi^* \nu)(Z,v_1,\cdots, v_{p-1}) \\
&= \left( \iota_Z(\psi^* \nu)\right) (v_1,\cdots,v_{p-1}).
\end{align*}
That is, $\psi^* \circ \iota_Y = \iota_Z \circ \psi^*$.
In particular, $\psi^* \alpha = \psi^* \iota_Z d \alpha = \iota_Y \psi^* d \alpha = \beta$.
\end{proof}

\begin{proposition}
A co-orientation preserving contactomorphism $\varphi: (N,\xi) \rightarrow (P,\eta)$ 
gives rise to a unique homogeneous symplectomorphism $\tilde{\varphi}: SN \rightarrow SP$ between their symplectizations.
On the other hand, a homogeneous symplectomorphism $\psi: (X, d \alpha) \rightarrow (Q, d \beta)$
induces a contactomorphism on the contact quotient in Proposition \ref{cq}.
These two constructions are inverse to each other if $X$ and $Q$ come from symplectization.
\end{proposition}

\begin{proof}
Assume $(N,\xi)$ and $(P,\eta)$ are co-oriented by $\alpha$ and $\beta$.
The equation $d \phi_x (\xi_x) = \eta_{\varphi(x)}$ implies that $\varphi^* \beta = h \alpha$ 
for some $h>0$. 
More precisely, let $R$ be the Reeb vector field of $\alpha$, then $h = \beta (\varphi_* R)$.
Define $\tilde{\varphi}: (N \times \Rp, d(t \alpha)) \rightarrow (P \times \Rp, d(s \beta))$ by 
$\tilde{\varphi}(x,t) = ( \varphi(x), (h(x))^{-1} t) $.
Then $\tilde{\varphi}^* s \beta = t \frac{1}{h} \varphi^* \beta = t \alpha$ 
so $\tilde{\varphi}$ is a homogeneous symplectomorphism.
Now assume there is another lifting $\tilde{\varphi}^\prime$.
Since they both descend to $\varphi$, there is $g > 0$ such that $\tilde{\varphi}^\prime(x,t) = g(x) \tilde{\varphi}(x,t)$.
But then $t \alpha = (\tilde{\varphi}^\prime)^* t \beta = g \varphi^* t \beta = g t \alpha$ so $g \equiv 1$.
Since $\psi$ preserves the Liouville form, it is clear that $\psi$ descends to a contactomorphism on the quotient.
And we also see the that two constructions are inverse to each other when the homogeneous symplectic manifolds
are given by symplectization.
\end{proof}

\begin{example}
We consider the case when $X = Q = \dT^* M$ is the cotangent bundle away from the zero section
for some manifold $M$.
One can identify it as the symplectization of $S^* M$ by picking a metric $g$.
Let $\varphi: S^* M \rightarrow S^* M$ be a co-orientation preserving 
contactomorphism and we would like to lift it to a homogeneous
symplectomorphism $\hat{\varphi}: \dT^* M \rightarrow \dT^* M$.

We describe here how the identification intertwines with the construction in the proposition.
Denote $s$ the section of $p: \dT^* M \rightarrow S^* M$ which is given by the unit covectors. 
We claim that there is a (unique) section $t: S^* M \rightarrow \dT^* M$ so that 
$\varphi^* (t^* \alpha_{can}) = s^* \alpha_{can}$.
(Note we cannot just require $t = s \circ \varphi^{-1}$ since this would implies $\id_{S^* M} = \varphi^{-1}$.)
If such $t$ exists, then $t^* \alpha_{can} = (\varphi^{-1})^* s^* \alpha_{can} = h s^* \alpha_{can}$
for some $h \in C^\infty(S^* M;\Rp)$ given by $\varphi$.
So we simply define $t:S^* M \rightarrow \dT^* M$ by $t = h \cdot s$ where $`\cdot'$ is the $\Rp$ action.
Then we can define $\hat{\varphi}:\dT^* M \rightarrow \dT^* M$ by 
$\hat{\varphi} = \sqrt{g} \cdot  (t \circ \phi \circ p)$.
One can compute that 
\begin{align*}
\hat{\varphi}^* \alpha_{can} 
&= \sqrt{g} \cdot ( p^* \circ \varphi^* \circ t^* \alpha_{can}) 
= \sqrt{g} \cdot (p^*  \circ s^* \alpha_{can}) \\
&= \sqrt{g} \cdot (s \circ p)^* \alpha_{can} 
= \sqrt{g} \cdot \frac{1}{\sqrt{g}} \alpha_{can}
= \alpha_{can}
\end{align*}
is a symplectomorphism.
Note that we use $s \circ p (x,\xi) = \left(1/\sqrt{g_x(\xi,\xi)} \right) (x,\xi)$ for the second to last equality.

Now consider a family of isotopy $\varphi_t :S^* M \rightarrow S^* M$ such that  $\varphi_0 = \id_{S^* M}$.
The requirement $(\varphi_t^{-1})^* s^* \alpha_{can} = h_t s^* \alpha_{can}$ ensures $h_t > 0$ since $h_0 \equiv 1$. 
We can then lift $\varphi_t$ to a family of homogeneous symplectomorphism 
$\hat{\varphi}_t: \dT^* M \rightarrow \dT^* M$ by the above process.
Since this process can be reserved, we see that there is a one-to-one correspondence between
contact isotropies on $S^* M$ and homogeneous isotropies on $\dT^* M$.
Note that the family version of isotopies works similarly.
\end{example}

Finally, we mention a lemma which we will use in the main text.
\begin{lemma}\label{homogeneous_hamiltonian_extends}
Let $H: \dT^* M \rightarrow \RR$ be a homogeneous Hamiltonian of degree $1$.
Then $H^2$ extends smooth to $T^* M$ by setting $H^2(x,0) = 0$. 
\end{lemma}

\begin{proof}
Being homogeneous of degree $1$ in this case means $H(x,r \xi) = r H(x,\xi)$ for $r \in \RR$.
Thus the quotient $H(x,\xi) / \abs{\xi}$ is homogeneous of degree $0$, i.e.,
$H(x,\xi) = f(x,[\xi]) \abs{\xi}$ for some smooth function $f \in C^\infty(S^* M)$ on the cosphere bundle.
The function $H^2(x,\xi)$ is then the product $f(x,[\xi])^2 \abs{\xi}^2$, which extends smoothly to $T^* M$.
\end{proof}


\subsection{Quotients of small stable categories}
We discuss quotients of (small) stable categories.
Let $\sC_0 \in \st$ be a small idempotent complete stable category and $K$ be a collection of objects in $\sC_0$ .
We would like to construct an associated localization $\sC_0 \rightarrow  \sC_0 / K$ so that
a morphism $f: X \rightarrow Y$ with $\cof(f) \in K$ becomes an isomorphism in $\sC_0 / K$.
This localization can be defined as a quotient in the following way. 
First, we take the stable subcategory $\langle K \rangle$ generated by $K$, and then take its idempotent completion
which we will denote it by $\sN = \sN(K)$.
Abstract argument implies $\sN$ is still stable and is embedded in $\sC_0$ as the subcategory
of retracts of objects in $\langle K \rangle$ because $\sC_0$ is idempotent complete.
Let $\iota_*: \sN \hookrightarrow \sC_0$ denote the inclusion and we recall that
colimits exist in $\st$ by Proposition \ref{clmbrad}.

\begin{definition}
We define the quotient $\sC_0/K$ of $\sC_0$ by $K$ as the cofiber $\cof{(\iota_*)}$ taken in $\st$
and use $j^*: \sC_0 \twoheadrightarrow \sC_0/K$ to denote the projection.
\end{definition}

We have the following description of morphisms in $\sC_0/K$.

\begin{proposition}\label{hombc}
Let $X$, $Y$ be objects in $\sC_0$. Then the $\Hom$ in $\sC_0/K$ can be computed as a colimit,
$$ \Hom_{\sC_0/K} (j^* X,j^* Y) = \clmi{Y \rightarrow Y^\prime} \Hom_{\sC_0} (X,Y^\prime)$$
where the colimit runs through the morphism $Y \rightarrow Y^\prime$ whose cofiber is in $\sN$.
Alternatively, we can compute the Hom-spaces by varying the first component, i.e.,
$$ \Hom_{\sC_0/K} (j^* X,j^* Y) = \clmi{X^\prime \rightarrow X} \Hom_{\sC_0} (X^\prime,Y)$$
with $\cof(X^\prime \rightarrow X) \in \sN$.
\end{proposition}

To prove the proposition, we first look more closely into the construction.
Begin with the inclusion $\sN \xhookrightarrow{\iota_*} \sC_0$, 
we translate to the category $\PrLcg$ by taking $\Ind$ by Proposition \ref{cs=st} and obtain
$$ \Ind (\sN) \xhookrightarrow{\Ind (\iota_*)} \Ind (\sC_0).$$
Because $\Ind (\iota_*)$ preserves small colimits, it admits a right adjoint $\Ind (\iota_*)^R$.

\begin{lemma}
For $X \in \sC_0 \hookrightarrow \Ind (\sC_0)$, the right adjoint of $\Ind (\iota_*)$ can be given by the formula
$$\Ind (\iota_*)^R (X) = ``\clmi{\alpha: Z \rightarrow X, \ Z \in \sN}" Z.$$
Here we use the quotation $``\colim"$ to emphasis the colimit is taken formally in $\Ind (\sN)$.
\end{lemma}

\begin{proof}
By definition, the formal colimit $``\clmi{\alpha: Z \rightarrow X, \ Z \in \sN}" Z$ is an object of $\Ind(\sN)$.
Because an object of $\Ind(\sN)$ is of the form $``\varinjlim" W$ over some filtered colimit by some objects in $\sN$,
it is sufficient to show, for all $W \in \sN$,
$$\Hom_{\sC_0}(\iota W, X) = \Hom_{\Ind(\sN)}( W, ``\clmi{\alpha: Z \rightarrow X, \ Z \in \sN}" Z).$$
Since $\Ind (\iota)$ preserves compact objects, we compute
$$\Hom_{\Ind(\sN)}( W, ``\clmi{\alpha: Z \rightarrow X, \ Z \in \sN}" Z)
= \clmi{\alpha: Z \rightarrow X, \ Z \in \sN} \Hom_\sN (W,Z)
= \Hom_{\sC_0}(W,X).$$
\end{proof}

\begin{proof}[Proof of the proposition \ref{hombc}]
By passing to right adjoints using Proposition \ref{cs=st}, the cofiber $\cof(\Ind(\iota_*))$ can be computed as, 
$$j_*: \fib( \Ind(\iota_*)^R) \subseteq \Ind (\sC_0),$$
the subcategory consisting of objects $Y$ such that $\Hom(X,Y) = 0$ for all  $X \in \sN$.
Because we are in the stable setting, this gives us a fiber sequence
$$ \Ind(\iota_*) \Ind(\iota_*)^R \rightarrow \id \rightarrow j_* j^*.$$
Thus, for $X, Y \in \sC_0$, one computes
\begin{align*}
\Hom_{\sC_0/K}(j^* X, j^* Y)
&= \Hom_{\sC_0}(X, j_* j^* Y) \\
&= \Hom_{\sC_0} \left(X, \cof( \Ind(\iota_*) \Ind(\iota_*)^R Y \rightarrow Y) \right) \\
&= \Hom_{\sC_0} \left(X, \cof( ``\clmi{\alpha: Z \rightarrow Y, \ Z \in \sN}"  Z \rightarrow Y) \right) \\
&= \clmi{\alpha: Z \rightarrow Y, \ Z \in \sN}  \Hom_{\sC_0} \left(X, \cof(Z \rightarrow Y) \right) \\
&= \clmi{Y \xrightarrow{\beta} Y^\prime, \ \cof(\beta) \in \sN}  \Hom_{\sC_0} (X, Y^\prime ). 
\end{align*}
Here, we notice the last equation is simply a change the expression for the same colimit.

To obtain the similar formula in which we vary the first component, we notice that there is an equivalence
$(\sC_0/K )^{op} = \sC_0^{op}/\sN^{op}$
because they satisfy the same universal property.
We thus compute
\begin{align*}
\Hom_{\sC_0/K}(j^* X, j^* Y)
&= \Hom_{ (\sC_0/K)^{op} }(j^* Y, j^* X) \\
&= \Hom_{\sC_0^{op} / \sN^{op} }(j^* Y, j^* X) \\
&= \clmi{X \xleftarrow{\gamma}X^\prime, \ \cof(\gamma) \in \sN^{op}}  \Hom_{ \sC_0^{op} } (Y, X^\prime) \\
&= \clmi{X^\prime \xrightarrow{\gamma} X, \ \cof(\gamma) \in \sN} \Hom_{ \sC_0 } (X^\prime, Y).
\end{align*}
\end{proof}

We will use the following "snake lemma" for categories.

\begin{lemma}\label{prlsn}
Consider the following diagram in $\PrLcs$ :

$$
\begin{tikzpicture}

\node at (0,2) {$\sC_0$};
\node at (3,2) {$\sC$};
\node at (6,2) {$\sC_1$};
\node at (0,0) {$\sD_0$};
\node at (3,0) {$\sD$};
\node at (6,0) {$\sD_1$};

\draw [right hook-latex, thick] (0.5,2) -- (2.5,2) node [midway, above] {$i$};
\draw [->>, thick] (3.5,2) -- (5.5,2) node [midway, above] {$p$};
\draw [right hook-latex, thick] (0.5,0) -- (2.5,0) node [midway, above] {$j$};
\draw [->>, thick] (3.5,0) -- (5.5,0) node [midway, above]  {$q$};

\draw [->, thick] (0,1.7) -- (0,0.3) node [midway,right] {$F_0$};
\draw [->, thick] (3,1.7) -- (3,0.3) node [midway, right] {$F$};
\draw [->, thick] (6,1.7) -- (6,0.3) node [midway,right] {$\bar{F}$};

\end{tikzpicture}
$$
where $p$ and $q$ are the quotient functor of the inclusion $i$ and $j$,
$F_0$ is the restriction of $F$ which factors through $\sD_0$ and $\bar{F}$ is the induced functor between the quotients.
Let $\iota: \fib(\bar{F}) \hookrightarrow \sC_1$ denote the fiber of $\bar{F}$,
$\pi: \sD_0 \twoheadrightarrow \cof(F_0)$ the cofiber of $F_0$
and $\partial: \fib(\bar(F)) \rightarrow \cof(F_0)$ the functor given by the composition 
$\partial = \pi \circ j^R \circ F \circ p^R \circ \iota$.
If $F$ is an equivalence, then $\partial$ is an equivalence.
\end{lemma}

\begin{proof}
For simplicity, we assume $\sC = \sD$ and $F$ is the identity so the diagram becomes,

$$
\begin{tikzpicture}

\node at (6,4) {$\fib( \bar{F})$};
\node at (0,2) {$\sC_0$};
\node at (3,2) {$\sC$};
\node at (6,2) {$\sC_1$};
\node at (0,0) {$\sD_0$};
\node at (3,0) {$\sC$};
\node at (6,0) {$\sD_1$};
\node at (0,-2) {$\cof(F_0)$};

\draw [right hook-latex, thick] (0.5,2) -- (2.5,2) node [midway, above] {$i$};
\draw [->>, thick] (3.5,2) -- (5.5,2) node [midway, above] {$p$};
\draw [right hook-latex, thick] (0.5,0) -- (2.5,0) node [midway, above] {$j$};
\draw [->>, thick] (3.5,0) -- (5.5,0) node [midway, above]  {$q$};

\draw [left hook-latex, thick] (5.7,1.6) to [out=210,in=-30] (3.3,1.6);
\node at (4.5,1.6) {$p^R$};
\draw [->>, thick] (2.7,-0.4) to [out=210,in=-30] (0.3,-0.4);
\node at (1.5,-0.4) {$j^R$};

\draw [right hook-latex, thick] (6,3.7) -- (6,2.3) node [midway,right] {$\iota$};

\draw [right hook-latex, thick] (0,1.7) -- (0,0.3) node [midway,right] {$F_0$};
\draw [double equal sign distance, thick] (3,1.7) -- (3,0.3) node [midway, right] {$ $};
\draw [->>, thick] (6,1.7) -- (6,0.3) node [midway,right] {$\bar{F}$};

\draw [->>, thick] (0,-0.3) -- (0,-1.7) node [midway,right] {$\pi$};

\end{tikzpicture}
$$

We will prove that the functor $\theta \coloneqq \iota^R \circ p \circ j \circ \pi^R$ is the inverse 
by showing that $\theta \circ \partial = \id_{ \fib(\bar{F}) }$.
The equation $\partial \circ \theta = \id_{ \cof(F_0) }$ can be proved similarly.
First write out $\theta \partial$ as $\iota^R p j \pi^R \pi j^R p^R \iota$.
Recall that in the stable setting, the sequence $\sC_0 \hookrightarrow \sC \twoheadrightarrow \sC_1$
comes with a fiber sequence of functors $$i i^R \rightarrow \id_{\sC} \rightarrow p^R p.$$
Apply this fact to $\sC_0 \hookrightarrow \sD_0 \twoheadrightarrow \cof(F_0)$,
we see there is a fiber sequence $$F_0 F_0^R \rightarrow  \id \rightarrow \pi^R \pi.$$
Apply $j \circ (-) \circ j^R$ and the fiber sequence becomes $i i^R \rightarrow  j j^R \rightarrow j \pi^R \pi j^R.$
Now apply $p \circ (- ) \circ p^R$ and the we see that $p j j^R p = p j \pi^R \pi j^R p^R$ since $p \circ i =0$.
Thus, we can simplify $\theta \partial$ to  $\iota^R p j j^R p^R \iota$.
Similar argument allows us to further simplify $\theta \partial$ to $\iota^R p p^R \iota 
= \iota^R \id_\sC \iota = \id_{\fib(\bar{F})}$.
\end{proof}

\bibliographystyle{plain}
\bibliography{MicrolocalSheaf_REF}

\end{document}